%% file: calcul-moulien.tex
\documentclass[a4paper,11pt]{smfart}

\usepackage[french]{babel}
\usepackage{amssymb}
\usepackage{latexsym}

\def\la{\langle}
\def\ra{\rangle}
\def\ni{\noindent}
\def\di{\displaystyle}

\def\un{\underline}

\def\bs{\underline{s}}
\def\bx{\underline{x}}

\def\bn{\underline{n}}
\def\ba{\underline{a}}

\def\blambda{\underline{\lambda}}
\def\bomega{\underline{\omega}}

\def\kK{\mathbb{K}}
\def\aA{\mathbb{A}}
\def\eE{\mathbb{E}}
\def\rR{\mathbb{R}}
\def\nN{\mathbb{N}}

\def\qQ{\mathbb{Q}}
\def\cC{\mathbb{C}}

\def\sh{\mbox{\rm sh}}
\def\csh{\mbox{\rm csh}}

\def\shi{\mbox{\rm shi}}

\def\td{\bigtriangledown}
\def\cal{\mathcal}
\def\Exp{\mbox{\rm Exp}}

\numberwithin{equation}{part}
\newtheorem{exemp}{Exemple}[part]
\newtheorem{thm}{Th\'eor\`eme}[part]
\newtheorem{defi}{D\'efinition}[part]
\newtheorem{lem}{Lemme}[part]
\newtheorem{pro}{Proposition}[part]
\newtheorem{cor}{Corollaire}[part]
\newtheorem{rema}{Remarque}[part]

\title{Calcul Moulien}
\alttitle{Mould calculus}
\author{Jacky Cresson}
\email{cresson@math.univ-fcomte.fr}

\address{Equipe de Math\'ematiques de Besan\c{c}on, CNRS-UMR
6623,\\ Universit\'e de Franche-Comt\'e, 16 route de Gray, 25030 Besan\c{c}on
cedex, France.}

\begin{document}
\baselineskip 6mm
\setcounter{tocdepth}{2}

\begin{abstract}
Ce texte est une introduction au calcul moulien, d\'evelopp\'e par Jean \'Ecalle. On donne une
d\'efinition pr\'ecise de la notion de moule ainsi que les principales propri\'et\'es de ces objets.
On interpr\`ete les diff\'erentes sym\'etries (alterna(e)l,symetra(e)l) des moules via les
s\'eries formelles non commutatives asso\c{c}i\'ees dans des big\`ebres gradu\'ees not\'ees
$\aA$ et $\eE$, correspodant aux deux types de colois \'etudi\'ees par Ecalle, \`a savoir
$\Delta (a)=a\otimes 1+1\otimes a$ et $\Delta_* (a_i)=\di\sum_{l+k=i} a_l \otimes a_k$.
On illustre en d\'etail l'application de ce formalisme dans le domaine de la recherche des formes normales de
champs de vecteurs et diff\'eomorphismes.
\end{abstract}

\begin{altabstract}
This paper is an introduction to mould calculus, as introduced by Jean \'Ecalle. We give a precise
definition of moulds and describe there main properties. We translate mould symmetries (alterna(e)l and
symetra(e)l) using non
commutative formal power series in two given bialgebras $\aA$ and $\eE$, corresponding to two coproducts structure
given by $\Delta (a)=a\otimes 1+1\otimes a$ and $\Delta_* (a_i)=\di\sum_{l+k=i} a_l \otimes a_k$.
We apply this formalism to the problem of normal forms for vector fields and diffeomorphisms.
\end{altabstract}


\maketitle

\newpage

\tableofcontents

\newpage
\include{avant-propos}

\newpage
\part*{Introduction}

\section*{Les raisons du texte}

Toute nouvelle notion, en math\'ematiques comme ailleurs, est
souvent source de difficult\'es. C'est notamment le cas des {\it
moules} et {\it comoules} introduits par Jean Ecalle au cours de
ses nombreux travaux, aussi bien dans le domaine des formes
normales de champs de vecteurs que dans celui de l'\'etude des
singularit\'es, ou plus recemment, des propri\'et\'es
alg\'ebriques des polyz\^etas. Par ailleurs, il n'est pas toujours
facile de situer ces nouvelles notions et ses \'eventuelles
connexions vis \`a vis des doma\^\i nes math\'ematiques
existant.\\

Il n'existe pas de texte introductif sur le calcul moulien, facile
d'acc\`es, et r\'epondant \`a ces interrogations. Or, ce travail
devient n\'ecessaire du fait de l'utilisation de ce formalismes
dans des domaines tr\`es diff\'erents de son champ d'application
initial\footnote{Je pense ici aux r\'ecentes contributions de Jean
Ecalle en th\'eorie des nombres sur la combinatoire des
polyz\^etas.}. L'objet de ce texte\footnote{Ce texte est bas\'e
sur mon cours de DEA ``Calcul moulien et s\'eries formelles non
commutatives" donn\'e \`a l'Universit\'e Semlalia de Marrakech en
avril 2002, et une s\'erie d'expos\'es intitul\'es ``Calcul
moulien et polyz\^etas" donn\'es en Mai 2002 \`a Paris VI en
compl\'ement du cours de DEA de Michel Waldschmidt ``Valeurs
z\^etas multiples".} est de combler ce vide. Nous allons d\'efinir
les moules, et les comoules, et expliciter certaines m\'ethodes
associ\'ees, comme la m\'ethode d'arborification. Mais surtout,
nous allons pr\'eciser le cadre th\'eorique de ces objets, \`a
savoir celle de la combinatoire des alg\`ebres de Lie libre et des
big\`ebres gradu\'ees. On verra notamment que certains des mots
nouveaux sont en fait connus depuis longtemps sous une autre
terminologie, comme par exemple, l'\'equivalence entre un moule
alternal et un \'el\'ement primitif
d'une alg\`ebre de Hopf.\\

La plupart du temps, le {\it calcul moulien}\footnote{Quand il n'est pas vu comme une ``th\'eorie" des moules.} est per\c{c}u
comme compliqu\'e et difficile \`a manier. Cette remarque permet aux sp\'ecialistes des divers champs d'applications du
calcul moulien de le laisser de c\^ot\'e. Une autre cons\'equence est que les apports de ce calcul, par exemple
dans l'\'etude des champs de vecteurs et des diff\'eomorphismes, sont souvent mal cern\'es voir ignor\'es.\\

Le but de ce texte est de montrer la simplicit\'e d'usage du
calcul moulien, et son apport aussi bien th\'eorique que pratique,
dans le probl\`eme classique des formes normales de champs de
vecteurs, via la d\'emonstration de plusieurs th\'eor\`emes
connus, notamment celui de Bryuno sur la lin\'earisation
analytique des champs de vecteurs non-r\'esonnants en pr\'esence
de petits diviseurs.

Bien entendu, les moules permettent {\it aussi} de d\'emontrer des r\'esultats
nouveaux, souvent difficiles d'atteinte pour les m\'ethodes
classiques\footnote{Par exemple, la d\'emonstration de l'analyticit\'e de la
correction dans \cite{EV} pour un champ de vecteur quelconque, qui correspond
dans le cas hamiltonien \`a la conjecture de Gallavotti.}. Mais il est plus
facile de s'appercevoir de la puissance de ce langage en comparant le degr\'e
de compr\'ehesion obtenu en suivant le chemin classique de d\'emonstration d'un
th\'eor\`eme bien \'etabli, et l'approche moulienne. Outre un apport conceptuel
\'evident,  il am\'eliore la compr\'ehension de r\'esultats existants, au point
parfois, de remettre en cause des ph\'enom\`enes pourtant bien
\'etablis\footnote{C'est le cas de l'\'etude de la convergence de la {\it
correction}, introduite par J. Ecalle et B. Vallet \cite{EV}, qui met en
\'evidence un ph\'enom\`ene {\it purement alg\'ebrique} de suppression des {\it
petits diviseurs surmultiples}. Or, ces derniers appara\^\i ssent dans les
manipulations classiques \cite{el}, et la convergence de la s\'erie se
d\'emontre alors au prix d'estimations tr\`es fines, connues sous le nom de
{\it
compensation d'Eliasson}.}\\

La puissance des moules s'illustre aussi par son vaste champ d'applications, notamment :\\

i - R\'esurgence (\'equation du pont),

ii - \'Equations diff\'erentielles (th\'eorie des invariants
holomorphes, th\'eorie des formes normales),

iii - \'Equations aux diff\'erences,

iv - \'Equations fonctionelles,

v - Analyse de Lie,

vi - Th\'eorie des nombres (multizetas symboliques).\\

J'en oublie surement et je renvoie aux travaux d'Ecalle, notamment son article de revue \cite{E6} pour plus de d\'etails.
Passons maintenant \`a la pr\'esentation standard des moules.

\section*{Pr\'esentation standard des moules}

La pr\'esentation des moules qui va suivre est essentiellement celle donn\'ee par Ecalle dans ses articles comme ``rappels"
sur les moules.\\

Soit $\Omega$ un semigroupe additif. On note $\Omega^*$
l'ensemble des
suites construites sur $\Omega$. Un \'el\'ements de $\Omega^*$
sera
not\'e $\un{\omega }= (\omega_1 ,\dots ,\omega_n )$ avec $\omega_i \in
\Omega$ pour tout $i$.
La longueur d'une suite $\un{\omega }$ se note $l(\un{\omega } )=n$.
On indicera les
\'el\'ements de $\Omega^*$ sous la forme $\un{\omega}^i$, et ses
composantes comme
$\un{\omega }^i =(\omega_1^i ,\dots ,\omega_n^i )$. On note
$\un{\omega}^1 \bullet \un{\omega}^2$ la
suite obtenue par {\it concat\'enation} des suites $\un{\omega}^1$ et
$\un{\omega}^2$ (ou $\un{\omega}^1 \un{\omega}^2$ si aucune confusion
n'est possible).\\

La plupart des textes sur le calcul moulien donne la ``d\'efinition" suivante: {\it les moules sont des fonctions \`a un nombre
variable de variables}. L'avantage de cette pseudo d\'efinition est qu'elle est ``parlante", et que l'on devine assez bien ce
qui se cache sous ce type d'objet. Dans la suite, nous adopterons la d\'efinition suivante:

\begin{defi}
Soit $\kK$ un anneau commutatif. On appelle moule une application de
$\Omega^*$ dans $\kK$, i.e.
$$
\left .
\begin{array}{lll}
\di\Omega^* & \di\stackrel{M}{\rightarrow} & \kK \\
\un{\omega } & \mapsto & M (\un{\omega } ) .
\end{array}
\right .
$$
\end{defi}

On notera un moule $M^{\bullet }$ et son \'evalutation sur une suite $\un{\omega}$ par $M^{\un{\omega}}$.\\

Cette notation pose de nombreux probl\`emes et peut pr\'eter \`a confusion\footnote{On peut n\'eanmoins lui trouver une
justification par la suite dans la relative simplicit\'e qu'elle induit sur les manipulations alg\'ebriques effectives des moules.}. N\'eanmoins, afin de faciliter la lecture des
textes d'Ecalle par la suite et de permettre une comparaison, nous allons conserver cette notation.\\

Les diff\'erentes propri\'et\'es des moules (alterna(e)lit\'e, sym\'etra(e)lit\'e) proviennent de la contraction
avec des {\it op\'erateurs non commutatifs} v\'erifiants certaines {\it colois sp\'ecifiques}. Ces
op\'erateurs interviennent naturellement dans l'\'etude des champs de vecteurs et
diff\'eo\-mor\-phismes locaux de $\cC^{\nu }$.

\begin{defi}
Soit $A$ une alg\`ebre sur un corps de caract\'eristique z\'ero $\kK$. Soit
$\cal B$ une alg\`ebre
d'op\'erateurs sur $A$, non commutatifs, munie de la composition
(usuelle). On
appelle comoule une application de $\Omega^*$ dans $\cal B$. On
note
\begin{equation}
\left .
\begin{array}{lll}
\di\Omega^* & \di\stackrel{B}{\rightarrow} & {\cal B} \\
\un{\omega } & \mapsto & B_{\un{\omega }} =B_{\omega_n } \dots
B_{\omega_1 } .
\end{array}
\right .
\end{equation}
\end{defi}

On notera un comoule $B_{\bullet }$.\\

Ces alg\`ebres seront toujours munies d'une coloi $\Delta : {\cal B}
\rightarrow {\cal B}\otimes_k
{\cal B}$, application lin\'eaire, compatible avec la loi de composition
sur $\cal B$. Dans ce
cas, $({\cal B} ,\Delta )$ forme une {\it big\`ebre} au sens de Bourbaki
(\cite{Bo1}). \'Ecalle
consid\`ere deux types de colois qui sont :\\

i - $\Delta \ :\  B_{\omega } \rightarrow B_{\omega } \otimes 1
+1\otimes B_{\omega } $ ,

ii - $\Delta_* \ :\  B_{\omega }  \rightarrow \di\sum_{\omega_1
+\omega_2 =\omega } B_{\omega_1 }
\otimes B_{\omega_2 } $\footnote{Dans ce cas, l'alphabet $\Omega$ doit poss\'eder une structure de semi-groupe pour donner
un sens \`a $\omega_1 +\omega_2$.}.\\

Le cas i) correspond aux alg\`ebres de d\'erivations sur $A$ et ii)
g\'en\'eralise
l'alg\`ebre des op\'erateurs diff\'erentiels d'ordre $p\geq 1$.\\

L'objet ``contract\'e" d'un moule et d'un comoule se note $\di
P_{M}=\di\sum_{\bullet } M^{\bullet } B_{\bullet }$, \'etant entendu
que la somme est faite sur toutes les suites possibles de
$\Omega^*$.\\

Deux questions se posent naturellement\footnote{Le naturel dont il est question ici peut prendre deux formes,
suivant l'int\'er\^et du lecteur:

- soit au niveau des applications, en montrant que la recherche des \'el\'ements primitifs et group-like est indispensable.
C'est ce qui est fait dans le dernier chapitre sur les formes normales de champs de vecteurs et diff\'eomorphismes.

- soit au niveau alg\'ebrique, ou ces deux notions correspondent aux d\'erivations et automorphismes de
l'alg\`ebre sous-jacente, ce qui est d\'etaill\'e dans la partie 2.} dans ce cadre: \\

i) $P_M$ est-il un \'el\'ement {\it primitif}, i.e. $\Delta (P_M )=P_M \otimes 1 +1\otimes P_M$ (m\^eme chose avec
$\Delta_*$) ? \\

ii) $P_M$ est-il un \'el\'ement {\it group-like}\footnote{On trouve aussi la terminologie moins courante
d'\'el\'ement {\it diagonal}.}, i.e. $\Delta (P_M )=P_M \otimes P_M$ (m\^eme chose avec
$\Delta_*$) ? \\

Cette information est toute enti\`ere contenue dans les propri\'et\'es
alg\'e\-bri\-ques du moule
$M^{\bullet }$. Autrement dit, les propri\'et\'es du moule dictent la
nature de l'objet formel
associ\'e\footnote{La terminologie de moule correspond bien \`a cette id\'ee.}.
Dans le cas ou les comoules appartiennent \`a une alg\`ebre
de d\'erivation, on obtient les \'el\'ements primitifs si le moule est
{\it alternal} et les \'el\'ements
``group-like" si le moule est {\it sym\'etral}. Si les comoules
proviennent d'un op\'erateur de
coloi ii), alors on a un \'el\'ement primitif si le moule est {\it
alternel} et ``group-like" si
le moule est {\it sym\'etrel}. \\

On note $M(\Omega )$ l'ensemble des moules sur $\Omega$, $M_{alta} (\Omega )$
(resp. $M_{alte} (\Omega )$) l'ensemble des moules alternals (resp.
alternels) et
$M_{syma} (\Omega )$ (resp. $M_{syme} (\Omega )$) l'ensemble des moules
sym\'etrals (resp. sym\'etrels). Nous avons le diagramme suivant :
$$
\left .
\begin{array}{lll}
\di M^{\bullet} \in M_{alta} (\Omega ) &
\di\stackrel{\exp}{\rightleftharpoons}\atop{\log} & N^{\bullet} \in
M_{syma} (\Omega ) \\
\Phi \, \downarrow & & \\
M^{\bullet} \in M_{alte} (\Omega ) &
\di\stackrel{\exp}{\rightleftharpoons}\atop{\log} & N^{\bullet } \in
M_{syme} (\Omega ) .
\end{array}
\right .
$$
Les diff\'erentes op\'erations sur les moules, notamment les
op\'erations d'addition, multiplication et
composition, se d\'eduisent des op\'erations correspondantes sur les
op\'erateurs formels associ\'es. On
d\'emontre ainsi que $M_{alta(e)}$ muni de la composition et
$M_{sima(e)}$ muni de la multiplication sont
des groupes. \\

Jean \'Ecalle d\'egage donc de ces questions, une alg\`ebre, appel\'ee
{\it alg\`ebre des moules}, et qui,
suivant les remarques pr\'ec\'edentes, interviendra dans toutes les
questions d'ordre {\it constructif}
sur les big\`ebres gradu\'ees $({\cal B},\Delta )$ et $({\cal
B},\Delta_* )$. L'alg\`ebre
des moules que nous pr\'esentons formellement ici, sera d\'etaill\'ee et
explicit\'ee dans le reste du texte. \\

{\bf Alg\`ebre des moules}. {\it Soit $M(\Omega )$ l'ensemble des moules
d\'efinis sur $\Omega$, muni des op\'erations\\

i) addition : $A^{\bullet} =M^{\bullet} +N^{\bullet}$ $\equiv$
$A^{\underline{\omega} }
=M^{\underline{\omega}} +N^{\underline{\omega}}$, \\

ii) multiplication : $A^{\bullet} =M^{\bullet} \times N^{\bullet}$
$\equiv$ $A^{\underline{\omega}}
=\di\sum_{\underline{\omega}^1 \bullet \underline{\omega}^2
=\underline{\omega}} N^{\underline{\omega}^1 }
M^{\underline{\omega}^2 }$, \\

\ni forme une alg\`ebre non commutative. \\

Si $\Omega$ a une structure de semigroupe, l'alg\`ebre des moules, peut
\^etre munie d'une loi de composition,
not\'ee $\circ$, compatible avec les lois $+$ et $\times$, et d\'efinie
par \\

iii) composition : $A^{\bullet} =M^{\bullet} \circ N^{\bullet}$
$\equiv$ $A^{\underline{\omega}} =
\di\sum_{s\geq 0, \omega^1 \dots \omega^s =\omega} M^{\parallel
\omega^1 \parallel ,\dots ,\parallel \omega^s \parallel}
N^{\omega^1 } \dots N^{\omega^s}$, o\`u $\parallel \underline{\omega}
\parallel =\di\sum_{i=1}^n \omega_i$. \\

L'alg\`ebre des moules muni des op\'erations $(+,\times ,\circ )$ est
une alg\`ebre \`a composition.} \\

Par ailleurs, on a les r\'esultats de stabilit\'e suivants sur les
moules : \\

{\bf Propri\'et\'es de stabilit\'e des moules}. {\it On note $S^{a(e)l}$
un moule sym\'etra(e)l et
$A^{a(e)l}$ un moule alterna(e)l. On a : \\

i) $S^{a(e)l} \times S^{a(e)l} =S^{a(e)l}$, \\

ii) $(S^{a(e)l} )^{-1} \times A^{a(e)l} \times S^{a(e)l} =A^{a(e)l}$, \\

iii) $A^{a(e)l} \circ A^{a(e)l} = A^{a(e)l}$.} \\

Il y a d'autres propri\'et\'es qui seront d\'etaill\'ees par la suite.

\section*{Plan du m\'emoire}

Ce m\'emoire est constitu\'e de $5$ parties.\\

La partie $1$ consiste en quelques rappels de la th\'eorie des alg\` ebres de lie,
des alg\`ebres enveloppantes d'alg\`ebres de lie, et les notions importantes de cog\`ebres
et big\`ebres, qui seront au coeur de ce travail.\\

La partie $2$ est enti\`erement consacr\'ee au calcul moulien, i.e. \`a la combinatoire de
deux big\`ebres, not\'ees $\aA$ et $\eE$, qui interviennent dans tous les travaux d'Ecalle.
On en donne deux sources naturelles, \`a savoir les d\'erivations sur une alg\`ebre, et
les op\'erateurs diff\'erentiels. On y d\'efinit aussi les principales sym\'etries des moules,
et leurs traductions dans la terminologie des alg\`ebres de lie libres.\\

La partie $3$ introduit l'op\'eration de {\it composition} des moules, qui est l'analogue non commutatif de la substitution
des s\'eries formelles. On d\'ecrit en d\'etail la structure d'{\it alg\`ebre \`a composition} ainsi obtenue, ainsi que
diverses structures associ\'ees, comme les groupes alerna(e)ls et symetra(e)ls.\\

La partie $4$, qui peut \^etre omise dans une premi\`ere lecture, tra\^\i te de questions th\'eoriques sur les
moules. On donne une interpr\'etation des sym\'etries secondaires alternil/symetril. On d\'emontre que certaines
sym\'etries de moules peuvent s'\'etablir sans calculs, si ces moules v\'erifient une \'equation diff\'erentielle
donn\'ee.\\

La partie $5$ enfin, discute la construction des formes normales d'objets analytiques locaux
(champs de vecteurs et diff\'eomorphismes).

\section*{Remerciements}

Je remercie Jean \'Ecalle pour ses commentaires, et les fichiers qu'il m'a envoy\'e
dans lesquels j'ai largement puis\'e des exemples de Moules. Je remercie
\'egalement Leila Schneps, pour sa relecture critique du manuscrit et son aide dans
la refonte de certains passages. Guilaume Morin m'a beaucoup facilit\'e la tache en reprenant la r\'edaction du
dernier chapitre et en r\'edigeant mes expos\'es oraux et mes notes sur la m\'ethode d'arborification pour son m\'emoire de DEA.
Enfin, j'exprime ma gratitude \`a tout ceux qui m'ont encourag\'e dans cette r\'edaction, notamment Jean-nicolas
D\'enari\'e, Herbert Gangle, Pierre Lochak, Michel Petitot, George Racinet, Pierre Cartier et Michel Waldschmidt.

\include{moule0}
\include{moule1}
\include{moule12}
\include{moule2}
\include{moule3}

\part*{Notations}

Soit $\kK$ un anneau. On note:\\

- $\kK [x]$ l'ensemble des polyn\^omes \`a coefficients dans $\kK$.

- $\kK [[x]]$ l'ensemble des s\'eries formelles \`a coefficients dans $\kK$.

- $\kK \{ x\}$ les s\'eries analytiques.\\

Soit X un alphabet.\\

- $A_X$ alg\`ebre libre sur $X$.

- ${\cal L}_X$ l'alg\`ebre de Lie libre construite sur $X$.

- ${\cal M}_X$ l'id\'eal de $\kK [[x]]$ form\'e des s\'eries formelles sans terme constant.

- ${\rm U}{\cal L}_X$ alg\`ebre enveloppante de l'alg\`ebre de Lie ${\cal L}_X$.

- $X^*$ l'ensemble des mots construits sur $X$.

- $X^{*,r}$ l'ensemble des mots de longueur $r \geq 0$.

- ${\cal M}_{\kK} (X)$ l'ensemble des moules sur $X$ \`a valeurs dans $\kK$.

- $\sh$: application de battage.

-$\csh$: application de battage contractant.

- $\Delta$: coproduit sur $\aA$.

- $\Delta_*$: coproduit sur $\eE$.\\

- $M^{\bullet}$: le moule $M$.

- $M^{\bs}$: le moule $M$ \'evalu\'e sur la suite $\bs$.

- $M^{\bullet^{<}}$: l'arborifi\'e du moule $M^{\bullet}$.

\end{document}

%% file: avant-propos.tex
\part*{Avant-propos sur le calcul moulien}

\centerline{{\it par} Jean Ecalle}
\vskip 1cm
Les {\it moules} sont des objets on ne peut plus concrets et banals : ce sont de
simples fonctions d'un ``nombre variable de variables"; ou si l'on
pr\'ef\`ere,  des fonctions d\'efinies sur un monoide. Mais
l\`a-dessus viennent se greffer:
\bigskip
\\ \noindent
(i) trois op\'erations de base, plus une douzaine de secondaires.
\\ \noindent
(ii) quatre grands types de ``sym\'etrie" ou d'``alternance", plus
une douzaine de secondaires
\\ \noindent
(iii) une batterie de r\`egles et recettes simples, qui disent
comment telle ou telle op\'eration affecte, conserve, transforme,
etc\dots, telle ou telle propri\'et\'e
\\ \noindent
(iv) une transformation de grande port\'ee, l'{\it arborification
}, qui sert  surtout \`a rendre {\it convergentes} des s\'eries
mouliennes divergentes, mais qui poss\`ede aussi la propri\'et\'e
inattendue   de ``respecter" l'expression analytique des
principaux moules utiles
\\ \noindent
(v) et enfin, bien s\^ur, un bestiaire de quelque trente {\it
moules fondamentaux}, qui surgissent et resurgissent un peu
partout, soit directement, soit comme {\it ingr\'edients} ou {\it
pi\`eces d\'etach\'ees} \`a partir desquelles sont construits les
{\it moules secondaires}, eux-m\^emes en quantit\'e ind\'efinie.
\bigskip
\\ \noindent
Aussi \'elabor\'e que puisse para\^\i tre cet appareil, il reste
malgr\'e tout
  d\'ecid\'ement \'el\'ementaire
dans ses ressorts. Aussi serait-il trompeur, \`a mon avis , de
parler d'une {\it th\'eorie des moules}. On serrerait sans doute
la v\'erit\'e de plus pr\`es en parlant \`a leur propos d'un {\it
syst\`eme de notations } doubl\'e d'un {\it mode d'emploi
sophistiqu\'e}, qui permet souvent de poursuivre les calculs
m\^eme l\`a o\`u la complexit\'e des expressions \`a manier
semble redhibitoire. On pourrait aussi parler d'un {\it \'etat
d'esprit moulien}\,: c'est la mentali\'e de celui qui   ne se
contente pas de th\'eor\`emes g\'en\'eraux d'existence,
d'unicit\'e, etc, nous laissant sur notre faim, mais
  qui d\'elib\'er\'ement {\it recherche l'explicite}, car il sait
par exp\'erience que c'est presque toujours possible,  toujours
payant, et souvent indispensable d\`es qu'on vise des r\'esultats
tant soit peu pr\'ecis. Et l'on pourrait ajouter que les moules
s'incrivent dans une d\'emarche {\it typiquement analytique} : ils
permettent en effet de cerner les difficult\'es, puis de les
s\'erier, puis de les vaincre , en les examinant  tour \`a tour
pour les composantes de longueur 1, 2, 3, etc,  jusqu'\`a ce que
les m\'ecanismes en jeu se d\'evoilent et livrent la solution
g\'en\'erale.
\bigskip
\\ \noindent
C'est pr\'ecis\'ement cette d\'emarche qui a permis, en th\'eorie
KAM, de dissiper la chim\`ere des {\it petits diviseurs
surmultiples}, qui n'ont aucune esp\`ece d'existence, mais qui
hantaient la th\'eorie depuis son origine.
\\ \noindent
La m\^eme d\'emarche s'applique, avec le m\^eme succ\`es, \`a
l'analyse des {\it objets analytiques locaux} ( champs de
vecteurs, diff\'eomorphismes, \'equations ou syst\`emes
diff\'erentiels ou fonctionnels \dots ) et en particulier \`a
l'\'etude de leurs {\it invariants holomorphes}. Ces derniers sont
souvent r\'eput\'es ``non calculables", alors qu'ils le sont
\'eminemment - gr\^ace aux moules.
\\ \indent
Les moules interviennent aussi en th\'eorie de la {\it
r\'esurgence}, o\`u d'ailleurs ils prennent  leur origine, car
c'est l\`a un contexte typiquement non commutatif, qui \`a chaque
pas requiert des indexations sur le monoide engendr\'e par
$\Bbb{C}$.
\\ \noindent
Il y a aussi tout le champ des {\it fonctions sp\'eciales} et sa
``compl\'etion naturelle", qui est justement le champ des {\it
moules sp\'eciaux}. Expliquons-nous. L'Analyse du $19^{eme}$
si\`ecle avait pour id\'eal la r\'esolution explicite des
\'equations (diff\'erentielles, etc) au moyen d'un certain nombre
de {\it fonctions sp\'eciales}, r\'epertori\'ees, d\'ecrites et
tabul\'ees une fois pour toutes. Mais cela s'est vite r\'ev\'el\'e
impraticable, car aucune collection de fonctions sp\'eciales n'y
suffisait. Aussi l'optique a-t-elle chang\'e et, pour la {\it
common wisdom} du $20^{eme}$ si\`ecle, le `but' \'etait au
contraire  de trouver des {\it algorithmes} de r\'esolution.
C'\'etait un progr\`es, mais un recul aussi : on perdait en
transparence ce qu'on gagnait en g\'en\'eralit\'e. Heureusement,
les deux choses sont conciliables : si le champ des fonctions
sp\'eciales est trop restreint pour ``tout exprimer", le champ des
moules sp\'eciaux, lui, y suffit - tout en incorporant l'aspect
{\it algorithmique}, vu le mode  de d\'efinition, par r\'ecurrence
sur la longueur, de la plupart des moules sp\'eciaux.
\\ \noindent
Qui dit {\it fonctions  sp\'eciales} dit aussi {\it constantes
trancendantes sp\'eciales} : les deux choses vont de pair. L\`a
aussi, les moules sont l'outil idoine : ils sont le langage
pr\'eadapt\'e dans lequel se construit et s'\'etudie le corps
d\'enombrable $\Bbb{N}\mathrm{a}$ des {\it naturels}, qui contient
(presque) toutes les {\it constantes transcendantes naturelles} -
\`a commencer par les {\it multizetas}, pour qui les principales
conjectures viennent justement d'\^etre r\'esolues, par une
d\'emarche qui, du d\'ebut \`a la fin, utilise les
notations et les op\'erations  mouliennes.

%% file: moule0.tex
\setcounter{section}{0}
\setcounter{thm}{0}
\setcounter{lem}{0}
\setcounter{defi}{0}
\setcounter{rema}{0}
\setcounter{equation}{0}

\part{Pr\'eliminaires}

Cette partie consiste en quelques rappels succincts sur les alg\`ebres de Lie libres, les notions de big\`ebres et
de cog\`ebres.

\section{Alg\`ebre associative et alg\`ebre de lie libre}
\label{alass}

On renvoie au livre de Serre (\cite{Se}, LA, Chap.4) et Bourbaki (\cite{Bo1}, chap.2) pour les
d\'emonstra\-tions des r\'esultats rappel\'es dans le $\S$\ref{alass}.\\

On rappelle qu'un magma libre est un ensemble $M$ muni d'une application \break $M\times M\rightarrow M$, not\'ee
$(x,y)\rightarrow xy$.

\begin{defi}
Soit $X$ un ensemble fini. On d\'efinit par r\'ecurrence une famille d'en\-sembles ${\mathcal X}_n$,
$n\geq 1$, tels que\\

i) ${\mathcal X}_1 = X$,

ii) pour $n\ge 2$,
\begin{equation}
{\mathcal X}_n =\coprod_{{{p+q =n}\atop{p,q\ge 1}}} {\mathcal X}_p
\times {\mathcal X}_q,
\end{equation}

o\`u $\coprod$ d\'enote la r\'eunion disjointe.
\end{defi}

\begin{rema}
Un \'el\'ement de ${\mathcal X}_2$ est un couple $(x,y)$ avec $x,y\in X$; un \'el\'ement de ${\mathcal X}_3$ sera
de la forme $(x,(y,z))$ ou bien $((x,y),z)$, etc.  De mani\`ere g\'en\'erale, l'ensemble ${\mathcal X}_n$ est
l'ensemble des mots parenth\'es\'es de longueur $n$.
\end{rema}

On note $M_X =\di\coprod_{n =1}^{\infty } {\mathcal X}_n$, et on
d\'efinit la multiplication
\begin{equation}
\left .
\begin{array}{lll}
M_X \times M_X & \rightarrow & M_X ,\\
x_p \times x_q & \rightarrow & (x_p,x_q) \in {\mathcal
X}_{p+q}\subset M_X  ,
\end{array}
\right .
\end{equation}
o\`u $x_p\in {\mathcal X}_p$, $x_q\in {\mathcal X}_q$ et la fl\`eche $\rightarrow$ d\'enote l'inclusion canonique
d\'efinie par ii). $M_X$ est un magma libre sur $X$. Un \'el\'ement $w$ de $M_X$ peut \^etre vu comme un mot non
commutatif et non associatif de $X$. Sa longueur $l(w)$ est l'unique $n$ tel que $w\in  {\mathcal X}_n$. \\

Soit $\kK$ un corps, et soit $A_X$ la $k$-alg\`ebre du magma libre
$M_X$; les \'el\'ements $\alpha \in  A_X$ sont les sommes finies
\begin{equation}
\alpha = \di\sum_{m\in  M_X} c_m \cdot m,\ \ c_m \in  \kK .
\end{equation}
La multiplication dans $A_X$ \'etend la  multiplication dans $M_X$; on continue donc \`a la noter $(\cdot,\cdot)$.
L'alg\`ebre $A_X$ est appel\'ee {\it  l'alg\`ebre libre} $A_X$ sur $X$. \\

Soit $I$ l'id\'eal de $A_X$ engendr\'e par les \'el\'ements de la forme $(a,a),\ a\in A_X$, et par ceux de la forme
$J(a,b,c)=((a,b),c)+((b,c),a)+((c,a),b)$ avec $a, b, c\in A_X$. L'alg\`ebre quotient $A_X /I$ est appel\'e
{\it  alg\`ebre de Lie libre} sur $X$ et se note ${\cal L}_X$. \\

Soit $U{\cal L}_X$ l'{\it alg\`ebre enveloppante universelle} de ${\cal L}_X$. Cette alg\`ebre est isomorphe \`a
l'alg\`ebre $Ass_X$ de ``polyn\^omes associatifs mais non commutatifs sur $X$''.\\

Nous utilisons la notation suivante pour ces polyn\^omes. L'ensemble $X$ \'etant suppos\'e fini, posons
$X=\{ X_1 ,\dots ,X_N\}$. Pour tout $r\geq 1$, on note $\Omega^*$ l'ensemble des suites
$\omega =(\omega_1,\ldots, \omega_r)$ avec $\omega_i\in  \{1,\ldots,N\}$ pour $1\le i\le r$. Pour chaque $r$,
on note $\Omega^{*,r}$ le sous-ensemble de $\Omega^*$ constitu\'e des suites de longueur $r$. A toute suite
$\omega=(\omega_1,\ldots,\omega_r)\in \Omega^*$, on associe un mot de $Ass_X$, \`a savoir
$X_\omega:=X_{\omega_1 } \dots  X_{\omega_r }$.

Un \'el\'ement $m\in  Ass_X$ s'\'ecrit alors
\begin{equation}
\label{polyno}
m=\di\sum_{\omega \in  \Omega} M^{\omega } X_{\omega},
\end{equation}
avec les $M^\omega\in \kK$ presque tous \'egaux \`a $0$.\\

L'inclusion de ${\cal L}_X$ dans son alg\`ebre enveloppante universelle $U{\cal L}_X$ donne une inclusion ${\cal
L}_X\rightarrow Ass_X$ puisque $U{\cal L}_X\simeq Ass_X$; ce morphisme est donn\'e explicitement par
\begin{equation}
\begin{array}{rcl}
{\cal L}_X &\hookrightarrow & Ass_{X}\cr X_i &\mapsto & X_i\cr
[X_i,X_j] & \mapsto & X_iX_j-X_jX_i.
\end{array}
\end{equation}

Le produit direct ${\cal L}_X\times {\cal L}_X$ est aussi une alg\`ebre de Lie, munie du produit de Lie donn\'e par
\begin{equation}
\Bigl[(x,x'),(y,y')\Bigr]=\Bigl([x,y],[x',y']\Bigr).
\end{equation}
On a l'homomorphisme diagonal d'alg\`ebres de Lie
\begin{equation}
\begin{array}{ccc}
{\cal L}_X & \rightarrow & {\cal L}_X\times {\cal L}_X\\
x & \mapsto & (x,x).
\end{array}
\end{equation}

Comme on a un isomorphisme
\begin{equation}
U({\cal L}_X\times{\cal L}_X)\buildrel\sim\over\rightarrow U{\cal L}_X\otimes U{\cal L}_X
\end{equation}

donn\'e par
\begin{equation}
(x,y)\mapsto x\otimes 1+1\otimes y,
\end{equation}

l'homomorphisme diagonal induit un homomorphisme appel\'e coproduit:
\begin{equation}
\begin{array}{ccc}
\Delta:U{\cal L}_X &\rightarrow &U{\cal L}_X\otimes U{\cal L}_X\\
x & \mapsto & x\otimes 1+1\otimes x.
\end{array}
\end{equation}

Comme $U{\cal L}_X\simeq Ass_X$, nous avons donc un coproduit
\begin{equation}
\Delta:Ass_X\rightarrow Ass_X\otimes Ass_X.
\end{equation}

Le r\'esultat important suivant caract\'erise de mani\`ere
naturelle la sous-alg\`ebre de Lie ${\cal L}_X$ \`a l'int\'erieur
de $Ass_X$:

\begin{thm}
\label{thm1}
Soit $X$ un ensemble fini; alors l'alg\`ebre de Lie libre ${\cal L}_X$ sur $X$ coincide avec l'ensemble des
\'el\'ements {\rm primitifs} de $Ass_X$, c'est \`a dire
\begin{equation}
{\cal L}_X =\{ m\in  Ass_X ,\  \Delta m =m \otimes  1 +1 \otimes m \}.
\end{equation}
\end{thm}

Soit ${\mathcal M}$ l'id\'eal de l'alg\`ebre $Ass_X$ engendr\'e par l'ensemble $X$, c'est-\`a-dire l'id\'eal de tous
les polynomes sans terme constant; il est engendr\'e par les mon\^omes (non commutatifs).\\

On d\'efinit une application
\begin{equation}
\psi:{\mathcal M} \rightarrow {\cal L}_X
\end{equation}
en posant
\begin{equation}
\psi (X_{\omega_1}\cdots X_{\omega_r})= {1\over r}
[[\cdots[[X_{\omega_1},X_{\omega_2}],X_{\omega_3}] \cdots,
X_{\omega_{r-1}}], X_{\omega_r}]
\end{equation}

pour les mon\^omes et en l'\'etendant par lin\'earit\'e \`a tout ${\mathcal M}$.\\

Comme ${\cal L}_X\subset Ass_X$, et d'ailleurs m\^eme ${\cal L}_X\subset{\mathcal M}$, on peut restreindre $\psi$ \`a ce
sous-ensemble, et on obtient le r\'esultat suivant, appel\'e {\it th\'eor\`eme de projection}:

\begin{thm}
\label{projec}
L'application $\psi$ est une retraction de ${\mathcal M}$ sur ${\cal L}_X$, i.e. on a $\psi_{\mid {\cal L}_X
} =id_{{\cal L}_X }$.
\end{thm}

Par exemple, on sait que l'\'el\'ement $X_1X_2-X_2X_1\in {\cal M}$ appartient en fait \`a ${\cal L}_{\cal X}$, puisqu'il
s'agit de $[X_1,X_2]$; on a bien
\begin{equation}
\psi(X_1X_2-X_2X_1)={{1}\over{2}}[X_1,X_2]-{{1}\over{2}}[X_2,X_1]=[X_1,X_2].
\end{equation}

Les alg\`ebres ${\cal L}_X$ et $Ass_X$ sont munies de graduations naturelles par la longueur; on appelle
${\cal L}_X^n$ resp. $Ass_X^n$ l'ensemble des combinaisons lin\'eaires de crochets (resp. de mots) homog\`enes de
longueur $n$.\\

On a donc
\begin{equation}
{Ass}_X =\bigoplus_{r=0}^{\infty } Ass_X^n ,\ \ \ {L}_X
=\bigoplus_{r=0}^{\infty}{\cal L}_X^n.
\end{equation}

On prend les compl\'et\'es par rapport \`a ces graduations, que l'on note
\begin{equation}
\hat{Ass}_X =\prod_{r=0}^{\infty } Ass_X^n ,\ \ \ \hat{L}_X
=\prod_{r=0}^{\infty}{\cal L}_X^n;
\end{equation}

un \'el\'ement $m$ d'un tel compl\'et\'e peut \^etre repr\'esent\'e par une serie formelle
\begin{equation}
m =\sum_{r=1}^\infty Y_r =\di\sum_{r=0}^\infty \Bigl(\sum_{\omega \in  \Omega_r}M^{\omega }
X_{\omega }\Bigr) ,
\end{equation}

avec $Y_r\in Ass_X^r$; nous noterons cette somme de mani\`ere condens\'ee
\begin{equation}
m=\di\sum_{\bullet } M^{\bullet } X_{\bullet } .
\end{equation}

Un \'el\'ement de $\hat{Ass}_X$ est enti\`erement d\'etermin\'e par
la donn\'ee de ses coefficients $\di M^{\bullet }$.\\

Soit $\hat{\mathcal M}\in \hat{Ass}_X$ l'id\'eal engendr\'e par l'ensemble fini $X$; c'est donc l'id\'eal des s\'eries
formelles sans terme constant. On d\'efinit les applications
\begin{equation}
\exp :\hat{\mathcal M} \rightarrow 1+\hat{\mathcal M},\ \ \ \ \ \log :
1+\hat{\mathcal M} \rightarrow \hat{\mathcal M },
\end{equation}

par les formules usuelles
\begin{equation}
\exp (x)=\di\sum {x^n \over n!},\qquad \ \log (1+x)=\di\sum_{n=1}^{\infty } (-1)^{n+1} \di {x^n \over n!} .
\end{equation}

On a
\begin{equation}
\exp \log =\log \exp =1 .
\end{equation}

On d\'efinit le {\it produit tensoriel compl\'et\'e}
\begin{equation}
\hat{Ass}_X\hat\otimes\hat{Ass}_X=\prod_{p,q}Ass_X^p\otimes = Ass_X^q,
\end{equation}

et on note que le coproduit $\Delta$ s'\'etend aux compl\'et\'es, ainsi que le th\'eor\`eme \ref{thm1}, c'est-\`a-dire que
${\cal L}_X$ s'identifie \`a l'int\'erieur de $\hat{Ass}_X$ avec l'ensemble des $m\in\hat{Ass}_X$ tels que
$\Delta(m)=m\otimes 1+1\otimes m$ dans $\hat{Ass}_X\hat\otimes\hat{Ass}_X$.  On a le r\'esultat important suivant:

\begin{thm}
\label{repres}
L'application $\exp$ d\'efinit une bijection de l'ensemble des \'el\'ements primitifs $\alpha \in \hat{\mathcal
M}$, c'est-\`a-dire tels que $\Delta \alpha =\alpha \otimes  1 + 1 \otimes  \alpha$, vers l'ensemble des
$\beta \in  1+ \hat{\mathcal M}$ tels que $\Delta \beta =\beta \otimes  \beta$.
\end{thm}

\section{Cog\`ebres et big\`ebres}

\subsection{Cog\`ebres et coproduit}
\label{cogebre}

Soit $\kK$ un corps.

\begin{defi}
\label{def4}
On appelle cog\`ebre sur $\kK$ un triplet $(E,\Delta,\varepsilon)$ ou $E$ est un $\kK$-espace vectoriel,
$\Delta$ est une application $\kK$-lin\'eaire $\Delta:E\rightarrow E\otimes_k E$, dit coproduit de $E$,
et $\varepsilon:E\rightarrow \kK$ est une application $\kK$-lin\'eaire, tels que
\begin{eqnarray}{c}
({\rm id}_E\otimes\Delta)\circ\Delta=(\Delta\otimes{\rm id}_E)\circ\Delta \\
({\rm id}_E\otimes \varepsilon)\circ\Delta=(\varepsilon\otimes{\rm id}_E)\circ\Delta={\rm id}_E.
\end{eqnarray}
\end{defi}

On donne un exemple dans la prochaine section.

\subsection{Big\`ebres et graduations}

Une big\`ebre est un quintuplet $(E,\mu,\eta,\Delta,\varepsilon)$ tel que \\

\noindent (1) $(E,\mu,\eta)$ est une alg\`ebre, i.e. $E$ est un
$\kK$-espace vectoriel, et $\mu:E\otimes_{\kK} E\rightarrow E$ et
$\eta:\kK\rightarrow E$ sont des applications $\kK$-lin\'eaires
telles que
\begin{eqnarray}
\mu\circ(\mu\otimes 1)=\mu\circ(1\otimes \mu)&\!\!:E\otimes_{\kK} E\otimes_{\kK} E\rightarrow E ,\\
\mu\circ (1\otimes \eta)=\mu\circ (\eta\otimes 1)=1&\!\!:E\rightarrow E ,
\end{eqnarray}

en identifiant $\kK\otimes_{\kK} E=E\otimes_{\kK} \kK=E$.\\

\noindent (2) $(E,\Delta,\varepsilon)$ est une cog\`ebre (voir \S \ref{cogebre}) \\

\noindent (3) $\Delta$ et $\varepsilon$ sont des homomorphismes d'alg\`ebres. \\

Une big\`ebre est {\it gradu\'ee} si l'alg\`ebre sous-jacente est munie d'une graduation.

\begin{defi}
\label{defiprimigroup}
Dans toute la suite, un \'el\'ement $x$ d'une big\`ebre $E$ munie d'un coproduit $\Delta$ sera dit {\it primitif} si
\begin{equation}
\Delta x = x\otimes 1 +1\otimes x ,
\end{equation}
et ``{\it group-like}" si
\begin{equation}
\Delta x =x\otimes x.
\end{equation}
\end{defi}

\subsection{Exemples}
\label{exemples}

\subsection{Alg\`ebre de d\'erivations}
\label{algederi}

On renvoie au livre de Jacobson (\cite{Ja}, chap.1, \S 2, p.7-8) pour plus de d\'etails. \\

Soit $A$ une $\kK$-alg\`ebre non associative quelconque. Une {\it d\'erivation} $D$ (ou {\it op\'erateur
diff\'erentiel d'ordre 1}) sur $A$ est une  application lin\'eaire de $A$ dans $A$ satisfaisant
\begin{equation}
D(xy)=(Dx)y + x(Dy).
\end{equation}

On note ${\rm Der} (A)$ l'ensemble des d\'erivations sur $A$.\\

Soient $D_1$, $D_2 \in {\rm Der} (A)$; on a
\begin{equation}
\left .
\begin{array}{lll}
(D_1 +D_2 ) (xy) & = & D_1 (xy)+D_2 (xy) , \\

& = & (D_1 x) y +x(D_1 y) +(D_2 x)y +x (D_2 y) , \\

& = &  (D_1 +D_2)x y +x(D_1 +D_2 )y .
\end{array}
\right .
\end{equation}

Par cons\'equent $D_1 +D_2 \in {\rm Der} (A)$. De m\^eme, on a $\alpha= D_1 \in {\rm Der} (A)$ si $\alpha \in \kK$.
On voit donc que ${\rm Der}(A)$ est un $\kK$-espace vectoriel.\\

On peut composer les d\'erivations; la composition de $D_2$ et $D_1$ est donn\'ee par
\begin{equation}
\left .
\begin{array}{lll}
D_2 D_1 (xy) & = & D_2 ( (D_1 x)y +x (D_1 y)) , \\
& = & (D_2 D_1 x) y + D_1 x D_2 y +D_2 x D_1 y +x D_2 D_1 y .
\end{array}
\right .
\end{equation}

On peut m\^eme donner la formule g\'en\'erale pour l'action d'un op\'erateur diff\'erentiel
$B=D_1\circ\cdots \circ D_r$ d'ordre $r$ sur $xy$; elle n'est pas bien belle, mais elle peut \^etre
utile.

\begin{defi}
\label{defipr}
Pour chaque entier $r\ge 1$, nous introduisons l'ensemble $P_r$ de paires
\begin{equation}
\Bigl((i_1,\ldots,i_n),(j_1,\ldots,j_m)\Bigr)
\end{equation}

de suites d'entiers, sous-ensembles de la suite $(1,\ldots,r)$, telles que
\begin{equation}
\begin{array}{l}
n+m=r\\
1\le i_1<i_2<\cdots<i_n\le r\\
1\le j_1<j_2< \cdots <j_m\le r\\
\{i_1,\ldots,i_n\}\cap\{j_1,\ldots,j_m\}=\emptyset\\
\{i_1,\ldots,i_n\}\cup\{j_1,\ldots,j_m\}=\{1,\ldots,r\}.
\end{array}
\end{equation}

Les paires $\Bigl(\emptyset,(1,\ldots,r)\Bigr)$ et $\Bigl((1,\ldots,r),\emptyset\Bigr)$ sont incluses dans $P_r$.
Nous \'ecrirons souvent $\underline{i}$ pour $(i_1,\ldots,i_n)$ et $\underline{j}$ pour $(j_1,\ldots,j_m)$), donc
$(\underline{i},\underline{j})\in P_r$.
\end{defi}

L'action d'un op\'erateur diff\'erentiel $B$ est donn\'ee
explicitement par
\begin{equation}
\label{sum}
B(xy)=\sum_{(\underline{i},\underline{j})\in P_r}(B_{\underline{i}}x)(B_{\underline{j}}y)
\end{equation}

o\`u
\begin{equation}
B_{\underline{i}}=D_{i_1}D_{i_2}\cdots D_{i_n}\ \ {\mbox{et}}\ \
B_{\underline{j}}=D_{j_1}D_{j_2}\cdots D_{j_m}.
\end{equation}

La composition de $r$ d\'erivations forme un op\'erateur diff\'erentiel d'ordre $r$, qui n'est pas une d\'erivation si
$r>1$; l'ensemble des d\'erivations ne forme donc pas une alg\`ebre pour cette multiplication.\\

Par contre, en posant $[D_1,D_2]=D_1D_2-D_2D_1$, on obtient
\begin{equation}
[D_1,D_2]xy=\bigl([D_1,D_2]x\bigr)y +x\bigl([D_1,D_2]y\bigr) .
\end{equation}

Par cons\'equent, $[D_1,D_2 ] \in {\rm Der} (A)$ et $[\ ,\ ]$ fait de ${\rm Der} (A)$ une alg\`ebre de Lie appel\'ee
{\it alg\`ebre de Lie des d\'erivations} ou {\it alg\`ebre de  d\'erivations} de $A$.\\

Soit $F(A)={D_1 ,\dots ,D_i }$, une famille finie de d\'erivations de ${\rm Der}(A)$. On note
$\kK\la\la F(A)\ra\ra$ l'anneau des s\'eries formelles sur l'alphabet $F(A)$. Pr\'ecisons que nous
consid\'erons l'identit\'e ${\rm id}_A$ comme l'unique op\'erateur diff\'erentiel d'ordre $0$, et nous l'identifions avec
l'\'el\'ement $1\in \kK$.  De cette mani\`ere, tous les \'el\'ements de $\kK\la\la F(A)\ra\ra$ sont des
op\'erateurs sur $A$. \\

On peut munir l'alg\`ebre $\kK\la\la F(A)\ra\ra$ d'une structure de cog\`ebre.\\

{\it Counit\'e}. On note $\varepsilon$ l'homomorphisme
\begin{equation}
\varepsilon:\kK\la\la{\rm Der}(A)\ra\ra\rightarrow \kK
\end{equation}
d\'efini en associant \`a chaque s\'erie formelle de $\kK\la\la F(A)\ra\ra$ son terme constant.\\

{\it Coproduit}. A toute d\'erivation $D\in{\rm Der}(A)$ on associe un op\'erateur de $A\otimes A$ dans
$A\otimes A$, not\'e $\bar{D}$, par
\begin{equation}
\bar{D} (\phi \otimes \psi )=D\phi \otimes \psi +\phi \otimes D\psi .
\end{equation}

Si $D_1$ et $D_2$ sont deux d\'erivations, on construit une application naturelle de $A\otimes A$ dans $A\otimes A$
not\'ee $D_1 \otimes D_2$, d\'efinie par
\begin{equation}
(D_1 \otimes D_2 )(\phi \otimes \psi )=D_1 \phi \otimes D_2 \psi .
\end{equation}

Finalement, on d\'efinit une application lin\'eaire $\Delta$ de
${\rm Der} (A)$ dans ${\rm Der} (A)\otimes {\rm Der} (A)$ par
\begin{equation}
\left .
\begin{array}{lll}
{\rm Der} (A) & \di\stackrel{\Delta}{ \rightarrow} & {\rm Der}(A)
\otimes {\rm Der} (A) ,\\ D & \mapsto & D\otimes 1 + 1 \otimes D .
\end{array}
\right .
\end{equation}

On note que
\begin{equation}
\bar{D} =\Delta (D).
\end{equation}

De plus, en notant $\mu$ l'op\'erateur lin\'eaire de $A\otimes A$ dans $A$ d\'efini par
$\phi \otimes \psi\mapsto \phi \cdot\psi$, on a l'\'egalit\'e
\begin{equation}
\label{proj1}
D\circ \mu =\mu \circ \Delta (D),
\end{equation}

traduisant la commutativit\'e du diagramme
\begin{equation}
\left .
\begin{array}{lll}
A\otimes A & \stackrel{\Delta (D)}{\longrightarrow} & A\otimes A , \\
\downarrow^{\mu}  & & \downarrow^{\mu}\\
A & \stackrel{D}{\longrightarrow} &  A.
\end{array}
\right .
\end{equation}

On peut d\'efinir directement $\Delta$ en posant
\begin{equation}
\Delta(B)=\sum_{i+j=r} B_i\otimes B_j .
\end{equation}

Notons que $\Delta(1)=1\otimes 1$, et de plus, on a la formule utile
\begin{equation}
\label{prod}
\Delta(B_1B_2)=\Delta(B_1)\Delta(B_2).
\end{equation}

Pour le voir, on commence \`a l'ordre 2:
\begin{equation}
\begin{array}{rcl}\Delta(D_1D_2)&=&(D_1D_2)\otimes 1+D_1\otimes
D_2+D_2\otimes D_1
+1\otimes (D_1D_2)\\
&=&(D_1\otimes 1+1\otimes D_1)(D_2\otimes 1+1\otimes D_2)\\
&=&\Delta(D_1)\Delta(D_2).
\end{array}
\end{equation}

et on continue par r\'ecurrence.

\begin{lem}
\label{cogeb}
Le triplet ${\cal D}=(\kK\langle\langle {\rm Der}(A)\rangle\rangle, \Delta,\varepsilon)$ est une cog\`ebre.
\end{lem}

\begin{proof}
Il suffit de v\'erifier les axiomes de la d\'efinition \ref{def4}, i.e.
\begin{equation}
\begin{array}{rl}
({\rm id}_{\cal D}\otimes
\varepsilon)\circ\Delta=(\varepsilon \otimes{\rm id}_{\cal
D})\circ\Delta={\rm id}_{\cal D}&\!\!:{\cal
D}\rightarrow{\cal D}\\
({\rm id}_{\cal D}\otimes\Delta)\circ\Delta= (\Delta\otimes{\rm
id}_{\cal D})\circ\Delta&\!\!:{\cal D}\rightarrow{\cal D}
\otimes{\cal D}\otimes {\cal D}.
\end{array}
\end{equation}

o\`u ${\cal D}=\kK\la\la{\rm Der}(A)\ra\ra$.\\

Montrons d'abord qu'elles sont valables pour $1\in {\cal D}$, ensuite pour une d\'erivation $D$, puis pour un
op\'erateur diff\'erentiel d'ordre $r$, et finalement pour un objet quelconque de ${\cal D}$, i.e.
une combinaison lin\'eaire infinie des objets ci-dessus.\\

Pour $1$, en fait, c'est \'evident. Soit donc $D\in{\cal D}$ une
d\'erivation. On a
\begin{equation}
\begin{array}{rl}
({\rm id}_{\cal D}\otimes \varepsilon)\circ\Delta(D)&=
({\rm id}_{\cal D}\otimes \varepsilon)(D\otimes 1+1\otimes D)\\
 &=D\\
&= (\varepsilon \otimes{\rm id}_{\cal D})(D\otimes 1+1\otimes D)
\end{array}
\end{equation}

et
\begin{equation}
\begin{array}{rl}
({\rm id}_{\cal D}\otimes\Delta)\circ\Delta(D)&=
({\rm id}_{\cal D}\otimes\Delta)(D\otimes 1+1\otimes D)\cr
&=(D\otimes 1\otimes 1+1\otimes D\otimes 1+1\otimes 1\otimes D)\cr
&=(\Delta\otimes{\rm id}_{\cal D})(D\otimes 1+1\otimes D).
\end{array}
\end{equation}

Montrons maintenant que les relations sont satisfaites pour tout
op\'erateur diff\'erentiel $B\in \kK\la\la{\rm Der}(A)\ra\ra$.\\

Nous savons que $B=D_1\circ D_2\circ\cdots \circ D_r$; nous allons profiter de l'equation (\ref{prod}). Supposons donc
que les relations sont satisfaites pour tout op\'erateur diff\'erentiel d'ordre $< r$ (on vient de voir que c'est vrai
pour $r=1$). Posons $B'=D_1\circ D_2\circ \cdots\circ D_{r-1}$.  On a alors
\begin{equation}
\begin{array}{rl}
({\rm id}_{\cal D}\otimes
\varepsilon)\bigl(\Delta(B)\bigr)&= ({\rm id}_{\cal D}\otimes
\varepsilon)\bigl(\Delta(B')\Delta(D_r)\bigr)\\
&=({\rm id}_{\cal D}\otimes \varepsilon)\bigl(\Delta(B')\bigr)
({\rm id}_{\cal D}\otimes \varepsilon)\bigl(\Delta(D_r)\bigr)\\
&=(\varepsilon\otimes {\rm id}_{\cal D})\bigl(\Delta(B')\bigr)
(\varepsilon\otimes {\rm id}_{\cal D})\bigl(\Delta(D_r)\bigr)\\
&=(\varepsilon\otimes {\rm id}_{\cal D})\bigl(\Delta(B)\bigr),
\end{array}
\end{equation}

o\`u l'hypoth\`ese de r\'ecurrence a \'et\'e utilis\'ee dans la troisi\`eme \'egalit\'e et on utilise aussi le fait que
$(\varepsilon\otimes{\rm id}_{\cal D})$ et $({\rm id}_{\cal D}\otimes \varepsilon)$ sont des homomorphismes de
$\kK$-alg\`ebres.\\

Pour la m\^eme raison, pour la deuxi\`eme relation, on a simplement
\begin{equation}
\begin{array}{rl}({\rm id}_{\cal
D}\otimes\Delta)\bigl(\Delta(B)\bigr)&=
({\rm id}_{\cal D}\otimes\Delta)\bigl(\Delta(B')\Delta(D_r)\bigr)\\
&=({\rm id}_{\cal D}\otimes\Delta)\bigl(\Delta(B')\bigr)
({\rm id}_{\cal D}\otimes\Delta)\bigl(\Delta(D_r)\bigr)\\
&=(\Delta\otimes{\rm id}_{\cal D})\bigl(\Delta(B')\bigr)
(\Delta\otimes{\rm id}_{\cal D})\bigl(\Delta(D_r)\bigr)\\
&=(\Delta\otimes{\rm id}_{\cal D})\bigl(\Delta(B')\Delta(D_r)\bigr)\\
&=(\Delta\otimes{\rm id}_{\cal D})\bigl(\Delta(B)\bigr).
\end{array}
\end{equation}

Pour terminer la d\'emonstration du lemme, on \'etend le coproduit $\Delta$ \`a tout
$\kK\la\la{\rm Der}(A)\ra\ra$ par lin\'earit\'e.
\end{proof}

\subsection{Alg\`ebre d'op\'erateurs diff\'erentiels}

Un deuxi\`eme type de coproduit intervient naturellement dans l'\'etude des op\'erateurs diff\'erentiels d'ordre
$r\geq 1$, que nous notons $B$ (ou $B_r$) pour les distinguer des d\'erivations, habituellement not\'ees $D$.

Nous avons vu que l'ensemble des op\'erateurs diff\'erentiels d'ordre $r$ sur $A$ pour tout $r\ge 1$ forme une
$A$-alg\`ebre, not\'ee ${\rm Op}(A)$, sous la multiplication correspondant simplement \`a la composition des
op\'erateurs. A partir de $B_r$, on d\'efinit une application de $A\otimes A$ dans $A\otimes A$ par
\begin{equation}
B_r (\phi \otimes \psi )=\sum_{i+j =r} B_i \phi \otimes B_j \psi,
\end{equation}

o\`u le sens de cette somme est comme dans (\ref{sum}). On obtient un coproduit $\Delta_*$ d\'efini sur
$\kK \langle\langle{\rm Op} (A)\rangle\rangle$ par
\begin{equation}
\left .
\begin{array}{lll}
{\rm Op} (A) & \di\stackrel{\di\Delta_*}{\longrightarrow }& {\rm
Op}(A) \otimes {\rm Op} (A) ,\\ B_r & \longmapsto & \di\sum_{i+j
=r} B_i \otimes B_j  .
\end{array}
\right .
\end{equation}

On a bien s\^ur, pour tout $B\in{\rm Op}(A)$, en conservant les notations pr\'ec\'edentes
\begin{equation}
\label{proj2}
B \circ \mu =\mu \circ \Delta_* (B) .
\end{equation}

Notons que l'equation (\ref{prod}) reste valable pour $\Delta_*$.

\begin{lem}
L'espace vectoriel $\kK \la\la{\rm Op}(A)\ra\ra$, muni du coproduit $\Delta_*$ et de l'application lin\'eaire
$\varepsilon$ du lemme \ref{cogeb}, forme une cog\`ebre.
\end{lem}

\begin{proof}
Pour tout op\'erateur diff\'erentiel $B$ d'ordre $r\ge 0$, nous v\'erifions directement que
$$
\begin{array}{rl}
({\rm id}_E\otimes\Delta_*)\bigl(\Delta_*(B)\bigr)&= ({\rm id}_E\otimes\Delta_*)\bigl(\di\sum_{i+j=r}(B_i\otimes
B_j)\bigr)\\
&=\di\sum_{i+j=r}({\rm id}_E\otimes\Delta_*)\bigl(B_i\otimes B_j\bigr)\\
&=\di\sum_{i+j=r}(B_i\otimes \Delta_*(B_j)\bigr)\\
&=\di\sum_{i+j=r}(B_i\otimes \sum_{i'+j'=j}(B'_{i'}\otimes B'_{j'})\bigr)\\
&=\di\sum_{i+i'+j'=r}(B_i\otimes B'_{i'}\otimes B'_{j'}\bigr)
\end{array}
$$
et

$$\begin{array}{rl}(\Delta_*\otimes {\rm
id}_E)\bigl(\Delta_*(B)\bigr)&= (\Delta_*\otimes {\rm
id}_E)\bigl(\di\sum_{i+j=r}(B_i\otimes
B_j)\bigr)\\
&=\di\sum_{i+j=r}(\Delta_*\otimes {\rm id}_E)\bigl(B_i\otimes
B_j\bigr)\\
&=\di\sum_{i+j=r}(\Delta_*(B_i)\otimes B_j\bigr)\\
&=\di\sum_{i+j=r}\bigl(\sum_{i'+j'=i}(B'_{i'}\otimes
B'_{j'})\otimes B_j\bigr)\\
&=\di\sum_{i'+j'+j=r}(B'_{i'}\otimes B'_{j'}\otimes
B_j\bigr),
\end{array}$$

qui est \'evidemment la m\^eme chose. Pour la deuxi\`eme relation, on a
$$\begin{array}{rl}({\rm id}_E\otimes
\varepsilon)\bigl(\Delta(B)\bigr)&= ({\rm id}_E\otimes
\varepsilon)\bigl(\di\sum_{i+j=r}(B_i\otimes
B_j)\bigr)\\
&=\di\sum_{i+j=r}({\rm id}_E\otimes \varepsilon)\bigl(B_i\otimes
B_j\bigr)\\
&=B,\end{array}$$

puisque $\varepsilon(B_j)=0$ sauf si $j=0$, auquel cas $B_j=1$ et $\varepsilon(B_j)=1$.

De m\^eme, bien s\^ur, on a
$$\begin{array}{rl}
(\varepsilon\otimes {\rm id}_E)\bigl(\Delta(B)\bigr)&= (\varepsilon\otimes {\rm
id}_E)\bigl(\di\sum_{i+j=r}(B_i\otimes B_j)\bigr)\\
&=\di\sum_{i+j=r}(\varepsilon\otimes {\rm id}_E)\bigl(B_i\otimes B_j\bigr)=B,
\end{array}$$

ce qui termine la d\'emonstration.
\end{proof}

%% file: moule1.tex
\setcounter{section}{0}
\setcounter{thm}{0}
\setcounter{lem}{0}
\setcounter{defi}{0}
\setcounter{rema}{0}
\setcounter{equation}{0}

\part{Calcul Moulien}

On d\'efinit les moules et les principales op\'erations alg\'ebriques sur ces objets via les s\'eries formelles
non-commutatives associ\'ees. On donne la traduction des diverses propri\'et\'es alg\'ebriques des s\'eries
(primitive/group-like) sur les moules associ\'es, pour les big\`ebres $\aA$ et $\eE$, consuisant aux sym\'etries
alterna(e)l/symetra(e)l respectivement.

\section{Moules}

Soit $A$ un alphabet. On note $A^*$ l'ensemble des mots construits sur $A$. Les moules sont d\'efinis par:

\begin{defi}
Soit $A$ un alphabet et $\kK$ un anneau. Un moule sur $A$ \`a valeur dans $\kK$ est une application, not\'ee $M^{\bullet}$,
de $A^*$ dans $\kK$.
\end{defi}

La notation $M^{\bullet}$ pour d\'esigner une application de $A^*$ dans $\kK$ n'est pas usuelle, mais elle coincide avec la
notation utilis\'ee par Jean Ecalle.\\

Par convention, un moule $M^{\bullet}$ \'etant donn\'e, on note $M^{\ba}$ la quantit\'e $M^{\bullet} (\ba )$
pour tout $\ba \in A^*$.\\

Les moules sont \'evidemment en correspondance bi-univoque avec les s\'eries formelles non commutatives de
$\kK \langle\langle A\rangle\rangle$. En effet, a tout moule $M^{\bullet}$ sur $A$, on associe la s\'erie $S(M)=\di\sum_{\ba \in A^*}
M^{\ba} \ba$ et vice versa.\\

La structure d'alg\`ebre de $\kK \langle \langle A\rangle\rangle$ se transporte sur l'ensemble des moules sur $A$ \`a
valeur dans $\kK$, not\'e ${\cal M}_{\kK} (A)$. Nous allons voir que ce passage des s\'eries aux coefficients est souvent
int\'eressant dans de nombreuses applications. Par ailleurs, les diff\'erentes propri\'et\'es alg\`ebriques des s\'eries
$S(M)$ se lisent compl\`etement sur le moule associ\'e, ce qui justifie a posteriori la terminologie.

\section{Les big\`ebres $\aA$ et $\eE$}

Le travail d'\'Ecalle repose sur l'\'etude combinatoire de deux big\`ebres, not\'ees $\aA$ et
$\eE$\footnote{J'ai emprunt\'e cette notation \`a C. Even.}, et qui m\^\i ment pour
l'essentiel les propri\'et\'es des big\`ebres $\cal D$ et $\cal B$ des s\'eries formelles sur
les d\'erivations et les op\'erateurs diff\'erentiels respectivement.\\

Soit $\kK$ un corps. Soit $X$ un alphabet et $Y=(y_i)_{i\in \nN}$ un
alphabet cod\'e par un semi-groupe additif\footnote{Pour simplifier l'expos\'e, j'ai choisi de
coder l'alphabet par le semi-groupe $\nN$. N\'eanmoins, l'ensemble des calculs men\'es
sur $Y$ se transposent sans probl\`emes sur un alphabet $Y_{\Omega} =\{y_{\omega} \}_{\omega \in \Omega}$, o\`u $\Omega$ est
un semi-groupe quelconque. Dans ce cas d'ailleurs, on remplace la lettre $y_{\omega}$ par l'\'el\'ement du semi-groupe
$\omega$ si aucune confusion n'est possible. Dans la suite, lorsque des $\omega$ appara\^\i tront dans le context de l'alphabet $Y$, il faudra comprendre ``le
cas g\'en\'eral o\`u on a un semi-groupe pour coder l'alphabet".}.\\

On note $\kK\la\la X\ra\ra$ (resp. $\kK\la\la Y\ra\ra$), la $\kK$-alg\`ebre des
s\'eries formelles non commutatives form\'ee sur $X^*$ (resp. $Y^*$). \\

\begin{defi}
\label{d}
On note $\aA =(\kK\la\la X \ra\ra , \varepsilon ,\Delta )$ la cog\`ebre d\'efinie par:\\

i) le coproduit $\Delta$ est d\'efini sur les lettres $x$ de $X$ par
\begin{equation}
\Delta (x)=x \otimes 1+1\otimes x,\ \ \  \forall x \in X.
\end{equation}
Le coproduit d'un mot $\un{x}=x_1 \dots x_r$ se calcule via la
propri\'et\'e
\begin{equation}
\label{mor}
\Delta(\un{x}^1 \un{x}^2)=\Delta (\un{x}^1) \Delta
(\un{x}^2) ,\ \forall \un{x}^1 , \un{x}^2 \in X^*.
\end{equation}
On \'etend $\Delta$ \`a $\kK\la\la X \ra\ra$ par lin\'earit\'e, en particulier
$$\Delta(\emptyset)=\Delta(1)=1\otimes 1.$$

ii) On d\'efinit la counit\'e $\varepsilon$ par :
\begin{equation}
\left .
\begin{array}{lll}
\kK\la\la X \ra\ra & \rightarrow & \kK, \\
\\
\di\sum_{\un{x }\in X^* } M^{\un{x}}
 \un{x} & \mapsto & M^{\emptyset} .
\end{array}
\right .
\end{equation}
\end{defi}

Le coproduit $\Delta$ de $\aA$ est le coproduit naturel. La propri\'et\'e (\ref{mor}) traduit le fait que $\Delta$ est un
morphisme d'alg\`ebre.

\begin{defi}
\label{b}
Rappelons que $Y$ est un alphabet cod\'e par un semi-groupe. Notons $\eE=(\kK \la\la Y \ra\ra ,\Delta_*, \varepsilon)$ la cog\`ebre d\'efinie
par :\\

i) Le coproduit $\Delta_*$ est d\'efini sur les lettres de $Y$ par
\begin{equation}
\Delta_* (y_r) =\di\sum_{i+j=r} y_i\otimes y_j.
\end{equation}
On \'etend $\Delta_*$ \`a $Y^*$ par la propri\'et\'e:
\begin{equation}
\Delta_*(\un{y}^1 \un{y}^2 )= \Delta_* (\un{y}^1 ) \Delta_*
(\un{y}^2 ),\ \ \forall \un{y}^1 ,\un{y}^2 \in Y^* .
\end{equation}
On \'etend $\Delta_*$ \`a $\kK\la\la Y \ra\ra$ par lin\'earit\'e.

ii) la counit\'e sur $\eE$ est d\'efinie comme dans la
d\'efinition \ref{d}, ii).
\end{defi}

On voit que les cog\`ebres $\aA$ et $\eE$ sont en fait des big\`ebres:  en effet, $\kK\langle\langle X\rangle\rangle$ (resp.
$\kK\langle\langle Y\rangle\rangle$) est une alg\`ebre, $\Delta$ et $\Delta_*$ sont des homomorphismes
d'alg\`ebres, et la condition sur $\varepsilon$ est \'evidente.\\

Encore une fois, $\aA$ et $\eE$ sont des big\`ebres gradu\'ees, car les alg\`ebres
$\kK\la\la X\ra\ra$ et $\kK\langle\langle Y\rangle\rangle$ sont munies de la graduation
naturelle donn\'ee par la longueur des mots dans l'alphabet $X$ (resp. $Y$).

\section{Combinatoire sur $\aA$ et $\eE$, battage et battage contractant}

\subsection{Combinatoire sur $\aA$}
\label{combinaA}

On prend comme ci-dessus un corps $\kK$ de caract\'eristique z\'ero et un alphabet $X$. On consid\`ere la $\kK$-alg\`ebre
$\kK\langle\langle X\rangle\rangle$. Pla\c cons-nous dans la big\`ebre gradu\'ee $\aA$ obtenu en munissant
$\kK\langle\langle X \rangle\rangle$ du coproduit $\Delta$ de la d\'efinition \ref{d}.\\

Soit $X^*$ l'ensemble des mots (ou suites) de lettres de $X$, y compris le mot vide $\emptyset$.  Nous \'ecrivons
$\underline{x}=(x_1,\ldots,x_r)$ pour un \'el\'ement de $X^*$. L'entier $r$ s'appelle la longueur de la suite; si
$r=1$ alors $\un{x}\in X$, et l'unique \'el\'ement de $X^*$ de longueur $0$ est $\emptyset$. \\

Notons $C(X)$ l'ensemble des couples de suites de $X^*$; on note un couple de suites $\langle\un{x}^1;\un{x}^2\rangle
=\langle (x_1^1,\ldots,x_n^1);(x_1^2,\ldots,x_m^2)\rangle$.

\begin{defi}
\`A chaque $\un{x}\in X^*$, soit $C_{\un{x}}$ le sous-ensemble de $C(X)$ de couples de suites ``engendr\'e par
$\Delta(\un{x})$'', i.e. apparaissant dans l'expression de $\Delta(\un{x})$.
\end{defi}

\noindent $\bullet$ Pour $r=0$, on a
$$\Delta(\emptyset)=\Delta(1)=1\otimes 1=\emptyset\otimes \emptyset$$
donc
$$C_{\emptyset} = \Bigl\{\langle\emptyset;\emptyset\rangle\Bigr\}.$$

\noindent $\bullet$ Pour $r=1$, on a
$$\Delta(x)=x\otimes 1+1\otimes x$$
donc
$$C_{x} = \Bigl\{\la x;\emptyset\ra,\la\emptyset;x\ra\Bigr\}.$$
\noindent $\bullet$ Pour $r=2$, on a $\un{x}=(x_1,x_2)$ et
$$
\left .
\begin{array}{lll}
\Delta ({x_1}{x_2}) & = & \Delta({x_1})\Delta({x_2}) \\
  & = & ({x_1}\otimes 1+1\otimes {x_1})({x_2}\otimes 1 + 1 \otimes {x_2}) \\
  & = & {x_1}{x_2}\otimes 1 + {x_2}\otimes {x_1}
+ {x_1} \otimes {x_2} + 1\otimes {x_1}{x_2}\\
  & = & {\un{x}}\otimes \emptyset + {x_2}\otimes {x_1}
+ {x_1} \otimes {x_2} + \emptyset\otimes {\un{x}}.
\end{array}
\right .
$$
donc
$$C_{\un{x }=(x_1,x_2 )} =\{ \la\un{x};\emptyset\ra,
\la x_1;x_2\ra, \la x_2;x_1\ra, \la\emptyset; \un{x}\ra\},$$

\noindent $\bullet$ Pour tout $r\ge 0$, on a $|C_{\un{x}}|=2^r$ si $\un{x}$ contient $r$ composantes distinctes.
De mani\`ere g\'en\'erale, on voit que l'ensemble des couples apparaissant dans l'expression de $\Delta({\un{x}})$,
c'est-\`a-dire les couples de $C_{\un{x}}$, est donn\'e par
\begin{equation}
\label{Cw}
C_{\un{x}}=\Bigl\{\Bigl\la\underline{x}_{\underline{i}};\underline{x}_{\underline{j}}\Bigr\ra
\mid (\underline{i},\underline{j})\in P_r\Bigr\},
\end{equation}

o\`u $\un{x}=(x_1,\ldots,x_r)$ est une suite de $X^*$ de longueur $r$, la d\'efinition de l'ensemble $P_r$ des paires
de suites d'entiers associ\'e \`a un entier $r\geq 1$ est donn\'ee dans la d\'efinition \ref{defipr} du \S I.\ref{algederi}, et pour
$$\underline{i}=(i_1,\ldots,i_n),\ \ \  \underline{j}=(j_1,\ldots,j_m)$$
on a
$$\underline{x}_{\underline{i}}=(x_{i_1},\ldots,x_{i_n}),\ \ \ \
\underline{x}_{\underline{j}}=(x_{j_1},\ldots,x_{j_m}).$$

Nous arrivons maintenant \`a la d\'efinition et aux r\'esultats
principaux de cette section.

\begin{defi}
\label{battage}
Soit $\la\un{x}^1;\un{x }^2\ra$ un couple de suites. On appelle battage de $\la\un{x}^1;\un{x }^2\ra$ et on
note $\sh ( \un{x}^1 ,\un{x }^2 )$ l'ensemble des suites obtenues en m\'elangeant les \'el\'ements des deux suites
$\un{x}^1$ et $\un{x}^2$ en pr\'eservant l'ordre interne de chacune d'elles.
\end{defi}

{\bf Exemple}. Soit $\un{x}^1 =(a,b)$ et $\un{x}^2 =c$, alors
\begin{equation}
\sh (\un{x}^1 ,\un{x}^2 )= \{ (a,b,c),\, (a,c,b),\, (c,a,b) \}.
\end{equation}

\begin{pro}
\label{seque}
Soit $\la\un{x}^1;\un{x}^2\ra$ un couple de suites. L'ensemble des mots $\un{s}$ telles que $\la\un{x}^1;\un{x}^2\ra$
soit dans $C_{\un{s}}$ est donn\'e par l'ensemble de battage $\sh (\un{x}^1 ,\un{x}^2 )$.
\end{pro}

\begin{proof}
Pour d\'emontrer cette proposition, nous introduisons une action de $X$ sur $C(X)$ qui nous aidera \`a
classer les couples apparaissant dans $C_{\un{x}}$.\\

Soit $x\in X$, et $\underline{x}^1$, $\underline{x}^2\in X^*$. La
multiplication dans $\kK \langle\langle X\rangle\rangle\otimes_{\kK}
k\langle\langle X \rangle\rangle$ des deux \'el\'ements
$\Delta(X)$ et ${\underline{x}^1}\otimes {\underline{x}^2}$ donne
\begin{equation}
\label{deriva}
\Delta(x)\cdot ({\un{x}^1}\otimes {\un{x}^2}) = x {\un{x}^1} \otimes {\un{x}^2} + \un{x}^1 \otimes x {\un{x}^2}.
\end{equation}

Cette formule peut se traduire en une action, ou plut\^ot une somme de deux actions de $X$ sur $C(X)$. En effet,
d\'efinissons deux actions de $X$ sur $C(X)$ comme suit:
$$
\left .
\begin{array}{llll}
a^+ : & X \times C(X) & \rightarrow
& C(X) ,\\
& ({x} , \la\un{x}^1 ; \un{x}^2 \ra) & \mapsto & \la\un{x}^1 ;
(x,\un{x}^2) \ra,
\end{array}
\right .
$$
$$
\left .
\begin{array}{llll}
a^- : & X \times C(X) & \rightarrow
& C(X) ,\\
& ({x} , \la\un{x}^1 ; \un{x}^2 \ra) & \mapsto & \la(x, \un{x}^1);
\un{x}^2\ra,
\end{array}
\right .
$$
o\`u $(x,\un{x}^i)$ d\'enote la concat\'enation des suites $x$ et $\un{x}^i$ dans $X^*$.

Pour tout $x \in \Omega $, les op\'erateurs $a^-_{x}$ et $a^+_{x}$ v\'erifient
$$a^+_{x_1 } a^-_{x_2 }  = a^-_{x_2 } a^+_{x_1 } .$$
Les actions de $a^+_{x}$ et de $a^-_x$ sont {\it non
commutatives}:\\

On a $a^+_{x_1 } a^+_{x_2 } = a^+_{x_2} a^+_{x_1 }$ (resp.  $a^-_{x_1 }
a^-_{x_2 } = a^-_{x_2} a^-_{x_1 }$) si et seulement si $x_1 =x_2$.\\

L'\'egalit\'e (\ref{deriva}) se traduit par la relation r\'ecursive suivante sur les ensembles $C_{\un{x}}$ :
\begin{equation}
\di C_{\un{x}} =(x_1 ,\dots ,x_r) =\di a^{+}_{x_1} \left ( \di
C_{(x_2,\dots ,x_r)} \right ) + a^{-}_{x_1} \left ( \di
C_{(x_2,\dots ,x_r)} \right),
\end{equation}
o\`u l'on \'etend les $a^+$ et $a^-$ aux ensembles de paires en
prenant la somme.

Autrement dit, on a l'\'enonc\'e suivant qui d\'ecoule en fait
imm\'ediatement de cette remarque :

\begin{lem}
\label{sopder}
Soit $\un{x} =(x_1 ,\dots ,x_r)$. Alors les \'el\'ements de $C_{\un{x}}$ sont donn\'es par
\begin{equation}
\label{opfor}
a^{\sigma_1}_{x_1} a^{\sigma_2}_{x_2} \dots
a^{\sigma_r}_{x_r} \la\emptyset ;\emptyset \ra,
\end{equation}
o\`u $\sigma_i = {\pm 1}$, $i=1,\dots ,r$.
\end{lem}

Nous terminons maintenant la d\'emonstration de la proposition \ref{seque}.\\

Une direction est facile: si $\un{x}\in \sh (\un{x}^1,\un{x}^2)$, alors (\ref{Cw}) montre que
$\la\un{x}^1;\un{x}^2\ra\in C_{\un{x}}$.

Supposons donc que $\un{x}$ est une suite telle que $\la\un{x}^1;\un{x}^2\ra$ est dans $C_{\un{x}}$. Alors
la longueur $r$ de $\un{x}$ est \'egale \`a la somme des longueurs
de $\un{x}^1$ et de $\un{x}^2$ (disons $n$ et $m$
respectivement avec $n+m=r$). Ecrivons $\un{x}=(x_1,\ldots,x_r)$.
On a vu au lemme \ref{sopder} que les $2^r$ \'el\'ements de
$C_{\un{x}}$ sont donn\'es par
$$a_{x_1}^{\sigma_1}\cdots a_{x_r}^{\sigma_r}\la\emptyset;\emptyset\ra$$
o\`u $\sigma_i=\pm 1$.

Si $\la\un{x}^1;\un{x}^2\ra\in C_{\un{x}}$, i.e. si $\un{x}$ est un $r$-uple tel que
\begin{equation}
\label{small}
a_{x_1}^{\sigma_1}\cdots a_{x_r}^{\sigma_r}\la\emptyset;\emptyset\ra=\la\un{x}^1;\un{x}^2\ra,
\end{equation}
alors il y a $n$ des $\sigma_i$ qui sont \'egaux \`a $-1$, disons $\sigma_{i_1}=\cdots= \sigma_{i_n}=-1$, et
les $m$ autres sont \'egaux \`a $+1$, disons $\sigma_{j_1}=\cdots=\sigma_{j_m}=+1$,
o\`u $\{i_1,\ldots,i_n\}\cup \{j_1,\ldots,j_m\}=\{1,\ldots,r\}$ et $1\le i_1<\cdots i_n\le r$, $1\le j_1<\cdots <j_m\le r$.

Alors (\ref{small}) implique que $(s_{i_1},\ldots,s_{i_n})= (\omega_1^1,\ldots,\omega_n^1)$ et
$(s_{j_1},\ldots,s_{j_m})=(x_1^2,\ldots,x_m^2)$, ce qui \'equivaut \`a dire que $\un{x}\in
\sh (\un{x}^1,\un{x}^2)$.
\end{proof}

\subsection{Combinatoire sur $\eE$}

Pla\c cons-nous maintenant dans la big\`ebre $\eE$, c'est-\`a-dire dans l'alg\`ebre $\kK\langle\langle Y \rangle\rangle$
munie du coproduit $\Delta_*$ de la d\'efinition \ref{b}. Rappelons que $Y$ est un alphabet cod\'e par un semi-groupe, et que $Y^*$
d\'enote l'ensemble des suites $(y_{s_1},\ldots,y_{s_r})$ de $r$ \'el\'ements de $Y$; la suite de longueur $r=0$ est
$\emptyset$. On consid\`ere la concat\'enation $\un{y}\un{y}'$ de deux suites dans $Y^*$; il
est entendu que $\un{y}\emptyset= \emptyset \un{y}=\un{y}$; la longueur de la suite concat\'en\'ee de deux suites de
longueurs $n$ et $m$ respectivement est $n+m$.

\begin{defi}
Comme au \S \ref{combinaA}, nous associons \`a chaque $\un{y}\in Y^*$ un ensemble de couples de suites $C^*_{\un{y}}$
``engendr\'e par $\Delta_*(\un{y})$'', i.e. contenant tous les couples apparaissant dans la somme
$\Delta_*(\un{y})$.
\end{defi}

\noindent $\bullet$ Pour $r=0$, d'apr\`es la d\'efinition \ref{b},
on a
$$\Delta_*(\emptyset)=1\otimes 1=,$$
donc
\begin{equation}
C^*_\emptyset=\{\la\emptyset;\emptyset\ra\}.
\end{equation}

\noindent $\bullet$ Pour $r=1$ et $y_s\in Y$, on a
\begin{equation}
\Delta_*(y_s)=\sum_{k+l=s} y_k \otimes y_l .
\end{equation}
On a donc
\begin{equation}
C^*_{y_s} =\Bigl\{\la y_k, y_l \ra\mid k+l=s, 0\leq k,l\leq s \Bigr\} .
\end{equation}

\noindent $\bullet$ Pour $r=2$, si $y_i$ et $y_j$ sont deux lettres de $Y$ et on pose $\un{y} = y_i y_j$, on a
\begin{equation}
\label{opdif}
\left .
\begin{array}{lll}
\Delta_*(\un{y}) & = & \Delta_*(y_i)\Delta_*(y_j)\\
 & = & \biggl(\di\sum_{k+l=i} y_k \otimes y_l\biggr) \biggl(\di\sum_{k'+l'=j} y_{k'} \otimes y_{l'} \biggr)\\
 & = & \di\sum_{k+l=i} \di\sum_{k'+l'=j} y_k y_{k'} \otimes y_l y_{l'} .\\
\end{array}
\right .
\end{equation}

On obtient donc
\begin{equation}
C^*_{\un{y}} =
\Bigl\{\Bigl\la y_k y_{k'} ; y_l y_{l'} )\Bigr\ra\mid
k+k'=i ,\ l+l'=j,\ 0\leq k,k'\leq i ,\ 0\leq l,l'\leq j \Bigr\} .
\end{equation}

On arrive \`a la d\'efinition et aux r\'esultats principaux de cette section.

\begin{defi}
Soit $\la\un{y}^1 ;\un{y}^2 \ra$ un couple de suites. On appelle battage contractant de $\la\un{y}^1 ;\un{y}^2 \ra$, et
on note $\csh (\un{y}^1 ,\un{y}^2 )$, l'ensemble des suites obtenues par battage de $\la\un{y}^1 ;\un{y}^2 \ra$ suivi de
la contraction \'eventuelle
\begin{equation}
(y_{s^1_i} ,y_{s^2_j} )\stackrel{*}{\mapsto} \di y_{s^1_i +s^2_j} ,
\end{equation}
d'une ou plusieurs paires $(y_{s^1_i} ,y_{s^2_j} )$ d'\'el\'ements cons\'ecutifs provenant de $\un{y}^1$ et $\un{y}^2$.
\end{defi}

\begin{pro}
\label{sotdif}
Soit $\la\un{y}^1 ;\un{y}^2 \ra$ un \'el\'ement de $Y^* \times Y^*$. L'ensemble des suites $\un{y}$ telles que
$\la\un{y}^1 ;\un{y}^2 \ra$ soit dans $C^*_{\un{y}}$ est donn\'e par le battage contractant de $(\un{y}^1 ,\un{y}^2 )$.
\end{pro}

\begin{proof}
Soient $a^+$ et $a^-$ les op\'erateurs introduits au \S \ref{combinaA}.  Notons les faits suivants.\\

\noindent $\bullet$ Pour $y_s\in Y$, les \'el\'ements de $C^*_{y_s}$ s'obtiennent en appliquant au couple
$\la\emptyset;\emptyset\ra$, seul \'el\'ement de $C^*_\emptyset$,
tous les op\'erateurs
\begin{equation}
a_{y_{s_1} }^- a_{y_{s_2}}^+\ \ {\mbox{avec}}\ \ s_1 +s_2 =s ,\ 0\leq s_1 ,s_2 \leq s \}.
\end{equation}

\noindent $\bullet$ Pour $\un{y}=y_{s_1} \ldots y_{s_r} \in Y^*$, $r\ge 2$, les \'el\'ements de $C^*_{\un{y}}$
s'obtiennent en appliquant tous les op\'erateurs
\begin{equation}
a_{y_{k_1} }^- a_{y_{k_2}}^+\ \ {\mbox{avec}}\ \  k_1 +k_2 =s_1 ,\  0\leq k_1 ,k_2 \leq s_1 \}
\end{equation}

aux couples de $C^*_{y_{s_2} \ldots y_{s_r}}$. On r\'esume ce r\'esultat dans l'\'enonc\'e suivant.

\begin{lem}
\label{sopdif}
Soit $\un{y} =y_{s_1}\dots y_{s_r}$ une suite de $Y^*$, et $C^*_{\un{y}}$ l'ensemble des couples de
suites engendr\'e par $\Delta_*(\un{y})$. Tout \'el\'ement de $C^*_{\un{y}}$ est de la forme
\begin{equation}
b_{y_{s_1}} b_{y_{s_2}} \dots b_{y_{s_r}} \la\emptyset ;\emptyset \ra,
\end{equation}

o\`u chaque op\'erateur $b_{y_{s_i}}$ est de la forme
\begin{equation}
b_{y_{s_i}}=a^-_{y_{s_1^i}} a^+_{y_{s_2^i}}\ \ \mbox{avec}\ \ s_1^i +s_2^i =s_i ,
\ 0\leq s_1^i ,s_2^i \leq s_i \} .
\end{equation}
\end{lem}

Terminons maintenant la d\'emonstration de la proposition \ref{sotdif}.\\

Pour la premi\`ere direction, on suppose que $\un{y}=y_{s_1}\ldots y_{s_t} \in \csh (\un{y}^1 ,\un{y}^2 )$, $t\leq r$.

Ceci veut dire que $\un{y}$ est obtenu d'un \'el\'ement de $\sh (\un{y}^1,\un{y}^2)$, donn\'e par une partition de
$\{1,\ldots,r\}$ comme d'habitude en deux sous-ensembles $\{i_1,\ldots,i_n\}$ et $\{j_1,\ldots,j_m\}$ avec $n+m=r$, par
additions successives de paires d'\'el\'ements adjacents provenant de $\un{y}^1$ et $\un{y}^2$ respectivement.

En d'autres termes, on a pour $1\le i\le t$, $y_{s_i}$ est \'egal soit \`a une composante de $\un{y}^1$,
soit \`a une composante de $\un{y}^2$, soit \`a la contraction d'une composante de $\un{y}^1$ et une composante
de $\un{y}^2$.

Pour $1\le i\le t$, on pose
\begin{equation}
b_{y_{s_i}}=
\left\{
\begin{array}{ll}
a^-_{y_{s_j^1}} & {\rm si}\ y_{s_i}=y_{s_j^1}\ {\rm est\ une\ composante\ de\ }\un{y}^1\\
a^+_{y_{s_k^2}} & {\rm si}\ y_{s_i}=y_{s_k^2} {\rm \ est\ une\ composante\ de \ }\un{y}^2\\
a^-_{y_{s_j^1}} a^+_{y_{s_k^2}} & {\rm si}\ y_{s_i}=y_{s_{j^1} +s_{k^2}}.
\end{array}
\right.
\end{equation}

Alors $b_{y_{s_1}}\cdots b_{y_{s_t}}\la\emptyset;\emptyset\ra\in C^*_{\un{y}}$ par le lemme \ref{sopdif}, et
$b_{y_{s_1}}\cdots b_{y_{s_t}}\la\emptyset;\emptyset\ra=\la\un{y}^1;\un{y}^2\ra$, ce qui prouve que si
$\un{y}\in \csh (\un{y}^1,\un{y}^2)$, alors
\begin{equation}
\la\un{y}^1;\un{y}^2\ra\in C_{\un{y}}.
\end{equation}

On d\'emontre maintenant la direction inverse, \`a savoir: si $\un{y}=y_{s_1}\ldots y_{s_t}$ est tel que
$\la\un{y}^1;\un{y}^2\ra\in C_{\un{y}}$, alors $\un{y}\in \csh (\un{y}^1,\un{y}^2)$.

Par le lemme \ref{sopdif}, on sait que
\begin{equation}
\la\un{y}^1;\un{y}^2\ra= b_{y_{s_1}}\cdots b_{y_{s_t}}\la\emptyset;\emptyset\ra
=\Bigl\la(y_{s^1_1} \ldots y_{s^1_t};y_{s^2_1} \ldots y_{s^2_t} \Bigr\ra ,
\end{equation}
avec $b_{y_{s_i}}=a^-_{y_{s^1_i}} a^+_{y_{s^2_i}}$ pour $1\le i\le t$, et $s^1_i +s^2_i=s_i$.

Supposons que $\un{y}^1 =y_{k_1} \ldots y_{k_n}$ et $\un{y}^2 =y_{l_1} \ldots y_{l_m}$. On constate donc que
\begin{equation}
\un{y}^1=y_{k_1} \ldots y_{k_n} =y_{s^1_1} \ldots y_{s^1_t} ,
\end{equation}
et
\begin{equation}
\un{y}^2=y_{l_1} \ldots y_{l_m} =y_{s^2_1} \ldots y_{s^2_t} .
\end{equation}

Avec les conditions $s^1_i +s^2_i =s_i$ pour $1\le i\le t$, ceci force la valeur de chaque $y_{s^1_i}$ et $y_{s^2_i}$.
En effet, on a forc\'ement
\begin{equation}
(y_{s^1_i} ,y_{s^2_i} )=
\left\{
\begin{array}{ll}
(y_{k_j} ,\emptyset)&{\rm si}\ y_{s_i}=y_{k_j} \ {\rm est\ une\ composante\ de\ }\un{y}^1\\
(\emptyset, y_{l_j} )&{\rm si}\ y_{s_i}=y_{l_j}\  {\rm est\ une\ composante\ de\ }\un{y}^2\\
(y_{k_j},y_{l_{j'}} )&{\rm si}\ s_i =k_j +l_{j'} \ {\rm est\ une\ somme}.
\end{array}
\right.
\end{equation}
Ceci termine la d\'emonstration de la proposition \ref{sotdif}.
\end{proof}

\section{\'El\'ements primitifs de $\aA$ et $\eE$}

\subsection{\'El\'ements primitifs de $\aA$}

On se place de nouveau dans la cog\`ebre $\aA$, d'alg\`ebre sous-jacente $\kK\la\la X\ra\ra$, sur l'alphabet $X$.
Rappelons qu'un {\it moule} est naturellement associ\'e \`a une s\'erie non
commutative dans les variables $x\in X$, i.e. un \'el\'ement
$$\di\sum_{\un{x}\in X^*} M^{\un{x}} \un{x}\in \kK \la\la X\ra\ra.$$

\begin{defi}
Un moule $M^{\bullet }$ est dit {\it alternal} si
\begin{equation}
\label{alter}
\di\sum_{\un{x} \in  \sh (\un{x}^1 ,\un{x}^2 ) } \di
M^{\un{x} } =0\ \ \ \ \forall \un{x}^1 ,\un{x}^2\in X^*\setminus\{
1\}.
\end{equation}
\end{defi}

Soit $\un{x} =(x_1 ,\dots ,x_n )$ une suite, et $l(\un{x})=r$ sa
longueur. L'action de $\Delta$ sur $\un{x}$ donne une expression
de la forme
\begin{equation}
\label{leib}
\Delta(\un{x}) =\un{x} \otimes 1 +1\otimes \un{x} +
\di\sum_{(\un{x}^1 ,\un{x}^2 )\in \widetilde{C}_{\un{x}} }
\un{x}^1 \otimes \un{x}^2 ,
\end{equation}
o\`u $\di\widetilde{C}_{\un{x}} =C_{\un{x}} \setminus \{ \la\un{x}
;\emptyset \ra ,
\la\emptyset ;\un{x} \ra \}$.\\

Le but des deux th\'eor\`emes suivants est de d\'emontrer que
l'{\it alternalit\'e} d'un moule exprime sa {\it primitivit\'e} en
tant qu'\'el\'ement de la cog\`ebre $\aA$ (on rappelle que $P\in
\aA$ est {\it primitif} \ si $\Delta(P)=P\otimes 1+1\otimes P$).\\

Une propri\'et\'e essentielle des \'el\'ements primitifs est la
suivante :

\begin{lem}
Un \'el\'ement primitif $P\in \aA$ a un terme constant \'egal \`a $0$.
\end{lem}

\begin{proof}
En effet, si le terme constant de $P$ (i.e. le coefficient de $1$) est \'egal \`a $a\in \kK$, alors le terme constant
de $\Delta(P)$ est \'egal \`a $a(1\otimes 1)$ alors que le terme constant de $P\otimes 1+1\otimes P$ est \'egal \`a
$a\otimes 1+1\otimes a$. Si $P$ est primitif on a donc $a\otimes 1 +1\otimes a =a\otimes a$. On en d\'eduit
$$
\left .
\begin{array}{lll}
(a+1)\otimes (a+1) & = & a\otimes a +a\otimes 1 +1\otimes a +1\otimes 1 ,\\
  & = & (2a+1)\otimes (2a +1) ,
\end{array}
\right .
$$
soit, $a+1=2a+1$, et donc $a=0$.
\end{proof}

Ceci correspond au fait que les \'el\'ements de Lie dans une alg\`ebre associative libre appartiennent tous \`a l'id\'eal
${\cal M}$ engendr\'e par les s\'eries sans terme constant (voir la partie I).

\begin{pro}
\label{homog}
On note $X^{*,r}$ l'ensemble des mots de longueur $r$ construits sur $X$. L'\'el\'ement $\di\sum_{\un{x} \in X^{*,r}}
M^{\un{x}} \un{x}$ est un \'el\'ement primitif si et seulement si pour tout couple de suites $\la\un{x}^1 ; \un{x}^2 \ra$
avec $l(\un{x}^1)=n>0$, $l(\un{x}^2)=m>0$, $n+m=r$, on a
\begin{equation}
\label{schu}
\di \sum_{\un{x} \in \sh (\un{x}^1 ,\un{x}^2 ) } M^{\un{x} } =0 .
\end{equation}
\end{pro}

\begin{proof}
Posons
\begin{equation}
P=\di\sum_{\un{x} \in X^{*,r}} M^{\un{x}} \un{x} .
\end{equation}
En utilisant la formule (\ref{leib}), on voit que
\begin{equation}
\label{forme}
\left .
\begin{array}{lll}
\Delta(P)=\Delta \Bigl ( \di\sum_{\un{x} \in X^{*,r} }
M^{\un{x} } \un{x} \Bigr) & = & \di\sum_{\un{x}\in X^{*,r} } M^{\un{x}}\ \Delta(\un{x})\\
&=& \di\sum_{\un{x} \in X^{*,r} } M^{\un{x} } (\un{x} \otimes
1 + 1\otimes \un{x} ) +R ,
\\
& = & P\otimes 1 +1\otimes P +R,
\end{array}
\right .
\end{equation}
o\`u le reste $R$ est \'egal \`a
\begin{equation}
\label{reste}
\di\sum_{\un{x}  \in X^{*,r} } M^{\un{x} }
\Biggl( \sum_{\la\un{x}^1 ;\un{x}^2 \ra\in = \widetilde{C}_{\un{x}} }
\un{x}^1 \otimes \un{x}^2 \Biggr).
\end{equation}

Par d\'efinition, $\di\sum_{\un{x} \in X^{*,r}} M^{\un{x}} \un{x}$ est un \'el\'ement primitif si et seulement si
$R=0$.\\

Soit $\la\un{x}^1 ;\un{x}^2 \ra$ un couple de suites de (\ref{reste}). Par la proposition \ref{seque}, l'ensemble des mots
$\un{x} \in X^{*,r}$ engendrant ce couple, i.e. tel que ce couple appara\^\i t dans $C_{\un{x}}$, est obtenu par battage de
$\la\un{x}^1  ;\un{x}^2 \ra$. En regroupant dans (\ref{reste}) les termes avec le m\^eme couple $\la \un{x}^1  ;\un{x}^2 \ra$, on
r\'eecrit (\ref{reste}) comme somme d'\'el\'ements de la forme
\begin{equation}
\left ( \di\sum_{\un{x} \in \sh (\un{x}^1 ,\un{x}^2 ) } M^{\un{x} }
\right ) \un{x}^1 \otimes \un{x}^2 .
\end{equation}
Le reste dans (\ref{forme}) est donc nul si et seulement si on a la condition (\ref{schu}).
\end{proof}

\begin{thm}
\label{gener}
Un \'el\'ement $P=\di\sum_{X^*} M^{\un{x}} \un{x} \in \aA$ est un \'el\'ement primitif si et seulement si le moule
associ\'e $M^{\bullet}$ est un moule alternal.
\end{thm}

\begin{proof}
On peut \'ecrire $P$ comme somme de composantes homog\`enes $P=\di\sum_{n\le 0} P_n$.  Or, si $P_n$
est primitif pour chaque $n\le 0$, on a
\begin{equation}
\Delta(P)=\di\sum_{n\le 0} \Delta(P_n)=\di\sum_{n\le 0}(P_n\otimes 1+1\otimes P_n)= P\otimes 1
+1\otimes P_n,
\end{equation}

donc $P$ est primitif. Inversement, si $P$ est primitif, on a

\begin{equation}
\begin{array}{rcl}
\Delta(P)&=& P\otimes 1+1\otimes P=\di\sum_{n\le 0} \Delta(P_n)\\
&=&\di\sum_{n\le 0}(P_n\otimes 1+1\otimes P_n +R_n)\\
&=&P\otimes 1+1\otimes P+\di\sum_{n\le 0}R_n.
\end{array}
\end{equation}

Donc $\sum_{n\le 0} R_n=0$, et comme il s'agit d'une s\'erie formelle non commutative, chaque partie homog\`ene
$R_n=0$.  On obtient donc $\Delta(P_n)=P_n\otimes 1+1\otimes P_n$ pour tout $n\ge 0$, i.e. $P_n$ est primitif.

On peut donc raisonner composante homog\`ene par composante homog\`ene.\\

Supposons donc $P$ primitif; alors chaque $P_n$ est primitif, et par le th\'eor\`eme \ref{homog}, la partie homog\`ene
$M^\bullet_n$ du moule $M^\bullet$ est alors alternal.  Donc $M^\bullet$ est alternal.  Inversement, si $M^\bullet$ est
alternal, alors chaque $M^\bullet_n$ l'est, donc chaque $P_n$ est primitif, donc $P$ est primitif.
\end{proof}

\subsection{\'Ecriture dans l'alg\`ebre de Lie}

Rappelons que l'alg\`ebre associative $\kK \la\la X\ra\ra$ s'identifie avec l'alg\`ebre enveloppante universelle de
l'alg\`ebre de Lie libre ${\cal L}_X$ sur l'alphabet $X$; on a un homomorphisme injectif
\begin{equation}
\begin{array}{ccc}
{\cal L}_X & \longrightarrow & \kK \la\la X\ra\ra \\
x & \mapsto & x , \\
\mbox{[} x ,x' \mbox{]} & \mapsto & x\, x'-x'\, x.
\end{array}
\end{equation}
Par le th\'eor\`eme \ref{thm1}, les \'el\'ements primitifs de $\kK\la\la X\ra\ra$ sont exactement les \'el\'ements dans
l'image de ${\cal L}_X$. Les moules alternaux de $\aA$ peuvent donc \^etre vu comme des \'el\'ements de l'alg\`ebre de
Lie libre ${\cal L}_X$. Nous en donnons l'expression explicite dans le th\'eor\`eme suivant.\\

Rappelons d'abord quelques faits de la section I.\ref{alass}. Soit ${\cal M}$ l'id\'eal ${\cal M}$ de $\kK\la\la X \ra\ra$
engendr\'e par les s\'eries formelles sans terme constant. Alors nous avons un homomorphisme, que l'on d\'efinit sur les
mon\^omes et \'etend par lin\'earit\'e:
\begin{equation}
\label{homLie}
\begin{array}{ccc}
\psi:{\cal M}&\rightarrow& {\cal L}_X \\
x_1 \cdots x_r & \mapsto & {{1}\over{r}}[[\cdots[[x_1,x_2], x_3],\ldots,x_{r-1}],x_r].
\end{array}
\end{equation}
Or, l'inclusion naturelle ${\cal L}_X \subset k\la\la X\ra\ra$ donne en fait ${\cal L}_X\subset {\cal M}$.  Le
th\'eor\`eme de projection \ref{projec} (voir section I.\ref{alass}) dit que
$$\psi|_{{\cal L}_X}={\rm id}_{{\cal L}_X}.$$

\begin{thm}
\label{repre}
Soit $M^{\bullet }$ un moule alternal et $\un{x}$ une suite de longueur $r>0$. Soit $\sigma(\un{x})$ l'ensemble
des suites de longueur $r$ d\'eduites de $\un{x}$ par permutation des composantes.

Alors l'\'el\'ement $\di\sum_{\un{u} \in \sigma (\un{x} )} M^{\un{u} } \,\un{u}$ appartient \`a ${\cal L}_{\cal X}
\subset{\cal M}\subset \kK\la\la X \ra\ra$, et s'\'ecrit
\begin{equation}
\di {{1}\over{r}} \di\sum_{\un{u} \in \sigma (\un{x})} M^{\un{u} } \, [\un{u}]
\in {\cal L}_{\cal X},
\end{equation}
o\`u
\begin{equation}
[x_1 \cdots x_r] = [[\cdots[[x_1,x_2],x_3],\ldots,x_{r-1}], x_r] .
\end{equation}
\end{thm}

\begin{proof}
Comme $M^{\bullet }$ est alternal, il est primitif, c'est-\`a-dire $\di\sum_{\un{u} \in \sigma (\un{x} )}
M^{\un{u} } \un{u} \in {\cal L}_X$, et le r\'esultat d\'ecoule alors imm\'ediatement de (\ref{homLie}) et du th\'eor\`eme
de projection.
\end{proof}

\begin{rema}
Dans les applications du calcul moulien, ce n'est pas ce
th\'eor\`eme que nous utiliserons, mais sa contrepartie dans le
cas o\`u l'alg\`ebre sous-jacente est li\'ee. Par exemple, soit
$\kK \la \la {\rm Der}(A) \ra\ra$ la big\`ebre introduit au \S
I.\ref{exemples}. Un \'el\'ement $P\in \kK \la \la {\rm Der}(A) \ra\ra$
s'\'ecrit
\begin{equation}
P=\di\sum_{\un{a}} P^{\un{a}} D_{\un{a}},
\end{equation}
o\`u $D_{\un{a}} =D_{a_n} \circ D_{a_{n-1}} \circ \dots \circ
D_{a_1}$.

Pour savoir si $P$ est primitif pour le coproduit $\Delta$ de $\kK \la \la {\rm Der}(A) \ra\ra$, on associe \`a $P$ son
{\it ``repr\'esentant" libre}, not\'ee $P_l$, d\'efini par
\begin{equation}
P_l =\di\sum_{\un{a}} P^{\un{a}} \un{a} \in \kK \la\la X\ra\ra .
\end{equation}
Le th\'eor\`eme \ref{gener} s'applique et on a: {\it si le moule $M^{\bullet}$ est alternal, alors l'\'el\'ement
$P=\sum M^\bullet X_\bullet$ est primitif}.\\

Evidemment, on ne capte pas de cette fa\c{c}on tous les \'el\'ements primitifs.

Consid\'erons les s\'eries formelles non commutatives construites sur l'alphabet \`a trois lettres $\{ D_1
, D_2 ,D_3 \}$, o\`u $D_1$, $D_2$ et $D_3$ sont trois d\'erivations sur une alg\`ebre $A$ telles que
\begin{equation}
D_2 =D_1 +D_2 .
\end{equation}
Alors, l'\'el\'ement
\begin{equation}
D_2 D_1 -D_1 D_1 -D_1 D_3
\end{equation}
est primitif, sans que le moule associ\'e soit alternal.\\

En effet, cet \'el\'ement est associ\'e au moule $M^{\bullet}$ d\'efini par
\begin{equation}
M^{2,1} =1,\ M^{1,1}=-1,\ M^{1,3} =-1\ \mbox{\rm et}\ M^{\bullet} =0\ \mbox{\rm sinon}.
\end{equation}
Ce moule n'est pas alternal car on a
\begin{equation}
M^{2,1}+M^{1,2} =1.
\end{equation}
Pourtant, comme
\begin{equation}
D_2 D_1 -D_1 D_1 -D_1 D_3 =D_3 D_1 -D_1 D_3,
\end{equation}
c'est bien un \'el\'ement primitif.
\end{rema}

Jean \'Ecalle utilise constamment ce va-et-vient entre les objets construits sur des alg\`ebres li\'ees et
leurs repr\'esentations libres.

\subsection{\'El\'ements primitifs de $\eE$}

\begin{defi}
Un moule $M^{\bullet }$ est dit alternel si
\begin{equation}
\di\sum_{\un{y} \in \csh (\un{y}^1 ,\un{y}^2 )} M^{\un{y}} =0\ \ \ \
\forall \un{y}^1 ,\un{y}^2 \in Y^*\setminus\{\emptyset\}.
\end{equation}
\end{defi}

Pla\c cons-nous maintenant dans la cog\`ebre $\eE$, d'alg\`ebre sous-jacente \'egale \`a $\kK\la\la Y\ra\ra$ sur
l'alphabet libre $Y$ index\'e par le semigroupe $\nN$. On a alors le th\'eor\`eme suivant, analogue du th\'eor\`eme
\ref{repre}:

\begin{thm}
\label{thm7}
Un \'el\'ement $\di\sum_{\un{y}\in Y^*} M^{\un{y}} \un{y} \in \eE$ est primitif si et seulement si $M^{\bullet }$ est
un moule alternel.
\end{thm}

\begin{proof}
On note $P=\di\sum_{\un{y} \in Y^*} M^{\un{y}} \un{y}$. $P$ primitif signifie que
\begin{equation}
\Delta_*(P) = P\otimes 1+1\otimes P.
\end{equation}
Or, on a
\begin{equation}
\label{long}
\begin{array}{lcl}
\Delta_*(P) & = & \di\sum_{\un{y}\in Y^*} M^{\un{y}}\ \Delta_*(\un{y})\\
 & = & \di\sum_{y_{s_1} \ldots y_{s_r} \in Y^*}M^{y_{s_1} \ldots y_{s_r} }
\ \Delta_*(y_{s_1} \ldots y_{s_r} )\\
 & = & \di\sum_{y_{s_1} \ldots y_{s_r} \in Y^*}M^{y_{s_1} \ldots y_{s_r} }
\ \Delta_*(y_{s_1}) \ldots \Delta_* (y_{s_r} )\\
 & = & \di\sum_{y_{s_1} \ldots y_{s_r} \in Y^*}M^{y_{s_1} \ldots y_{s_r} }
\ \left ( \di\sum_{k_1 +l_1 =s_1} y_{k_1} \otimes y_{l_1} \right )  \ldots
\left ( \di\sum_{k_r +l_r =s_r} y_{k_r} \otimes y_{l_r} \right ) \\
 & = & \di\sum_{y_{s_1} \ldots y_{s_r} \in Y^*}M^{y_{s_1} \ldots y_{s_r} }
\ \di\sum_{k_i +l_i =s_i , \ 1\leq i\leq r} y_{k_1} \ldots y_{k_r} \otimes y_{l_1} \ldots y_{l_r} .
\end{array}
\end{equation}

Ici, par la d\'efinition de $\Delta_*$, les $y_{k_i}$ et $y_{l_i}$ appartiennent \`a $Y\cup\{\emptyset\}$,
c'est-\`a-dire que la somme porte sur l'ensemble des $r$-uples de couples
$\Bigl((y_{k_1},y_{l_1}),\ldots,(y_{k_r},y_{l_r})\Bigr)$ avec $k_i +l_i =s_i$ restants.

Soit $A$ l'ensemble de ces $r$-uples de couples, moins les deux \'el\'ements suivants:
\begin{equation}
\Bigl((y_{k_1} , y_{l_1} ),\ldots,( y_{k_r} ,y_{l_r} )\Bigr)
=\Bigl((y_{s_1} ,\emptyset),\ldots,(y_{s_r} ,\emptyset)\Bigr)
\end{equation}
et
\begin{equation}
\Bigl((y_{k_1} ,y_{l_1} ),\ldots,(y_{k_r},y_{l_r} )\Bigr)
=\Bigl((\emptyset,y_{s_1}),\ldots,(\emptyset,y_{s_r} )\Bigr).
\end{equation}

L'expression de $\Delta_*(P)$ donn\'ee dans la derni\`ere ligne de (\ref{long}) devient alors
\begin{equation}
\begin{array}{lcl}
\Delta_*(P) & = & \di\sum_{\un{y} \in Y^*} M^{\un{y}} \un{y} \otimes 1+
\di\sum_{\un{y} \in Y^*} M^{\un{y}} 1\otimes \un{y} \\
 & & + \di\sum_{\un{y} \in Y^*} \sum_{A} y_{k_1} \ldots y_{k_r} \otimes y_{l_1} \ldots y_{l_r} , \\
 & = & P\otimes 1+1\otimes P+ R ,
\end{array}
\end{equation}
o\`u
\begin{equation}
R=\di\sum_{\un{y}\in Y^*} M^{\un{y}} \sum_{A} y_{k_1} \ldots y_{k_r} \otimes y_{l_1} \ldots y_{l_r} .
\end{equation}
$P$ est donc primitif si et seulement si $R=0$.\\

Posons $\un{y}^1=y_{k_1} \dots y_{k_r} $ et $\un{y}^2=y_{l_1} \dots y_{l_r}$; notons que, contrairement aux apparences,
les longueurs de $\un{y}^1$ et $\un{y}^2$ ne sont pas forc\'ement \'egales \`a $r$ puisque l'on ignore les composantes
\'egales \`a $\emptyset$. La proposition \ref{sotdif} dit qu'un couple $\la\un{y}^1;\un{y}^2\ra$ appartient \`a
$C^*_{\un{y}}$ si et seulement si $\un{y}\in \csh (\un{y}^1,\un{y}^2)$.

On voit donc que le coefficient de chaque expression $\un{y}^1 \otimes \un{y}^2$ dans la somme $R$ est \'egal \`a
\begin{equation}
\di\sum_{\un{y}\in \csh (\un{y}^1,\un{y}^2 )} M^{\un{y}}.
\end{equation}
$R$ s'annule donc si et seulement si $M^\bullet$ est alternel.
\end{proof}

\section{\'El\'ements ``group-like" de $\aA$ et $\eE$}

\subsection{\'El\'ements ``group-like" de $\aA$}

Rappelons de la d\'efinition \ref{defiprimigroup} qu'un \'el\'ement $P\in \kK\la\la X\ra\ra$ est dit ``group-like'' pour
le coproduit $\Delta:\kK\la\la X\ra\ra \rightarrow \kK\la\la X\ra\ra\otimes_{\kK} \kK\la\la X\ra\ra$ s'il v\'erifie
\begin{equation}
\label{aut}
\Delta(P)=P \otimes P.
\end{equation}

Deux observations permettent de pr\'eciser la nature des \'el\'ements group-like de $\eE$.\\

- Aucun polyn\^ome ne peut v\'erifier (\ref{aut}). En effet, d\'efinissons la ``longueur'' d'un produit tensoriel
$P_1\otimes P_2$ de deux mon\^omes comme \'etant la somme $n+m$ o\`u $n$ est la longueur de $P_1$ en tant que mon\^ome
et $m$ la longueur de $P_2$. Soit $P$ un polyn\^ome de $\kK\la X \ra$, et soit $M$ le mon\^ome le plus long apparaissant
dans $P$.  Alors on voit que la longueur de chaque terme apparaissant dans $\Delta(P)$ est inf\'erieure ou \'egale \`a
la longueur de $M$, alors que le terme $M\otimes M$, deux fois trop long, appara\^\i t dans $P\otimes P$.
Ceci montre que si un \'el\'ement de $\kK\la\la X\ra\ra$ a une chance d'\^etre group-like, il doit
s'agir d'une s\'erie formelle.\\

- De m\^eme, $\Delta(P)$ fait apparaitre des couples $(1,\un{x})$
ou $(\un{x}, 1)$, ce qui n'est possible dans $P\otimes P$ que si
on fait intervenir la suite vide.  La s\'erie $P$ doit donc avoir
un terme constant $a\ne 0$. La condition $\Delta(P)=P\otimes P$
implique alors que ce terme constant v\'erifie
$a(1\otimes 1)=a\otimes a$, donc $a=1$.\\

Ces deux observations nous conduisent \`a consid\'erer les moules de la forme
\begin{equation}
P=\sum_{\un{x}\in X^*} M^{\un{x}} \un{x} ,
\end{equation}
avec $M^{\emptyset}=1$.\\

\begin{rema}
Si les d\'erivations proviennent d'un champ de vecteurs, le groupe
des automorphismes ainsi d\'efini est isomorphe au groupe des
diff\'eomorphismes tangents \`a l'identit\'e.
\end{rema}

\begin{defi}
Un moule $M^{\bullet }$ est dit {\it sym\'etral} si
\begin{equation}
\label{syme} \di\sum_{\un{x} \in \sh (\un{x}^1 , \un{x}^2 )} \di
M^{\un{x}} = M^{\un{x}^1 } M^{\un{x}^2 }\ \ \ \ \forall \un{x}^1
,\un{x}^2 \in X^* .
\end{equation}
\end{defi}

\begin{thm}
\label{thm8}
Un \'el\'ement $\di\sum_{\un{x}} M^{\un{x}} \un{x} \in \aA$ est ``group-like" si et seulement si $M^{\bullet}$ est un
moule sym\'etral.
\end{thm}

\begin{proof}
Soit $P=1+Q=1+\di\sum_{\un{\omega}\ne\emptyset} M^{\un{x}}
\un{x}$. On a
\begin{equation}
\label{egali}
\left .
\begin{array}{lll}
\Delta(P)& = & 1\otimes 1 +\di\sum_{\un{x} \in X^* } M^{\un{x}} \Delta(\un{x})  \\
& = & 1\otimes 1 + Q\otimes 1 + 1\otimes Q +\di\sum_{\un{x}\in X^{*,r}} M^{\un{x}} \ \ \cdot\!\! \di\sum_{\la
\un{x}^1 ;\un{x}^2 \ra \in \widetilde C_{\un{x}} } \un{x}^1 \otimes \un{x}^2  ,
\end{array}
\right .
\end{equation}
o\`u l'on rappelle que $C_{\un{x}}$ d\'enote l'ensemble de
couples $\la\un{x}^1; \un{x}^2\ra$ apparaissant dans la
somme $\Delta(\un{x})$, et $\widetilde C_{\un{x}}$
d\'enote le sous-ensemble des couples avec
$\un{x}^1\ne\emptyset$, $\un{x}^2\ne\emptyset$.

On a aussi
$$P\otimes P =1\otimes 1 +Q\otimes 1+ 1\otimes Q +Q\otimes Q,$$
o\`u
\begin{equation}
\label{direct}
Q\otimes Q=
\sum_{\un{x}^1\ne\emptyset,\un{x}^2\ne\emptyset}
M^{\un{x}^1} M^{\un{x}^2} \un{x}^1 \otimes
\un{x}^2  .
\end{equation}
Pour que $P$ soit group-like, il faut donc que le reste de (\ref{egali}) soit \'egal \`a
(\ref{direct}).

Soit $\la\un{x}^1 ;\un{x}^2 \ra$ un couple de suites intervenant dans (\ref{direct}). Rappelons de la proposition
\ref{seque}, \S 3.1, que ce couple $\la\un{x}^1;\un{x}^2\ra$ appartient \`a $C_{\un{x}}$ si et seulement si
$\un{x}$ appartient \`a $\sh (\un{x}^1,\un{x}^2)$. Pour un couple donn\'e $\la\un{x}^1;\un{x}^2\ra$, le coefficient du
terme $\un{x}^1\otimes \un{x}^2$ dans le reste de (\ref{egali}) est donc donn\'e par
$$\di\sum_{\un{x}\in \sh(\un{x}^1,\un{x}^2)} M^{\un{x}}.$$
Le reste de (\ref{egali}) est donc \'egal \`a (\ref{direct}), i.e. $P$ est group-like, si et seulement si
$$\di\sum_{\un{x}\in \sh(\un{x}^1 ,\un{x}^2 )} M^{\un{x}} = M^{\un{x}^1 } M^{\un{x}^2} ,$$
d'o\`u le th\'eor\`eme.
\end{proof}

\subsection{\'El\'ements ``group-like" de $\eE$}

La condition de sym\'etrie d'un moule correspondant au fait d'\^etre group-like dans $\eE$ est donc :

\begin{defi}
Un moule $M^{\bullet}$ est dit sym\'etrel si
\begin{equation}
\di\sum_{\un{y} \in \csh (\un{y}^1 ,\un{y}^2 )}
M^{\un{y}} =M^{\un{y}^1} M^{\un{y}^2}\ \ \ \
\forall \un{y}^1 ,\un{y}^2 \in Y^*.
\end{equation}
\end{defi}

\begin{thm}
\label{thm9}
Un \'el\'ement $\di\sum_{\un{y} \in Y^*} M^{\un{y}} {\un{y}} \in \eE$ est ``group-like" si et seulement si
$M^{\bullet }$ est un moule sym\'etrel.
\end{thm}

\begin{proof}
Elle est analogue au cas sym\'etral.  En effet, il suffit d'adapter la d\'emonstration du th\'eor\`eme \ref{thm8} en
rempla\c cant (\ref{egali}) par

\begin{equation}
\label{egali2}
\left .
\begin{array}{lll}
\Delta_*(P)& = & 1\otimes 1 +\di\sum_{\un{y} \in Y^*} M^{\un{y}} \Delta_*( \un{y})  \\
 & = & 1\otimes 1 + Q\otimes 1 + 1\otimes Q +\di\sum_{\un{y}\in Y^*}
M^{\un{y}} \ \ \cdot \di\sum \un{y}^1 \otimes \un{y}^2  ,
\end{array}
\right .
\end{equation}

o\`u la deuxi\`eme somme porte sur tous les couples $\langle\un{y}^1 ;\un{y}^2 \ra$ tels que

-- $\un{y}^1=y_{k_1} \ldots y_{k_r}$, $\un{y}^2= y_{l_1} \ldots y_{l_r}$ avec
\'eventuellement certains $y_{k_i}$ ou $y_{l_i}$ \'egaux \`a $\emptyset$;

-- $k_i +l_i =s_i$ pour $1\le i\le r$ si $\un{y} =y_{s_1} \dots y_{s_r}$

-- $\un{y}^1 \ne\emptyset$ et $\un{y}^2 \ne\emptyset$.

Par la proposition \ref{sotdif}, le coefficient d'un terme donn\'e $\un{y}^1 \otimes \un{y}^2$ est
donn\'e par
\begin{equation}
\di\sum_{\un{y}\in \csh (\un{y}^1 ,\un{y}^2 )} M^{\un{y}}.
\end{equation}
Comparant donc le reste de (\ref{egali2}) avec (\ref{direct}), on voit que $P$ est group-like si et seulement si le moule
associ\'e $M^\bullet$ est sym\'etrel.
\end{proof}

\section{Exemples de moules alterna(e)l, sym\'etra(e)l}

\subsection{Un moule alternal}

Dans cet exemple, nous prenons pour $\Omega$ un ensemble
d\'enombrable d'ind\'etermin\'ees, et pour le corps $\kK$ le corps
$\qQ(\Omega)$ des fractions rationnelles dans les \'el\'ements de
$\Omega$.\\

On d\'efinit le moule \'el\'ementaire $T^{\bullet}$ par
\begin{equation}
\left .
\begin{array}{l}
\di T^{\emptyset} =0 ,\ \ \ T^{\omega} =0\ \ \forall \omega \in\Omega, \\
\di T^{(\omega_1,\ldots,\omega_r)} =\di {1\over (\omega_2
-\omega_1 )} ..{1\over (\omega_3 -\omega_2 )} \dots .. {1\over
(\omega_r -\omega_{r-1} )} ,\ r\geq 2 .
\end{array}
\right .
\end{equation}
On a

\begin{lem}
Le moule $\di T^{\bullet}$ est alternal.
\end{lem}

\begin{proof}
Elle se fait par r\'ecurrence sur la longueur des suites. On doit
v\'erifier pour toutes suites $\un{\omega}^1\ne\emptyset$,
$\un{\omega}^2\ne\emptyset$, la propri\'et\'e
$$\di\sum_{\un{\omega} \in \sh (\un{\omega}^1 ,\un{\omega}^2 )}
T^{\un{\omega}} =0.$$
La propri\'et\'e est trivialement vraie si $l(\un{\omega} )=2$.\\

Soient $\un{\omega}^1 =(\omega_1^1 ,\dots ,\omega_n^1 )$ et
$\un{\omega}^2 =(\omega_1^2 ,\dots ,\omega_m^2 )$ telles que
$n+m=r>2$. Supposons que la propri\'et\'e d'alternalit\'e soit
vraie pour toutes les suites de longueur $\leq
r-1$.\\

On commence par noter que
\begin{equation}
T^{(\omega_1,\ldots,\omega_r)} =\di {1\over (\omega_r
-\omega_{r-1} )} T^{(\omega_1,\ldots,\omega_{r-1})}.
\end{equation}
On a de plus la propri\'et\'e suivante du battage de deux suites :

\begin{lem}
\label{petithorreur}
L'ensemble $sh(\un{\omega}^1,\un{\omega}^2)$
est la r\'eunion disjointe des quatre ensembles suivants:
$$\Bigl(\sh (\un{\omega}^1
,\di\un{\omega}^{2}_{\le m-2 } ),\omega_{m-1}^2, \omega_m^2\Bigr)
\coprod \Bigl(\sh (\un{\omega}^1_{\le n-1 } ,\un{\omega}^2_{\le m-1}
),\omega_n^1 , \omega_m^2\Bigr) \coprod$$
$$\Bigl(\sh (\un{\omega}^1_{\le n-1} ,\un{\omega}^2_{\le m-1} ),\omega_m^2 ,
\omega_n^1\Bigr) \coprod \Bigl(\sh (\di\un{\omega}^1_{\le n-2}
,\un{\omega}^2 ),\omega_{n-1}^1 , \omega_n^1\Bigr) ,$$ o\`u
$\un{\omega}_{\le j}$ d\'enote la sous-suite
$(\omega_1,\ldots,\omega_j)$ des $j$ premi\`eres composantes de
$\un{\omega}$.
\end{lem}

\begin{proof}
Soit $\un{\omega}\in \sh (\un{\omega}^1,\un{\omega}^2)$.  Alors
$l(\un{\omega})=r$, et on peut se demander que peuvent \^etre les
deux derni\`eres composantes de $\omega$.  Comme le battage ne
m\'elange pas l'ordre interne des composantes de $\un{\omega}^1$
et de $\un{\omega}^2$, on voit que les deux derni\`eres
composantes de $\un{\omega}$, en l'ordre, doivent former l'un des
couples suivants: $(\omega_{n-1}^1,\omega_n^1)$,
$(\omega_{m-1}^2,\omega_m^2)$, $(\omega_n^1,\omega_m^2)$,
$(\omega_m^2,\omega_n^1)$.  Les $r-2$ premi\`eres composantes de
$\un{\omega}$ sont donc forc\'ement obtenues par battage des
composantes restantes de $\un{\omega}^1$ et $\un{\omega}^2$, ce
qui d\'emontre le r\'esultat.
\end{proof}

On a donc
\begin{equation}
\label{37} \left .
\begin{array}{lll}
\di\sum_{\un{\omega} \in \sh (\un{\omega}^1 ,\un{\omega}^2 ) }
T^{\un{\omega}}  & = & \di {1\over (\omega_m^2 -\omega_{m-1}^2 )}
\di\sum_{\un{s} \in \bigl(\sh (\un{\omega}^1 ,\un{\omega}^{2}_{\le m-2 }),\omega_{m-1}^2\bigr)} T^{\un{s}}\\
&+& \di {1\over (\omega_m^2 -\omega_n^1 )} \di\sum_{\un{s} \in
\bigl(\sh (\un{\omega}^1_{\le n-1 } ,\un{\omega}^2_{\le
m-1}),\omega_n^1\bigr)}
T^{\un{s}} \\
&+& \di {1\over (\omega_n^1 -\omega_m^2 ) } \di\sum_{\un{s} \in
\bigl(\sh (\un{\omega}^1_{\le n-1 } , \un{\omega}^2_{\le m-1}),\omega_m^2\bigr)} T^{\un{s}}\\
&+&  \di {1\over (\omega_n^1 -\omega_{n-1}^1 )} \di\sum_{\un{s}
\in \bigl(\sh (\un{\omega}^1_{\le n-2}
,\un{\omega}^2),\omega_{n-1}^1\bigr) } T^{\un{s}} .
\end{array}
\right .
\end{equation}
Par ailleurs, le moule $T^{\bullet }$ \'etant alternal jusqu'\`a
l'ordre $r-1$, on a les \'egalit\'es :
\begin{equation}
\left .
\begin{array}{l}
\di\sum_{\un{s} \in \bigl(\sh (\un{\omega}^1 ,\un{\omega}^{2}_{\le
m-2 }),\omega_{m-1}^2\bigr)} T^{\un{s}} + \di\sum_{\un{s} \in
\bigl(\sh(\un{\omega}^1_{\le n-1 } ,\un{\omega}^2_{\le
m-1}),\omega_n^1\bigr)}
T^{\un{s}} =0 , \\
\di\sum_{\un{s} \in \bigl(\sh (\un{\omega}^1_{\le n-1 } ,
\un{\omega}^2_{\le m-1}),\omega_m^2\bigr)} T^{\un{s}} +
\di\sum_{\un{s} \in \bigl(\sh (\un{\omega}^1_{\le n-2}
,\un{\omega}^2),\omega_{n-1}^1\bigr) } T^{\un{s}} =0 .
\end{array}
\right .
\end{equation}
L'\'equation (\ref{37}) peut donc s'\'ecrire sous la forme
\begin{equation}
\label{37bis} \left .
\begin{array}{l}
\di\sum_{\un{\omega} \in \sh (\un{\omega}^1 ,\un{\omega}^2 ) }
T^{\un{\omega}} = \left [ \di {1\over \omega_m^2 -\omega_n^1 }
-\di {1\over \omega_m^2 -\omega_{m-1}^2} \right ] \di\sum_{\un{s}
\in \bigl(\sh(\un{\omega}^1_{\le n-1 } ,\un{\omega}^2_{\le
m-1}),\omega_n^1\bigr)}
T^{\un{s}} \\
+ \left [ \di {1\over \omega_n^1 -\omega_m^2} -\di {1\over
\omega_n^1 -\omega_{n-1}^1} \right ] \di\sum_{\un{s} \in \bigl(\sh
(\un{\omega}^1_{\le n-1 } , \un{\omega}^2_{\le
m-1}),\omega_m^2\bigr)} T^{\un{s}} .
\end{array}
\right .
\end{equation}
On d\'ecompose les deux sommes sous la forme
\begin{equation}
\left .
\begin{array}{lll}
\di\sum_{\un{s} \in \bigl(\sh(\un{\omega}^1_{\le n-1 }
,\un{\omega}^2_{\le m-1}),\omega_n^1\bigr)} T^{\un{s}} & = &
\di\sum_{\un{s} \in \bigl(\sh(\un{\omega}^1_{\le n-2 }
,\un{\omega}^2_{\le m-1}),\omega_{n-1}^1 \omega_n^1\bigr)}
T^{\un{s}} + \di\sum_{\un{s} \in \bigl(\sh(\un{\omega}^1_{\le n-1
} ,\un{\omega}^2_{\le m-2}),\omega_{m-1}^2 \omega_n^1\bigr)}
T^{\un{s}} , \\
\di\sum_{\un{s} \in \bigl(\sh (\un{\omega}^1_{\le n-1 } ,
\un{\omega}^2_{\le m-1}),\omega_m^2\bigr)} T^{\un{s}} & = &
\di\sum_{\un{s} \in \bigl(\sh (\un{\omega}^1_{\le n-2 } ,
\un{\omega}^2_{\le m-1}),\omega_{n-1}^1 \omega_m^2\bigr)}
T^{\un{s}} + \di\sum_{\un{s} \in \bigl(\sh (\un{\omega}^1_{\le n-1
} , \un{\omega}^2_{\le m-2}),\omega_{m-1}^2 \omega_m^2\bigr)}
T^{\un{s}} ..
\end{array}
\right .
\end{equation}
On note
\begin{equation}
A=\di\sum_{\un{s} \in \bigl(\sh(\un{\omega}^1_{\le n-2 }
,\un{\omega}^2_{\le m-1}),\omega_{n-1}^1 \bigr)} T^{\un{s}}, \ \
B=\di\sum_{\un{s} \in \bigl(\sh(\un{\omega}^1_{\le n-1 }
,\un{\omega}^2_{\le m-2}),\omega_{m-1}^2 \bigr)} T^{\un{s}} .
\end{equation}
L'\'equation (\ref{37bis}) s'\'ecrit donc
\begin{equation}
\label{37ter} \left .
\begin{array}{l}
\di\sum_{\un{\omega} \in \sh (\un{\omega}^1 ,\un{\omega}^2 ) }
T^{\un{\omega}} =A \left [ \di {1\over (\omega_{m-1}^2 -\omega_m^2
)(\omega_{n-1}^1 -\omega_m^2 ) } +\di
{1\over (\omega_{n-1}^1 -\omega_n^1 ) (\omega_{n-1}^1 -\omega_n^1 )} \right ] \\
+B \left [ \di {1\over (\omega_{m-1}^2 -\omega_m^2
)(\omega_{m-1}^2 -\omega_m^2 )}  + \di {1\over  (\omega_{n-1}^1
-\omega_n^1 ) (\omega_{m-1}^2 -\omega_n^1 ) } \right ] .
\end{array}
\right .
\end{equation}
Par hypoth\`ese d'alternalit\'e de $T^{\bullet}$ jusqu'\`a l'ordre
$r-1$, on a
\begin{equation}
A+B=0 ,
\end{equation}
soit, en notant $C_A (\un{\omega }^1 ,\un{\omega}^2 )$ et $C_B
(\un{\omega }^1 , \un{\omega}^2 )$ les coefficients de $A$ et $B$
dans (\ref{37ter} ),
\begin{equation}
\di\sum_{\un{\omega} \in \sh (\un{\omega}^1 ,\un{\omega}^2 ) }
T^{\un{\omega}} = A (C_A (\un{\omega }^1 ,\un{\omega}^2 ) -C_B
(\un{\omega}^1 ,\un{\omega}^2 )).
\end{equation}
Un simple calcul donne $C_A (\un{\omega}^1 ,\un{\omega}^2 ) -C_B
(\un{\omega}^1 ,\un{\omega}^2 ) =0$, d'o\`u le r\'esultat.
\end{proof}

\subsection{Un moule alternel}

Le moule alternel le plus simple est d\'efini par $J^{\emptyset }
=0$ et pour toute suite $\un{\omega }=(\omega_1 ,\dots ,\omega_r
)$ par
$$J^{\un{\omega }} =\di {(-1)^{r+1} \over r} .$$

\begin{lem}
\label{propbatco} Le moule $J^{\bullet }$ est alternel.
\end{lem}

\begin{proof}
Elle se fait par r\'ecurrence sur la longueur des suites. Pour
$\un{\omega } =(\omega_1 , \omega_2 )$, on a
$$J^{\omega_1 ,\omega_2 } +J^{\omega_2 ,\omega_1 } +J^{\omega_1
+\omega_2 } =-\di{1\over 2} -{1\over 2} +1 =0.$$ La propri\'et\'e
d'altern\'elit\'e est donc v\'erifi\'ee pour les suites de
longueur 2. Pour les suites de longueur $\geq 3$, on note la
propri\'et\'e suivante du battage contractant :

\begin{lem}
Pour toutes suites $\un{\omega}^1$ et $\un{\omega}^2$ avec
$l(\un{\omega}^1 )=n>0$ et $l(\un{\omega}^2 )=m>0$, l'ensemble
$csh(\un{\omega}^1,\un{\omega}^2)$ est la r\'eunion disjointe des
trois ensembles suivants:
$$\bigl(csh (\un{\omega}^1_{n-1},\un{\omega}^2),\omega^1_n\bigl)
\coprod \bigl(csh (\un{\omega}^1,\un{\omega}^2_{m-1} ),\omega^2_m
\bigr) \coprod \bigl(csh(\un{\omega}^1_{\le n-1}
,\un{\omega}^2_{\le m-1} ),(\omega^1_n +\omega^2_m)\bigr).$$
\end{lem}

\begin{proof}
Comme pour le lemme \ref {petithorreur}, il suffit de constater
que toute suite appartenant \`a $csh(\un{\omega}^1,\un{\omega}^2)$
a comme derni\`ere composante soit $\omega_n^1$, soit
$\omega_m^2$, soit la somme des deux.
\end{proof}

On a aussi
$$\di\sum_{\un{s}\in csh(\un{\omega}^1 ,\un{\omega}^2 )} J^{\un{s}}
=-\di {(r-1)\over r} \di\sum_{\un{s}\in csh(\un{\omega}^1
,\un{\omega}^2 )} J^{\un{s}_{\le r-1}}.$$ On d\'eduit donc du
lemme \ref{propbatco} l'\'egalit\'e
$$
\left .
\begin{array}{lll}
\di\sum_{\un{s}\in csh(\un{\omega}^1 ,\un{\omega}^2 )} J^{\un{s}}
& = & -\di {(r-1)\over r} \left ( \di\sum_{\un{s}\in
csh(\un{\omega}^1 ,\un{\omega}^{2,<n-1} )} J^{\bullet } +
\di\sum_{\un{s}\in csh(\un{\omega}^{1,<m-1} ,\un{\omega}^2 )}
J^{\bullet } \right ) \\
& & + \di {(r-2) \over r-1} \di\sum_{\un{s}\in
csh(\un{\omega}^{1,<m-1} ,\un{\omega}^{2,<n-1} )} J^{\bullet } .
\end{array}
\right .
$$
Par hypoth\`ese de r\'ecurrence, ces trois sommes sont nulles,
d'o\`u le r\'esultat.
\end{proof}

\subsection{Un moule sym\'etral}

On d\'efinit le moule $S^{\bullet}$ par $\di S^{\emptyset } =1$ et
\begin{equation}
S^{\omega_1 \dots \omega_r} =\di {(-1)^r \over \omega_1 (\omega_1
+\omega_2 )\dots (\omega_1 +\dots +\omega_r )} .
\end{equation}
On a

\begin{lem}
Le moule $S^{\bullet }$ est sym\'etral.
\end{lem}

\begin{proof}
Elle se fait par r\'ecurrence sur la longueur des suites. On doit
v\'erifier pour toutes suites $\un{\omega}^1$, $\un{\omega}^2$, la
propri\'et\'e
$$\di\sum_{\un{\omega} \in sh(\un{\omega}^1 ,\un{\omega}^2 )}
S^{\un{\omega}} = S^{\un{\omega}^1}  S^{\un{\omega}^2 } .$$ La
propri\'et\'e est vraie si $l(\un{\omega} )=2$. En effet, on a
\begin{equation}
\left .
\begin{array}{lll}
\di\sum_{\un{\omega}\in sh(\omega_1 ,\omega_2 )} S^{\un{\omega}} &
= &
S^{\omega_1 \omega_2} +S^{\omega_2 \omega_1 } , \\
& = & \di {1\over \omega_1 (\omega_1 +\omega_2 )} +\di {1\over
\omega_2 (\omega_1 +\omega_2 )} , \\
& = & \di {1\over \omega_1 \omega_2 } =S^{\omega_1 } S^{\omega_2 }
.
\end{array}
\right .
\end{equation}
Soient $\un{\omega}^1 =(\omega_1^1 ,\dots ,\omega_n^1 )$ et
$\un{\omega}^2 =(\omega_1^2 ,\dots ,\omega_m^2 )$ deux suites
telles que $l(\un{\omega}^1 )+l(\un{\omega}^2 )=r$, $r>2$.
Supposons que la propri\'et\'e de sym\'etralit\'e soit vraie pour
toutes les suites de longueur
$\leq r-1$.\\

On commence par noter que
\begin{equation}
\label{relS} S^{\un{\omega} =\omega_1 \dots \omega_r} =\di {1\over
\parallel \un{\omega } \parallel } S^{\un{\omega}_{<r-1} } ,
\end{equation}
o\`u $\parallel \un{\omega}\parallel =\omega_1 +\dots +\omega_r$.\\

On a de plus la propri\'et\'e suivante du battage de deux suites :

\begin{lem}
Pour toutes suites $\un{\omega}^1$ et $\un{\omega}^2$ avec
$l(\un{\omega}^1 )=n>0$ et $l(\un{\omega}^2 )=m>0$, l'ensemble
$\sh (\un{\omega}^1 ,\un{\omega}^2 )$ est la r\'eunion disjointe
des deux ensembles suivants:
\begin{equation}
\bigl( \sh (\un{\omega}^1 ,\un{\omega}^{2}_{<m-1 } ) \omega_m^2
\bigl) \coprod \bigl( \sh (\un{\omega}^1_{<n-1 } ,\un{\omega}^2 )
\omega_n^1 \bigl) .
\end{equation}
\end{lem}

La d\'emonstration est laiss\'ee au lecteur. On a donc
\begin{equation}
\left .
\begin{array}{lll}
\di \sum_{\un{\omega} \in \sh (\un{\omega}^1 ,\un{\omega}^2 )}
S^{\un{\omega}} & = & \di {-1\over \parallel \omega \parallel}
\sum_{\un{\omega} \in \sh (\un{\omega}^1 ,\un{\omega}^2 )}
S^{\un{\omega}_{<r-1}} , \\
& = & -\di {1\over \parallel \omega \parallel} \sum_{\un{\un{s}}
\in \sh (\un{\omega}^1_{<n-1} ,\un{\omega}^2 )} S^{\un{s}} - \di
{1\over \parallel \omega \parallel} \sum_{\un{s} \in
\sh (\un{\omega}^1 ,\un{\omega}^2_{<m-1} )} S^{\un{s}} .
\end{array}
\right .
\end{equation}
Comme $l(\un{s})=r-1$, on en d\'eduit par hypoth\`ese de
r\'ecurrence,
\begin{equation}
\left .
\begin{array}{lll}
\di \sum_{\un{\omega} \in \sh (\un{\omega}^1 ,\un{\omega}^2 )}
S^{\un{\omega}} & = & -\di {1\over \parallel \omega \parallel}
S^{\un{\omega}^1_{<n-1}} S^{\un{\omega}^2} -\di {1\over \parallel
\omega
\parallel}
S^{\un{\omega}^1} S^{\un{\omega}^2_{<m-1}}  ,\\
& = & \di {\parallel \un{\omega}^1 \parallel \over \parallel
\omega
\parallel}
S^{\un{\omega}^1 } S^{\un{\omega}^2} +\di {\parallel \un{\omega}^2
\parallel \over \parallel \omega
\parallel}
S^{\un{\omega}^1} S^{\un{\omega}^2} ,
\end{array}
\right .
\end{equation}
en utilisant la relation (\ref{relS}). En simplifiant, on obtient
finalement,
\begin{equation}
\left .
\begin{array}{lll}
\di \sum_{\un{\omega} \in \sh (\un{\omega}^1 ,\un{\omega}^2 )}
S^{\un{\omega}} & = & \di {\parallel \un{\omega}^1 \parallel
+\parallel \un{\omega}^2
\parallel \over \parallel \omega
\parallel}
S^{\un{\omega}^1 } S^{\un{\omega}^2 } , \\
& = &  S^{\un{\omega}^1 } S^{\un{\omega}^2 } ,
\end{array}
\right .
\end{equation}
ce qui termine la preuve.
\end{proof}

\subsection{Un moule sym\'etrel}

On d\'efinit le moule
\begin{equation}
Se^{\omega_1 \dots \omega_r} =\di {e^{\parallel \un{\omega }
\parallel} \over (e^{-\omega_1} -1 ) \dots
(e^{-\parallel \un{\omega}_{<i} \parallel} -1) (e^{-\parallel
\un{\omega} \parallel} -1)} ,
\end{equation}
o\`u $\un{\omega}_{<i} =\omega_1 \dots \omega_i $.

\begin{lem}
Le moule $Se^{\bullet}$ est sym\'etrel.
\end{lem}

\begin{proof}
Elle se fait par r\'ecurrence sur la longueur des s\'equences. On
commence par noter que
\begin{equation}
\label{petres} Se^{\un{\omega }} =\di {e^{\omega_r} \over
(e^{-\parallel \un{\omega}
\parallel} -1 )} Se^{\un{\omega}_{<r-1}} .
\end{equation}
D'apr\`es le lemme \ref{propbatco}, on a
$$
\left .
\begin{array}{lll}
\di\sum_{\un{\omega}\in csh(\un{\omega}^1 ,\un{\omega}^2 ) }
Se^{\un{\omega}} & = & \di\sum_{\un{\omega} \in csh (\un{\omega}^1
,\un{\omega}^2_{<r-1} )\omega_r^2 } Se^{\un{\omega}}
+\di\sum_{\un{\omega} \in csh (\un{\omega}^1_{<m-1} ,\un{\omega}^2
)\omega_1 } Se^{\un{\omega}} \\
& & +\di\sum_{\un{\omega} \in csh(\un{\omega}^1_{<m-1}
,\un{\omega}^2_{<r-1} ) (\omega_1^m +\omega_2^r ) }
Se^{\un{\omega}} ,
\end{array}
\right .
$$
d'o\`u, en utilisant (\ref{petres}),
$$
\left .
\begin{array}{lll}
\di\sum_{\un{\omega}\in csh(\un{\omega}^1 ,\un{\omega}^2 ) }
Se^{\un{\omega}} & = & \di {e^{\omega_r^2 } \over (e^{-\parallel
\un{\omega} \parallel } -1)} Se^{\un{\omega}^1}
Se^{\un{\omega}^2_{<r-1}} + \di {e^{\omega_m^1 } \over
(e^{-\parallel \un{\omega} \parallel } -1)}
Se^{\un{\omega}^{1,<m-1} }
Se^{\un{\omega}^2 } \\
& & + \di {e^{\omega_m^1 +\omega_r^2 } \over (e^{-\parallel
\un{\omega}
\parallel } -1)} Se^{\un{\omega}^{1,<m-1} }
Se^{\un{\omega}^{2,<r-1}} , \\
& = & \di {e^{-\parallel \omega^2 \parallel} -1 \over
(e^{-\parallel \un{\omega} \parallel } -1)} Se^{\un{\omega}^1}
Se^{\un{\omega}^2 }  + \di {e^{-\parallel \omega^1 \parallel} -1
\over (e^{-\parallel \un{\omega}
\parallel } -1)} Se^{\un{\omega}^1}
Se^{\un{\omega}^2 } \\
& & + \di {(e^{-\parallel \omega^1 \parallel} -1 ) (e^{-\parallel
\omega^2 \parallel} -1 ) \over (e^{-\parallel \un{\omega}
\parallel } -1)} Se^{\un{\omega}^1}
Se^{\un{\omega}^2 } , \\
& = &  Se^{\un{\omega}^1} Se^{\un{\omega}^2 } .
\end{array}
\right .
$$
On a donc le lemme.
\end{proof}

%% file: moule12.tex
\part{Alg\`ebre \`a composition des moules}
\setcounter{section}{0}
\setcounter{thm}{0}
\setcounter{lem}{0}
\setcounter{defi}{0}
\setcounter{rema}{0}
\setcounter{equation}{0}

La correspondance entre moules et s\'eries formelles non commutatives permet de
munir l'ensemble des moules d'une structure d'alg\`ebre non commutative. On
introduit aussi une op\'eration, appell\'ee {\it composition}, et qui est
l'analogue {\it non commutatif} de la substitution des s\'eries formelles.
Cette op\'eration n'existe pas dans les travaux combinatoires usuels, comme par
exemple dans l'\'etude des alg\`ebres de Hopf. On introduit aussi les groupes
alterna(e)l et symetra(e)l.

\section{Structure d'alg\`ebre}

Soit $X$ un ensemble d'\'el\'ements $X_{\omega}$ indic\'es par un
semi-groupe $\Omega$. On suppose $X$ muni d'un coproduit not\'e
$\Delta$. On note $\kK \langle\langle X\rangle\rangle$ l'alg\`ebre des s\'eries formelles non
commutatives form\'ees sur $X$ muni du coproduit $\Delta$. On note
${\cal M}_{\kK} (\Omega )$ l'ensemble des moules sur $\Omega$. La structure
d'alg\`ebre de $\kK \langle\langle X\rangle\rangle$ se traduit directement  sur les moules.\\

Soient $\di\sum_{\bullet} M^{\bullet } D_{\bullet}$ et
$\di\sum_{\bullet} N^{\bullet} D_{\bullet }$  deux \'el\'ements de
$\kK \langle\langle X\rangle\rangle$. On d\'efinit l'addition et la multiplication de deux
moules via les relations
\begin{equation}
\left .
\begin{array}{lll}
\di\sum_{\bullet} M^{\bullet } D_{\bullet} + \di\sum_{\bullet}
N^{\bullet } D_{\bullet} & = &
\di\sum_{\bullet} (M^{\bullet } +N^{\bullet } ) D_{\bullet} ,\\
\left ( \di\sum_{\bullet} N^{\bullet } D_{\bullet} \right ) \left
( \di\sum_{\bullet} M^{\bullet } D_{\bullet} \right ) & = &
\di\sum_{\bullet} (M^{\bullet } \times N^{\bullet } ) D_{\bullet}
.
\end{array}
\right .
\end{equation}
On peut donc munir ${\cal M}_{\kK} (\Omega )$ de la structure d'alg\`ebre
suivante :

\begin{thm}
L'ensemble des moules ${\cal M}_{\kK} (\Omega )$ muni des op\'erations
$$
\left .
\begin{array}{lll}
A^{\bullet } = M^{\bullet} +N^{\bullet } & \Longleftrightarrow &
A^{\un{\omega} } =M^{\un{\omega } } +
N^{\un{\omega } } \\
A^{\bullet } = M^{\bullet } \times N^{\bullet }  &
\Longleftrightarrow & A^{\un{\omega} } = \di\sum_{\un{\omega}^1
\bullet \un{\omega}^2 =\un{\omega} }  M^{\un{\omega}^1 }
N^{\un{\omega}^2 } ,
\end{array}
\right .
$$
est une alg\`ebre non commutative.

L'\'el\'ement neutre pour la multiplication est le moule
$1^{\bullet}$  d\'efini par
$$1^{\omega } =1\ \mbox{\rm si}\ \omega =\emptyset\  \mbox{\rm et}\ 1^{\omega } =0\
\mbox{\rm sinon.}$$
\end{thm}

On peut pr\'eciser la relation entre les moules alternaux et
sym\'etraux via le th\'eor\`eme \ref{repres}.

\begin{defi}
Soit $M^{\bullet }$ un moule, on appelle exponentielle de
$M^{\bullet }$ et on note $\exp (M^{\bullet } )$ la s\'erie $\exp
(N^{\bullet } )= \di\sum \di {(M^{\bullet } )^n \over n!}$, avec
la convention $(M^{\bullet } )^0 =1^{\bullet }$.
\end{defi}

Une simple application de la formule de Leibniz donne

\begin{lem}
\label{exist} Pour tout moule sym\'etral $M^{\bullet }$ il existe
un moule alternal $N^{\bullet }$ tel que $M^{\bullet } =\exp
N^{\bullet }$.
\end{lem}

\begin{proof}
Il suffit de voir que $\exp (\di\sum_{\bullet } N^{\bullet }
D_{\bullet } ) = \di\sum_{\bullet } \exp N^{\bullet } D_{\bullet
}$ par d\'efinition du moule exponentielle. Le th\'eor\`eme
\ref{repres} permet de conclure.
\end{proof}

\section{Composition}

On peut munir l'alg\`ebre des moules d'une composition. Cette loi
de composition est l'analogue de la notion de {\it substitution}
dans l'alg\`ebre des s\'eries formelles (voir \cite{D}, annexe 21,
p.398-400).

\begin{defi}
Soit $\Omega$ un semi-groupe. On note $\parallel .\parallel$ l'application de $\Omega$ dans $\Omega$ d\'efinie pour tout
$\un {\omega} =\omega_1 \dots \omega_r$ par
\begin{equation}
\parallel \omega_1 \dots \omega_r \parallel =\omega_1 +\dots +\omega_ r .
\end{equation}
Soit $M^{\bullet}$ et $N^{\bullet}$ deux moules dans ${\cal M}_{\kK} (\Omega )$. Le moule
compos\'e $A^{\bullet} =M^{\bullet} \circ N^{\bullet}$  est
d\'efini par
\begin{equation}
\label{comp}
A^{\omega} =\di\sum_{s\geq 0, \omega^1 \dots \omega^s
=\omega ,\ \omega^i \not= \emptyset} \di M^{\parallel \omega^1
\parallel ,\dots ,\parallel \omega^s \parallel } N^{\omega^1 }
\dots N^{\omega^s } ,
\end{equation}
L'\'el\'ement neutre pour la composition est le moule $I^{\bullet }$ d\'efini par
\begin{equation}
I^{\omega } =1\ \mbox{\rm si}\ l(\omega )=1\ \mbox{\rm et}\ I^{\omega } =0\
\mbox{\rm sinon}.
\end{equation}
\end{defi}

Cette op\'eration d'apparence compliqu\'ee, peut s'expliquer de la mani\`ere suivante:\\

Soient $M$ et $N$ les deux s\'eries formelles associ\'ees aux moules $M^{\bullet}$ et $N^{\bullet}$:
\begin{equation}
\left .
\begin{array}{lll}
M & = & \di\sum_{\un{\omega} \in \Omega^*} M^{\un{\omega}} \un{\omega} ,\\
N & = & \di\sum_{\un{\omega} \in \Omega^*} N^{\un{\omega}} \un{\omega} .
\end{array}
\right .
\end{equation}
L'application $\parallel .\parallel$ permet de construire un nouvel alphabet \`a partir de $M$:\\

\`A tout $\omega \in \Omega$, on associe la lettre
\begin{equation}
m_{\omega} =\di\sum_{\un{\omega} \in \Omega^*,\ \parallel \un{\omega}\parallel =\omega} M^{\un{\omega}} \un{\omega} .
\end{equation}

La s\'erie $N\circ M$ est alors d\'efinie comme suit:
\begin{equation}
\label{expcomp}
N\circ M=\di\sum_{\un{\omega}} N^{\un{\omega}} m_{\un{\omega}} .
\end{equation}

Avec ces notations, nous avons le r\'esultat suivant :

\begin{lem}
\label{compop}
La s\'erie $N\circ M$ poss\`ede un d\'eveloppement moulien de la forme
\begin{equation}
N\circ M=\di\sum_{\un{\omega} \in \Omega^*} (N^{\bullet }
\circ M^{\bullet } )^{\un{\omega}} \un{\omega}.
\end{equation}
\end{lem}

\begin{proof}
Il suffit de d\'evelopper l'expression (\ref{expcomp}). On a:
\begin{equation}
\label{expcomp2}
N\circ M=\di\sum_{\un{\omega}} N^{\omega_1 \dots \omega_r} \left (
\di\sum_{\un{\omega}^1 \in \Omega^*,\ \parallel \un{\omega}^1 \parallel =\omega_1} M^{\un{\omega}^1} \un{\omega}^1
\right ) \dots \left (
\di\sum_{\un{\omega}^r \in \Omega^*,\ \parallel \un{\omega}^r\parallel =\omega_r } M^{\un{\omega}^r} \un{\omega}^r
\right ) .
\end{equation}
Soit $\un{\omega} \in \Omega^*$ fix\'e. Le coefficient dans (\ref{expcomp2}) pour toute d\'ecomposition
\begin{equation}
\un{\omega}=\un{\omega}^1 \dots \un{\omega}^r  ,
\end{equation}
de $\un{\omega}$ est donn\'e par
\begin{equation}
N^{\parallel \un{\omega}^1 \parallel \dots \parallel \un{\omega}^r \parallel } M^{\un{\omega}^1} \dots M^{\un{\omega}^r} ,
\end{equation}
ce qui termine la preuve.
\end{proof}

On a le th\'eor\`eme principal de cette section:

\begin{thm}
L'alg\`ebre des moules muni des op\'erations $(+,\times ,\circ )$
est une alg\`ebre \`a composition, i.e.,
\begin{equation}
\left .
\begin{array}{lll}
i)\ (M^{\bullet } + N^{\bullet } )\circ A^{\bullet } & = &
(M^{\bullet } \circ A^{\bullet } ) +
(N^{\bullet } \circ A^{\bullet } ) , \\
ii)\ (M^{\bullet } \times N^{\bullet } )\circ A^{\bullet } & = &
(M^{\bullet } \circ A^{\bullet } )\times (N^{\bullet } \circ
A^{\bullet } ) .
\end{array}
\right .
\end{equation}
\end{thm}

La d\'emonstration repose sur des calculs \'el\'ementaires.

\section{Groupe alterna(e)l et sym\'etra(e)l}

Dans les calculs pratiques sur les moules, il est commode de
conna\^\i tre le comportement de la  propri\'et\'e d'alternalit\'e
(resp. sym\'etralit\'e) vis \`a vis
des lois de composition et  de multiplication. \\

On commence par noter le simple r\'esultat suivant :

\begin{lem}
L'ensemble des moules $M^{\bullet }$ tels que $M^{\emptyset }
\not= 0$ forment un groupe, not\'e $M_{\times} (\Omega )$,  dont
les \'el\'ements sont les moules poss\`edant un inverse
multiplicatif.
\end{lem}

\begin{proof}
Soient $M^{\bullet }$ et $N^{\bullet }$ deux moules dans
$M_{\times} (\Omega )$. On note $A^{\bullet } =M^{\bullet } \times
N^{\bullet }$. Par d\'efinition, on a $A^{\emptyset}
=M^{\emptyset} N^{\emptyset}$.  Comme $M^{\emptyset} \not=0$ et
$N^{\emptyset } \not= 0$, on en d\'eduit $A^{\emptyset } \not= 0$.

La caract\'erisation des moules ayant un inverse se fait par
r\'ecurrence sur la longueur des suites. Soit $M^{\bullet }$  un
moule poss\`edant un inverse multiplicatif not\'e $N^{\bullet }$,
alors on doit avoir $1^{\bullet } =M^{\bullet } \times N^{\bullet
}$, soit
$$1^{\omega } =\di\sum_{\omega_1 \omega_2 =\omega }
M^{\omega_1 } N^{\omega_2 }.$$ On a donc une r\'ecurrence pour
d\'eterminer le moule $N^{\bullet}$ via la relation
$$1^{\omega } =\di\sum_{\omega_1 \omega_2 =\omega,\ \omega_1 \not=\emptyset }
M^{\omega_1 } N^{\omega_2 } + M^{\emptyset} N^{\omega } .$$ La
condition $\omega_1 \not=0$ fait que les moules $N^{\omega_2 }$
intervenant dans la somme sont associ\'es \`a des suites de
longueur inf\'erieur \`a celle de $\omega$. On peut donc trouver,
de mani\`ere it\'erative, l'expression du moule $N^{\bullet}$ si
et seulement si $M^{\emptyset } \not=0$.
\end{proof}

Soit $M^{\bullet} \in M_{\times } (\Omega )$, on notera
$(M^{\bullet} )^{-1}$ le moule inverse.

\begin{lem}
L'ensemble des moules $M^{\bullet}$ tel que $M^{\emptyset } =0$ et
$M^{\omega } \not=0$ si $l(\omega )=1$ forment un  groupe, not\'e
$M_{\circ} (\Omega )$, dont les \'el\'ements sont les moules
poss\`edant un inverse pour la  composition.
\end{lem}

\begin{proof}
Soient $M^{\bullet }$ et $N^{\bullet }$ deux moules dans
$M_{\circ} (\Omega )$. On note $A^{\bullet } =M^{\bullet } \circ
N^{\bullet }$. Par d\'efinition, on a $A^{\emptyset}
=M^{\emptyset} N^{\emptyset}$.  Comme $M^{\emptyset} =0$ et
$N^{\emptyset } = 0$, on en d\'eduit $A^{\emptyset } = 0$. De
plus, $A^{\omega } = M^{\omega } N^{\omega }$ si $l(\omega )=1$.
Comme $M^{\omega } \not= 0$ et $N^{\bullet } \not= 0$, on a
$A^{\omega } \not =0$, et $A^{\bullet } \in M_{\circ} (\Omega )$.

La caract\'erisation des moules ayant un inverse pour la
composition se fait par r\'ecur\-ren\-ce sur la longueur des
suites. Soit $M^{\bullet }$ un moule poss\`edant un inverse de
composition not\'e $N^{\bullet }$, alors on doit avoir $I^{\bullet
} =M^{\bullet } \circ N^{\bullet }$, soit
$$I^{\omega } =\di\sum_{r\geq 1}
\sum_{\omega_1 \dots \omega_r =\omega }  M^{\parallel \omega_1
\parallel ,\dots ,\parallel \omega_r \parallel } N^{\omega_1 }
\dots N^{\omega_r }.$$ On a donc une r\'ecurrence pour
d\'eterminer le moule $N^{\bullet}$ via la relation
$$I^{\omega } =\di\sum_{r>1} \di\sum_{\parallel \omega_1 \parallel ,\dots ,\parallel
\omega_r \parallel } N^{\omega_1 } \dots N^{\omega_r } +
M^{\parallel \omega
\parallel } N^{\omega } .$$
La condition $r>1$ fait que les moules $N^{\omega_1 }$, $\dots$,
$N^{\omega_r}$ intervenant dans la somme sont associ\'es \`a des
suites de longueur inf\'erieur \`a celle de $\omega$. On peut donc
trouver, de mani\`ere it\'erative, l'expression du moule
$N^{\bullet}$ si et seulement si $M^{\parallel \omega \parallel }
\not=0$  pour tout $\omega$, soit $M^{\omega } \not= 0$ lorsque
$l(\omega )=1$.
\end{proof}

Soit $M^{\bullet} \in M_{\circ } (\Omega )$, on notera
$(M^{\bullet} )^{(-1)}$ le
moule inverse. \\

Le r\'esultat principal de cette section est :

\begin{thm}
Les ensembles $M_{sym} (\Omega )$ des moules sym\'etrals, muni de
la multiplication, et $M_{alt} (\Omega )$ des moules alternals,
muni de la composition, sont des sous-groupes non distingu\'es  de
$M_{\times} (\Omega )$ et $M_{\circ} (\Omega )$ respectivement.
\end{thm}

La d\'emonstration repose sur le lemme suivant :

\begin{lem}
Soient $A^{\bullet} \in M_{sym} (\Omega )$, $B^{\bullet } \in
M_{sym} (\Omega )$, $C^{\bullet}\in M_{alt} (\Omega )$  et
$E^{\bullet } \in M_{alt} (\Omega )$, on a
$$
\left .
\begin{array}{lll}
\mbox{\rm i)}\ (A^{\bullet} )^{-1} \in M_{sym} (\Omega ), & \ \ \
&  \mbox{\rm ii)} \ (C^{\bullet} )^{(-1)}
\in M_{alt} (\Omega ), \\
\mbox{\rm iii)}\ A^{\bullet} \times B^{\bullet }  \in M_{sym}
(\Omega ), & & \mbox{\rm iv)}\ C^{\bullet}
\circ E^{\bullet } \in M_{alt} (\Omega ) , \\
\mbox{\rm v)}\ B^{\bullet} \times C^{\bullet} \times (B^{\bullet}
)^{-1} \in M_{alt} (\Omega ) . & &
\end{array}
\right .
$$
\end{lem}

\begin{proof}
On peut d\'emontrer i) et iii) en utilisant la formule de
Campbell-Hausdorff. Soit $x$ et $y$ deux \'el\'ements de
$\hat{L}_X$, alors $\exp (x)\exp (y)=\exp z$, o\`u $z\in
\hat{L}_X$.  Tout moule sym\'etral s'obtient par le lemme
\ref{exist} comme exponentielle d'un moule alternal. La formule de
Campbell-Hausdorff assure que $(\di\sum_{\bullet} A^{\bullet }
D_{\bullet } ) (\di\sum_{\bullet} B^{\bullet } D_{\bullet } )
=(\di\sum_{\bullet} (A^{\bullet } \times B^{\bullet } ) D_{\bullet
} )$ s'\'ecrit sous la forme $(\di\sum_{\bullet} \exp F^{\bullet }
D_{\bullet } )$ o\`u $F^{\bullet }$ est un moule alternal. Par
d\'efinition, on a donc $\exp F^{\bullet} \in M_{sym} (\Omega )$.

De m\^eme, on obtient facilement $(A^{\bullet} )^{-1} \in M_{sym} (\Omega )$. \\

La propri\'et\'e v) exprime la stabilit\'e des d\'erivations via
une conjugaison par un automorphisme de $U$. Elle  est donc
\'evidente.

Les propri\'et\'es ii) et iv) repose sur le fait que toute
transformation associ\'ee \`a un moule alternal  transforme un
\'el\'ement de $\hat{L}_X$ en un \'el\'ement de $\hat{L}_X$. Elles
d\'ecoulent donc de la  d\'efinition m\^eme de la composition des
moules interpr\'et\'ee en terme d'op\'erateur dans le lemme
\ref{compop}.
\end{proof}

On note $ret\, :\ \Sigma_{\Omega } \rightarrow \Sigma_{\Omega }$
l'involution, qui a toute suite $\un{\omega} =(\omega_1 , \dots
,\omega_n )$ associe $ret {\un{\omega}} =(\omega_n ,\dots
,\omega_1 )$.

\begin{lem}
L'inverse multiplicatif d'un moule symetral $\di M^{\un{\omega}}$
est donn\'e par
$$\di (M^{\un{\omega}} )^{-1} =\di (-1)^{l(\un{\omega } )}
\di M^{ret{\un{\omega}}}.$$
\end{lem}

\begin{proof}
Elle se fait par r\'ecurrence sur la longueur de $\un{\omega }$.
Pour simplifier l'\'ecriture, on note $N^{\bullet}$ le moule inverse de $M^{\bullet}$. \\

Pour $\un{\omega} =\emptyset$, on a $1=M^{\emptyset}
N^{\emptyset}$, d'o\`u $N^{\emptyset } =1$. Pour $\un{\omega }
=\omega_1$, on a $1^{\omega_1} =M^{\omega_1} N^{\emptyset}
+M^{\emptyset} N^{\omega_1}$, d'o\`u $N^{\omega_1 } =-M^{\omega_1
}$. Un calcul int\'eressant est donn\'e par la longueur 2, qui
fournit la clef de la d\'emonstration. Pour $\un{\omega }
=\omega_1 \omega_2$, on a $1^{\omega_1 \omega_2} =M^{\omega_1
\omega_2 } N^{\emptyset}  +M^{\omega_1 } N^{\omega_2 }
+M^{\emptyset} N^{\omega_1 \omega_2 }$. En utilisant l'expression
de $N^{\omega_1 }$, on obtient $0=M^{\omega_1 \omega_2 }
-M^{\omega_1 } M^{\omega_2 } +N^{\omega_1 \omega_2 }$. La
propri\'et\'e de sym\'etralit\'e de $M^{\bullet }$ permet de
simplifier cette expression, et on obtient $0=M^{\omega_1 \omega_2
} -M^{\omega_1 \omega_2 } -M^{\omega_2 \omega_1 } +N^{\omega_1
\omega_2 }$, d'o\`u $N^{\omega_1 \omega_2 } =M^{\omega_2 \omega_1
}$.

Plus g\'en\'eralement, pour tout suite $\un{\omega}$ de longueur
$n \geq 1$, on a
\begin{equation}
1^{\un{\omega }} = \di\sum_{i=1}^n \di M^{\un{\omega}^{\leq i}}
N^{\un{\omega}^{>i}} + N^{\un{\omega}} ,
\end{equation}
o\`u $\un{\omega}^{<\leq i} =(\omega_1 ,\dots ,\omega_i )$ et
$\un{\omega }^{>i} =(\omega_{i+1} ,\dots ,\omega_n )$, avec les
conventions $\un{\omega} ^{\leq 0} =\emptyset$ et
$\un{\omega}^{>n} =\emptyset$.

Les suites $\di \un{\omega}^{>i}$ intervenant dans la somme sont
toutes de longueur $<n$. Par hypoth\`ese de r\'ecurrence, on a
$N^{\un{s}} =(-1)^{l(\un{s} )} M^{ret (\un{s} )}$ pour tout suite
$\un{s}$ telle que $l(\un{s} )< n$. On a donc
\begin{equation}
1^{\un{\omega} } = \di\sum_{i=1}^n \di (-1)^{l(\un{\omega}^{>i} )}
M^{\un{\omega}^{\leq i}} M^{ret (\un{\omega}^{>i} ) } +
N^{\un{\omega}} ,
\end{equation}
soit
\begin{equation}
1^{\un{\omega} } = M^{\un{\omega}} +\di\sum_{i=1}^{n-1} \di
(-1)^{l(\un{\omega}^{>i} )} \sum_{\un{s} \in \sh
(\un{\omega}^{\leq i} , ret (\un{\omega}^{>i} ) } M^{\un{s}} +
N^{\un{\omega}} .
\end{equation}

On note $S_i =sh(\di\un{\omega}^{\leq i} ,ret (\un{\omega }^{>i} )
)$, pour $i=1,\dots ,n-1$. On d\'efinit une suite $E_i$,
$i=1,\dots ,n-2$ d'ensembles par r\'ecurrence, tel que $S_1 =ret
(\un{\omega } ) \coprod E_1 $, $S_j =E_{j-1} \coprod E_j$, avec
$S_{n-1} =\un{\omega} \coprod E_{n-2}$.
Il n'est pas utile d'expliciter les ensembles $E_i$. \\

On a donc, pour toute suite $\un{\omega}$ de longueur $\geq 2$, la
relation
\begin{equation}
\left .
\begin{array}{lll}
0 & = & M^{\un{\omega}} +(-1)^{n-1} (M^{ret(\un{\omega} )}
+\di\sum_{\un{s} \in E_1} M^{\un{s}} ) +(-1)^{n-2}
(\di\sum_{\un{s} \in E_1} M^{\un{s}} +\di\sum_{\un{s} \in E_2}
M^{\un{s}} )
+\dots \\
 & & + (-1) (\di\sum_{\un{s} \in E_{n-2}} M^{\un{s}} +M^{\un{\omega}} ) +N^{\un{\omega}} .
\end{array}
\right .
\end{equation}
On en d\'eduit,
\begin{equation}
N^{\un{\omega}} =(-1)^n M^{ret(\un{\omega} )} ,
\end{equation}
ce qui termine la d\'emonstration du lemme.
\end{proof}

%% file: moule2.tex
\setcounter{section}{0}
\setcounter{thm}{0}
\setcounter{lem}{0}
\setcounter{defi}{0}
\setcounter{rema}{0}
\setcounter{equation}{0}

\part{Sym\'etries secondaires et d\'erivations}

On introduit la sym\'etrie sym\'etril/alternil qui intervient dans des travaux r\'ecents de Jean Ecalle sur la combinatoire
des polyz\^etas. On d\'emontre aussi que l'on peut lire {\it directement} la sym\'etrie de moules v\'erifiants certaines
\'equations diff\'erentielles. Au passage, on d\'efinit quelques d\'erivations et automorphismes importants sur
l'alg\`ebre des moules.

\section[Alternil et sym\'etril]{Sym\'etries alternil et sym\'etril}

Les sym\'etries alternal(el)/sym\'etral(el) sont uniquement
li\'ees \`a la traduction du caract\`ere primitif ou ``group like"
des s\'eries formelles non commutatives, les coproduits $\Delta$
et $\Delta_*$ \'etant donn\'es. Les sym\'etries
alternil/sym\'etril sont d'une autre nature : elles proviennent de
l'utilisation de {\it s\'eries g\'en\'eratrices} associ\'ees \`a
des moules alternel/sym\'etrel. Ces deux sym\'etries ont donc un
statut particulier vis \`a vis des quatre premi\`eres.

\begin{defi}
Soit $Se^{\bullet}$ un moule sym\'etral. La s\'erie
g\'en\'eratrice\footnote{Des s\'eries de m\^eme type interviennent
d\'eja dans le travail de Jean Ecalle, via la m\'ethode d'{\it
amplification} (voir \cite{ES}).} associ\'ee \`a $Se^{\bullet}$
est un moule $Sig^{\bullet}$ d\'efini par
\begin{equation}
Sig^{v_1 ,\dots ,v_r} =\di\sum_{1\leq s_i} Se^{s_1 ,\dots ,s_r}
v_1^{s_1 -1} \dots v_r^{s_r -1} .
\end{equation}
\end{defi}

La notation $Sig$ fait r\'ef\'erence \`a la sym\'etrilit\'e de
$Sig^{\bullet}$.

\begin{defi}
Un moule $M^{\bullet}$ est sym\'etril (resp. alternil) si
\begin{equation}
\di\sum_{\un{s} \in\shi (\un{x} ,\un{y} )} M^{\un{s}} =M^{\un{x}}
M^{\un{y}}\ (\mbox{\rm resp.}\ 0),
\end{equation}
o\`u $\shi (\un{x} ,\un{y} )$ s'obtient comme $\csh (\un{x}
,\un{y} )$ en rempla\c{c}ant l'addition des variables par un symbole
abstrait $*$ d\'ecrivant l'\'evaluation de $M^{\un{s}}$ suivant la
r\`egle
\begin{equation}
M^{\dots ,x_i*y_j ,\dots} =\di {1\over x_i -y_j} \left ( M^{\dots
,x_i ,\dots} -M^{\dots ,y_j ,\dots} \right ),
\end{equation}
les termes en pointill\'es pouvant comporter eux aussi le symbole
$*$.
\end{defi}

Le rapport aux s\'eries g\'en\'eratrices est donn\'e par:

\begin{lem}
Le moule $Sig^{\bullet}$ est sym\'etril.
\end{lem}

\begin{proof}
Il suffit de calculer $Sig^{v_1 ,\dots ,v_r} Sig^{v_{r+1} \dots
v_n}$. Par d\'efinition, on a
\begin{equation}
Sig^{v_1 ,\dots ,v_r} Sig^{v_{r+1} \dots v_n} =\di\sum_{1\leq s_i}
Se^{s_1 ,\dots s_r} Se^{s_{r+1} ,\dots ,s_n} v_1^{s_1 -1} \dots
v_n^{s_n -1} .
\end{equation}
Comme le moule $Se^{\bullet}$ est sym\'etrel, on a
$$Se^{s_1 ,\dots s_r} Se^{s_{r+1} ,\dots ,s_n} =\di\sum_{\un{s} \in\csh ((s_1 ,\dots ,s_r ),
(s_{r+1} ,\dots ,s_n ))} Se^{\un{s}},$$ soit
\begin{equation}
\label{tril} Sig^{v_1 ,\dots ,v_r} Sig^{v_{r+1} \dots v_n}
=\di\sum_{1\leq s_i} \di\sum_{\un{s} \in\csh ((s_1 ,\dots ,s_r ),
(s_{r+1} ,\dots ,s_n ))} Se^{\un{s}}\, v_1^{s_1 -1} \dots v_n^{s_n
-1} .
\end{equation}
Cette somme contient deux types de termes :\\

i - des termes de la forme $\di\sum_{1\leq s_i} Se^{\un{s}}\,
v_1^{s_1 -1} \dots v_n^{s_n -1}$, o\`u $\un{s} \in \csh
((s_1 ,\dots ,s_r ), (s_{r+1} ,\dots ,s_n ))$.\\

ii- des termes de la forme $\di\sum_{1\leq s_i}
Se^{\dots s_i +s_j \dots }\, v_1^{s_1 -1} \dots v_n^{s_n -1}$.\\

Pour i), on obtient, via une r\'eorganisation des $v_i$,
$$\di\sum_{\un{v} \in\csh ((v_1,\dots ,v_r ),(v_{r+1} ,\dots ,v_n ))} Sig^{\un{v}}.$$
Pour ii), c'est un peu plus compliqu\'e.\\

Nous allons le faire sur un exemple. Le cas g\'en\'eral s'en
d\'eduisant sans peine. Soit
$$Z(v_1 ,\dots ,v_n )= \di\sum_{1\leq s_i} Se^{s_1+s_2 ,s_3 ,\dots ,s_n} v_1^{s_1 -1}
\dots v_n^{s_n -1}.$$ On a
\begin{equation}
\left .
\begin{array}{lll}
v_1 Z(v_1 ,\dots ,v_n ) & = & \di\sum_{1\leq s_i} Se^{s_1+s_2 ,s_3
,\dots ,s_n} v_1^{s_1}
\dots v_n^{s_n -1} , \\
 & = & \di\sum_{1\leq s_i ,\, s_2=1} Se^{s_1+1 ,s_3 ,\dots ,s_n} v_1^{s_1}
\dots v_n^{s_n -1} \\
 & & + \di\sum_{1\leq s_i ,\, s_2 \geq 2} Se^{s_1+s_2 ,s_3 ,\dots ,s_n} v_1^{s_1}
\dots v_n^{s_n -1} , \\
 & = & Sig (v_1 ,v_3 ,\dots ,v_n ) -\di\sum_{1\leq s_i } Se^{1 ,s_3 ,\dots ,s_n} v_3^{s_3 -1}
\dots v_n^{s_n -1} \\
 & & +V_1 (v_1 ,\dots ,v_n ) ,
\end{array}
\right .
\end{equation}
o\`u $V_1 (v_1 ,\dots ,v_n ) =  \di\sum_{1\leq s_i ,\, s_2 \geq 2}
Se^{s_1+s_2 ,s_3 ,\dots ,s_n} v_1^{s_1} \dots v_n^{s_n -1}$.

Par ailleurs, on a
\begin{equation}
\left .
\begin{array}{lll}
v_2 Z(v_1 ,\dots ,v_n ) & = & \di\sum_{1\leq s_i} Se^{s_1+s_2 ,s_3
,\dots ,s_n} v_1^{s_1 -1}
v_2^{s_2} \dots v_n^{s_n -1} , \\
 & = & \di\sum_{1\leq s_i ,\, s_1=1} Se^{s_1+1 ,s_3 ,\dots ,s_n} v_1^{s_1 -1} v_2^{s_2}
\dots v_n^{s_n -1}\\
 & & + \di\sum_{1\leq s_i ,\, s_1 \geq 2} Se^{s_1+s_2 ,s_3 ,\dots ,s_n} v_1^{s_1 -1}
v_2^{s_2} \dots v_n^{s_n -1} , \\
 & = & Sig (v_2 ,\dots ,v_n ) -\di\sum_{1\leq s_i } Se^{1,s_3 ,\dots ,s_n} v_3^{s_3 -1}
\dots v_n^{s_n -1} \\
 & & +V_2 (v_1 ,\dots ,v_n ) ,
\end{array}
\right .
\end{equation}
o\`u $V_2 (v_1 ,\dots ,v_n ) =  \di\sum_{1\leq s_i ,\, s_1 \geq 2}
Se^{s_1+s_2 ,s_3 ,\dots ,s_n} v_1^{s_1 -1} v_2^{s_2} \dots
v_n^{s_n -1}$.

En posant $s_1 ' =s_1 -1$, $s_2 ' =s_2 +1$ dans $V_2$, on montre
que
\begin{equation}
V_1 =V_2 .
\end{equation}
Finalement, on a
\begin{equation}
Z(v_1 ,\dots ,v_n )=\di {1\over v_1 -v_2} (Sig^{v_1 v_3 \dots v_n}
-Sig^{v_2 \dots v_n} ).
\end{equation}
Le cas g\'en\'eral s'en d\'eduit sans peine.
\end{proof}

\begin{rema}
i. La sym\'etrie alternil ne provient pas de la sym\'etrie
alternel traduite sur les fonctions
g\'en\'eratrices. \\

ii. On peut se demander si la sym\'etrie symetral donne naissance
\`a une nouvelle sym\'etrie. En fait, il est clair, vue
l'\'equation (\ref{tril}), avec $\csh$ remplac\'e par $\sh$, que
la fonction g\'en\'eratrice d'un moule sym\'etral est {\it encore}
un moule sym\'etral.
\end{rema}

\section{D\'erivations et sym\'etries des moules}

\subsection{D\'erivations sur l'alg\`ebre des moules}

Dans les applications, on est souvent conduit \`a utiliser des
d\'erivations sur l'alg\`ebre des moules. Ce paragraphe donne, en
suivant Ecalle \cite{EV}, un proc\'ed\'e de construction d'une
grande quantit\'e de d\'erivations, suffisantes pour la plupart
des applications.\\

On commence par une d\'efinition:

\begin{defi}
Soit $D$ une application de ${\cal M}_{\kK} (\Omega )$ dans ${\cal M}_{\kK} (\Omega)$ lin\'eaire. Le moule image
d'un moule $M^{\bullet}$ par $D$ est not\'e $D(M)^{\bullet}$. L'application $D$
est une d\'erivation sur l'alg\`ebre $({\cal M}_{\kK} (\Omega ),\times )$ si elle v\'erifie
\begin{equation}
D(M\times N)^{\bullet} =D(M)^{\bullet}\times N^{\bullet} +
M^{\bullet} \times D(N)^{\bullet} ,
\end{equation}
pour tout moules $M^{\bullet}$ et $N^{\bullet}$ de ${\cal M}_{\kK} (\Omega )$.
\end{defi}

Certaines d\'erivations respectent des sym\'etries des moules.

\begin{defi}
On dit que la d\'erivation est alternale si elle pr\'eserve
l'alternalit\'e du moule sur lequel elle agit.
\end{defi}

On v\'erifie directement que les d\'erivations $\td_{\omega}$ et
$\mbox{\rm lang}$ d\'efinies ci apr\`es sont alternales.

\subsection{Construction de d\'erivations}

\subsubsection{D\'erivations simples}

On peut chercher \`a construire des d\'erivations simples de la
forme
$$(D_{\lambda} M^{\bullet} )^{\bs} =\lambda_{\bs} M^{\bs} ,$$
o\`u $\lambda_{\bullet}$ est un moule fix\'e.\\

Quelles sont les propri\'et\'es que doit v\'erifier
$\lambda_{\bullet}$ pour que
l'application $D_{\lambda}$ ci-dessus soit une d\'erivation ?\\

Il suffit de calculer $D_{\lambda} (M^{\bullet} \times N^{\bullet}
)^{\bullet}$. On a
$$(D_{\lambda} (M^{\bullet} \times M^{\bullet} )^{\bs} =\lambda_{\bs} \left ( \di\sum_{\bs^1 \bs^2 =\bs}
M^{\bs^1} N^{\bs^2} \right ).$$ Si $D_{\lambda}$ est une
d\'erivation, on a l'identit\'e de Leibniz qui impose
$$
\left .
\begin{array}{lll}
(D_{\lambda} (M^{\bullet} \times M^{\bullet} )^{\bs} & = &
((D_{\lambda} M^{\bullet} ) \times N^{\bullet} )^{\bs}
+(M^{\bullet} \times (D_{\lambda} N^{\bullet} ) )^{\bs} ,\\
 & = & \di\sum_{\bs^1 \bs^2 =\bs} \lambda_{\bs^1} M^{\bs1} N^{\bs^2} +
\di\sum_{\bs^1 \bs^2 =\bs} M^{\bs1} \lambda_{\bs^2} N^{\bs^2} ,\\
 & = & \di\sum_{\bs^1 \bs^2 =\bs} (\lambda_{\bs^1} +\lambda_{\bs^2} ) M^{\bs1} N^{\bs^2} .
\end{array}
\right .
$$
On a donc le th\'eor\`eme suivant:

\begin{thm}
L'application $D_{\lambda}$ es une d\'erivation si et seulement si
$\lambda_{\bullet}$ v\'erifie
$$\lambda_{\bs^1 \bs^2} =\lambda_{\bs^1} +\lambda_{\bs^2},\ \forall\ \bs^1 ,\bs^2 \in \Omega^*.$$
\end{thm}

Il existe deux exemples importants de d\'erivations simples:\\

- On note $\mbox{\rm lang}$ la d\'erivation simple d\'efinie par le moule
\begin{equation}
\lambda_{\bullet} =l(\bullet ) ,
\end{equation}
o\`u $l(\bullet )$ est l'application {\it longueur de la suite}, i.e.
\begin{equation}
(\mbox{\rm lang} M)^{\un{s} } =  l(\un{s} ) M^{\un{s} } .
\end{equation}

- On note $\td_{\omega }$ la d\'erivation simple d\'efinie par le moule
\begin{equation}
\lambda_{\bullet} =\parallel \bullet \parallel .
\end{equation}
Notons que cette expression a un sens si et seulement si $\Omega$ poss\`ede une structure de semi-groupe. Par ailleurs, comme
$\lambda_{\un{s}} M^{\un{s}}$ dit appartenir \`a $\kK$, cela impose $\un{s} \in \kK^*$.

On a donc
\begin{equation}
(\td M)^{\un{s} }  =  \parallel \un{s} \parallel M^{\un{s} } .
\end{equation}

Ces d\'erivations interviennent dans la th\'eorie des formes normales de champs de vecteurs ou diff\'eomorphismes locaux.

\subsubsection{Autres constructions}

Soit un moule $Dar^{\bullet }$ avec $Dar^{\emptyset} =0$. On
d\'efinit un op\'erateur $dar$ de ${\cal M}_{\kK} (\Omega )$ dans ${\cal M}_{\kK} (\Omega )$
en posant :
\begin{equation}
\label{derdar}
\left .
\begin{array}{llll}
dar\, & : \, {\cal M}_{\kK} (\Omega ) & \rightarrow & {\cal M}_{\kK} (\Omega ) , \\
 & M^{\bullet} & \mapsto & (dar M)^{\omega } =\di\sum_{\omega_1
\omega_2 \omega_3 =\omega} M^{\omega_1 ,\parallel \omega_2
\parallel ,\omega_3 } Dar^{\omega_2 } ,
\end{array}
\right .
\end{equation}
o\`u la somme est d\'efinie sur toutes les factorisations de
$\omega$.

\begin{lem}
L'op\'erateur $dar$ est une d\'erivation sur ${\cal M}_{\kK} (\Omega ,\times )$
\end{lem}

\begin{proof}
On note $E^{\bullet} =dar (M_1^{\bullet } \times M_2^{\bullet } )$
et $F^{\bullet } =M_1^{\bullet } \times M_2^{\bullet }$. On a
\begin{equation}
\label{dar} E^{\omega } =\di\sum_{\omega =\omega_1 \omega_2
\omega_3} F^{\omega_1 ,\parallel \omega_2 \parallel ,\omega_3}
Dar^{\omega_2 }.
\end{equation}
Comme $F^{\omega_1 ,\parallel \omega_2 \parallel ,\omega_3 }
=\di\sum_{\hat{\omega } \tilde{\omega} =\omega_1 ,
\parallel \omega_2 \parallel ,\omega_3 } M_1^{\hat{\omega
}} M_2^{\tilde{\omega }}$, on a
$$F^{\omega_1 ,\parallel \omega_2 \parallel ,\omega_3 } =  \di\sum_{\parallel \omega
\parallel \in
\hat{\omega }}  M_1^{\hat{\omega }} M_2^{\tilde{\omega }}
+\di\sum_{\parallel \omega_2 \parallel \in \tilde{\omega}}
M_1^{\hat{\omega }} M_2^{\tilde{\omega }}.$$ En rempla\c{c}ant dans
(\ref{dar}), on obtient
$$E^{\omega } = \di\sum_{\omega =\omega_1 \omega_2 \omega_3} =
(\di\sum_{\parallel \omega \parallel \in \hat{\omega }}
M_1^{\hat{\omega }} Dar^{\omega_2} ) M_2^{\tilde{\omega }}
+\di\sum_{\parallel \omega_2
\parallel \in \tilde{\omega}}  M_1^{\hat{\omega }} (M_2^{\tilde{\omega }}
Dar^{\omega_2 } ).
$$
Comme $\hat{\omega }$ est de la forme $\hat{\omega }
=\hat{\omega_1} ,\parallel \omega_2 \parallel , \hat{\omega_2 }$
dans le premier terme de la somme, et tel que $\hat{\omega_1}
\hat{\omega_2} \tilde{\omega} = \omega$ (une expression  analogue
pour $\tilde{\omega }$ dans le second terme), on a finalement
$$E^{\omega } =\di\sum_{\omega =\hat{\omega } \tilde{\omega }} (dar
M_1)^{\hat{\omega }} M_2^{\tilde{\omega}} + \di\sum_{\omega
=\hat{\omega } \tilde{\omega }}  M_1^{\hat{\omega}} (dar
M_2)^{\tilde{\omega }},$$ ce qui termine la d\'emonstration.
\end{proof}

On peut choisir le moule $Dar^{\bullet}$ pour que la d\'erivation
$dar$ soit alternale :

\begin{lem}
La d\'erivation $dar$ d\'efinie par (\ref{derdar}) est alternale
si et seulement si le moule $Dar^{\bullet }$ est alternal.
\end{lem}

Il existe d'autres d\'erivations construites avec l'op\'erateur
$diff$ d\'efini par
\begin{equation}
\left .
\begin{array}{lll}
\mbox{\rm diff}\, :\, M(\Omega ) & \rightarrow & M(\Omega ) ,\\
M^{\bullet} & \mapsto & (\mbox{\rm diff} M)^{\un{w}} =\di\sum_{w_i
\in \un{w}} \di {\partial M^{\un{w}} \over \partial_{w_i }}  ,
\end{array}
\right .
\end{equation}
ce qui suppose que $M^{\un{w}}$ soit diff\'erentiable en chacune
des variables
$w_i \in \un{w}$.

\begin{rema}
Une alg\`ebre de Hopf sur un corps $\kK$ est un sch\`ema en groupe
affine sur $K$ (voir \cite{wat}, p. 9). Il serait interessant
d'avoir une interpr\'eation fonctorielle sur le groupe d'une
d\'erivation moulienne.
\end{rema}

\subsection{D\'erivations et moules sym\'etrals}

Jean Ecalle m'a sugg\'er\'e le r\'esultat suivant permettant de
d\'emontrer sans trop d'efforts, la sym\'etrie de nombreux moules.

\begin{thm}
\label{dermou}
Soit $D_{\lambda}$ une d\'erivation et $M^{\bullet}$ un moule v\'erifiant l'\'equation diff\'erentielle:
\begin{equation}
\label{dm}
D_{\lambda} M^{\bullet} =A^{\bullet} \times M^{\bullet} ,
\end{equation}
avec:\\

i) Le moule $\lambda^{\bullet}$ v\'erifie $\lambda_{\bs} \not= 0$ pour tout $\bs$, $l(\bs )\geq 1$.

ii) On a $M^{\emptyset} =1$,

iii) le moule $A^{\bullet}$ est alternal,\\

alors, le moule $M^{\bullet}$ est symetral.
\end{thm}

\begin{proof}
Commen\c{c}ons par introduire une notion qui nous sera utile pour la suite:\\

Le moule $M^{\bullet}$ est sym\'etral jusqu'\`a l'ordre $r$ si il v\'erifie la condition de sym\'etralit\'e
$$\di\sum_{\un{s}\in \sh (\un{s}^1 ,\un{s}^2 )} M^{\un{s}} =M^{\un{s}^1} M^{\un{s}^2} ,$$
pour tout couple de suites $(\un{s}^1 ,\un{s}^2 )$ telles que $l(\un{s}^1 )+l(\un{s}^2 )=r$.\\

Supposons $M^{\bullet}$ sym\'etral jusqu'\`a l'ordre $r$. Soient $\un{s}^1$ et $\un{s}^2$ deux suites telles que
$$l(\un{s}^1 )+l(\un{s}^2 )=r+1 .$$
On a
\begin{equation}
\lambda_{\un{s}^1 \un{s}^2 } \di\sum_{\un{s} \in \sh (\un{s}^1 \un{s}^2 )} M^{\un{s}} =
\di\sum_{\un{s} \in \sh (\un{s}^1 \un{s}^2 )} \lambda_{\un{s}} M^{\un{s}} ,
\end{equation}
car
\begin{equation}
\lambda_{\un{s}} =\lambda_{\un{s}^1 \un{s}^2 },\ \forall \un{s} \in \sh (\un{s}^1 ,\un{s}^2 ) ,
\end{equation}
soit
\begin{equation}
\lambda_{\un{s}^1 \un{s}^2 } \di\sum_{\un{s} \in \sh (\un{s}^1 \un{s}^2 )} M^{\un{s}} =\di
\di\sum_{\un{s} \in \sh (\un{s}^1 \un{s}^2 )} \left (
\di\sum_{\un{u}\un{v} =\un{s}} A^{\un{u}} M^{\un{v}} \right ) .
\end{equation}
On obtient donc
\begin{equation}
\left .
\begin{array}{lll}
\lambda_{\un{s}^1 \un{s}^2 } \di\sum_{\un{s} \in \sh (\un{s}^1 \un{s}^2 )} M^{\un{s}} & = &
\di\sum_{i=0}^{l(\un{s}^1 )} A^{\un{s}^{1,\leq i}} \di\sum_{\un{v}\in\sh (\un{s}^{1,>i} ,\un{s}^2 )} M^{\un{v}} \\
 & & + \di\sum_{i=0}^{l(\un{s}^2 )} A^{\un{s}^{2,\leq i}} \di\sum_{\un{v}\in\sh (\un{s}^1 ,\un{s}^{2,>i} )} M^{\un{v}}
\end{array}
\right .
\end{equation}
Par sym\'etralit\'e de $M^{\bullet}$ d'ordre $r$, on en d\'eduit
\begin{equation}
\lambda_{\un{s}^1 \un{s}^2 } \di\sum_{\un{s} \in \sh (\un{s}^1 \un{s}^2 )} M^{\un{s}} =
\di\sum_{i=0}^{l(\un{s}^1 )} A^{\un{s}^{1,\leq i}} M^{\un{s}^{1,>i}} M^{\un{s}^2}
+ \di\sum_{i=0}^{l(\un{s}^2 )} A^{\un{s}^{2,\leq i}} M^{\un{s}^1} M^{\un{s}^{2,>i}},
\end{equation}
soit
\begin{equation}
\left .
\begin{array}{lll}
\lambda_{\un{s}^1 \un{s}^2 } \di\sum_{\un{s} \in \sh (\un{s}^1 \un{s}^2 )} M^{\un{s}} & = &
M^{\un{s}^2} \left ( \di\sum_{i=0}^{l(\un{s}^1 )} A^{\un{s}^{1,\leq i}} M^{\un{s}^{1,>i}} \right )
+M^{\un{s}^1} \left ( \di\sum_{i=0}^{l(\un{s}^2 )} A^{\un{s}^{2,\leq i}} M^{\un{s}^{2,>i}} \right ) , \\
 &= & \lambda_{\un{s}^1} M^{\un{s}^1} M^{\un{s}^2} +\lambda_{\un{s}^2} M^{\un{s}^1} M^{\un{s}^2} ,\\
 & =& M^{\un{s}^1} M^{\un{s}^2}(\lambda_{\un{s}^1} +\lambda_{\un{s}^2} ) , \\
  & =& M^{\un{s}^1} M^{\un{s}^2} \lambda_{\un{s}^1 \un{s}^2} .
\end{array}
\right .
\end{equation}
Comme $\lambda_{\un{s}} \not =0$ pour toute suite $\un{s} \in \Omega^* \setminus \{ \emptyset \}$, on obtient
\begin{equation}
\di\sum_{\un{s}\in \sh (\un{s}^1 ,\un{s}^2 )} M^{\un{s}} =M^{\un{s}^1 } M^{\un{s}^2} ,
\end{equation}
donc, la sym\'etralit\'e de $M^{\bullet}$ \`a l'ordre $r+1$. Une simple r\'ecurrence termine la preuve.
\end{proof}

\begin{rema}
L'\'equation (\ref{dm}) est sugg\'er\'ee par le
r\'esultat classique suivant. Pour tout $x\in \kK << \Omega >>$,
avec $x$ primitif, on a
$$D(e^x )=\di\sum_{k\geq 1} \di {1\over k!} ad(x)^{k-1} (Dx) e^x = \left ( \di{ e^{ad(x)} -1 \over e^{ad(x)}} \right )
(Dx) e^x ,$$ avec $ad(x)(y)=xy-yx$. Du point de vue moulien, on
obtient, en notant $x=\di\sum_{w\in \Omega^*} A^w w$, $e^x
=\di\sum_{w\in \Omega^*} exp(A)^w w=\di\sum_{w\in \Omega^*} M^w
w$,
$$D(M^{\bullet} )=\di\sum_{k\geq 1} \di {1\over k!} ad(A^{\bullet})(DA^{\bullet} ) M^{\bullet}.$$
\end{rema}

Ce lemme permet notamment de d\'emontrer, sans aucun calcul sur
les suites, la sym\'etralit\'e du moule
$$S^{\omega_1 ,\dots ,\omega_r } =\di { (-1)^r \over \omega_1 (\omega_1 +\omega_2 ) \dots (\omega_1 +\dots +
\omega_r )},$$ \'etudi\'e au \S 5.1. En effet, on a
$$\nabla S^{\bullet} =S^{\bullet } \times I^{\bullet },$$
avec $I^{\bullet }$ evidemment alternal.

Par ailleurs, le seul moule $M^{\bullet}$ v\'erifiant $\nabla
M^{\bullet} =0$ pour tout suite $\un{\omega }$ telle que
$\parallel \un{\omega } \parallel \not= 0$, est le moule
constant \'egal \`a $0$. On en d\'eduit donc que $S^{\bullet }$ est sym\'etral.

\section{Automorphismes et sym\'etries}

\subsection{D\'efinition}

Commen\c{c}ons par une d\'efinition:

\begin{defi}
Soit $A$ une application de ${\cal M}_{\kK} (\Omega )$ dans ${\cal M}_{\kK} (\Omega)$ lin\'eaire. Le moule image
d'un moule $M^{\bullet}$ par $A$ est not\'e $A(M)^{\bullet}$. L'application $A$
est un automorphisme de l'alg\`ebre $({\cal M}_{\kK} (\Omega ),\times )$ si elle v\'erifie
\begin{equation}
\label{defauto}
A(M\times N)^{\bullet} =A(M)^{\bullet}\times A(N)^{\bullet} ,
\end{equation}
pour tout moules $M^{\bullet}$ et $N^{\bullet}$ de ${\cal M}_{\kK} (\Omega )$.
\end{defi}

On peut chercher \`a construire des automorphismes simples de la forme
\begin{equation}
\label{autof}
A_f (M)^{\bullet} =f(\bullet )M^{\bullet} ,
\end{equation}
o\`u $f$ est une application de $\Omega^*$ dans $\kK$. L'application $f$ doit bien entendu satisfaire certaines
contraintes.

\begin{lem}
L'application $A_f$ de ${\cal M}_{\kK} (\Omega )$ dans ${\cal M}_{\kK} (\Omega)$ d\'efinie par (\ref{autof}) est un
automorphisme de ${\cal M}_{\kK} (\Omega )$ si et seulement si $f$ est un morphisme de $(\Omega^* ,\bullet )$ dans
$(\kK ,\times )$.
\end{lem}

\begin{proof}
Il suffit d'\'etudier la relation induite par la relation (\ref{defauto}). Soient $M^{\bullet}$ et $N^{\bullet}$
deux moules. On obtient pour toute suite $\un{\omega} \in \Omega^*$
\begin{equation}
f(\un{\omega} )\di\sum_{\un{\omega}^1 \un{\omega}^2 =\un{\omega}} M^{\un{\omega}^1} N^{\un{\omega}^2} =
\di\sum_{\un{\omega}^1 \un{\omega}^2 =\un{\omega}} f(\un{\omega}^1 ) M^{\un{\omega}^1} f(\un{\omega}^2 )
N^{\un{\omega}^2} ,
\end{equation}
soit
\begin{equation}
f(\un{\omega} )=f(\un{\omega}^1 )f(\un{\omega}^2 ) .
\end{equation}
L'application $f$ est donc un morphisme de $(\Omega^* ,\bullet )$ dans $(\kK ,\times )$.
\end{proof}

\subsection{Exemples}

Soit $\Omega$ un alphabet de lettres $\omega \in \rR$. Pour tout moule $M^{\bullet}$ sur $M(\Omega )$, et pour toute suite
$\un{\omega} \in \Omega^*$, on d\'efinit une application $e^{\td}$ par
\begin{equation}
e^{\td} M^{\un{\omega}} =e^{\parallel \un{\omega} \parallel}
M^{\un{\omega}} .
\end{equation}

\begin{rema}
Cet automorphisme intervient naturellement dans les probl\`emes de formes normales des diff\'eomorphismes locaux.
\end{rema}

\subsection{Automorphismes et moules symetrels}

Le lemme suivant est l'analogue de la relation liant derivations et moules symetrals, dans le cas des automorphismes et
des moules symetrels :

\begin{lem}
Soit $\Omega$ un alphabet, $M(\Omega )$ l'alg\`ebre des moules sur
$\Omega$, et $A_f$ un automorphisme admissible sur $M(\Omega )$.
Soit $M^{\bullet}$ un moule tel que
\begin{equation}
\label{eq}
A_f (M)^{\bullet}  =(1^{\bullet} +I^{\bullet} ) \times M^{\bullet} .
\end{equation}
Alors le moule $M^{\bullet}$ est sym\'etrel.
\end{lem}

\begin{proof}
Pour toute suite $\un{\omega}$ on a
\begin{equation}
f(\un{\omega}) M^{\un{\omega}} =M^{\un{\omega}} +M^{\un{\omega}^{\geq 2}} ,
\end{equation}
soit
\begin{equation}
M^{\un{\omega}} =\di {1\over f(\un{\omega} )-1} M^{\un{\omega}^{\geq 2}} .
\end{equation}
On a donc
\begin{equation}
\left .
\begin{array}{lll}
\di\sum_{\un{\omega} \in \csh (\un{u},\un{v} )} M^{\un{\omega}}
 & = & \di\sum_{\un{\omega} \in \csh (\un{u},\un{v} )} \di {1\over f(\un{\omega} )-1} M^{\un{\omega}^{\geq 2}} ,\\
 & = & \di {1\over f(\un{u}\un{v} )-1} \di\sum_{\un{\omega} \in \csh (\un{u},\un{v} )} M^{\un{\omega}^{\geq 2}} ,
\end{array}
\right .
\end{equation}
car la valeur de $f$ sur un mot ne d\'epend pas de l'ordre des lettres.\\

Les propri\'et\'es du battage contractant, via le lemme \ref{propbatco}, donnent
\begin{equation}
\left .
\begin{array}{lll}
\di\sum_{\un{\omega} \in \csh (\un{u},\un{v} )} M^{\un{\omega}}
 & = & \di {1\over f(\un{u}\un{v} )-1}
\left [
M^{\un{u}} M^{\un{v}^{\geq 2}} +M^{\un{u}^{\geq 2}} M^{\un{v}} +M^{\un{u}^{\geq 2}} M^{\un{v}^{\geq 2}}
\right ] ,\\
 & = & \di {1\over f(\un{u}\un{v} )-1}
\left [
f(\un{v} )-1 +f(\un{u})-1 +(f(\un{u} ) -1)(f(\un{v})-1)
\right ] ,\\
 & = & M^{\un{u}} M^{\un{v}} ,
\end{array}
\right .
\end{equation}
car $f(\un{u}\un{v} )=f(\un{u} )f(\un{v})$ par hypoth\`ese.
\end{proof}

On peut sans doute d\'emontrer un th\'eor\`eme analogue dans le cas de l'\'equation
\begin{equation}
A_f (M)^{\bullet} =(1^{\bullet} +E^{\bullet} ) M^{\bullet} ,
\end{equation}
o\`u $E^{\bullet}$ est un moule alternel, mais la complexit\'e des calculs est beaucoup plus importante.

\section{Dualit\'e des alg\`ebres $\aA$ et $\eE$}

Le th\'eor\`eme qui suit donne une application permettant de passer de
$\eE$ \`a $\aA$ (et vice versa). La sym\'etrie symetrel/alternel
est beaucoup plus complexe et couteuse en terme de calculs que la
sym\'etrie symetral/alternal. Il peut donc \^etre utile de passer
de $\eE$ dans $\aA$ pour v\'erifier des propri\'et\'es de
sym\'etrie.

\begin{thm}
\label{dualalg}
Soit $Se^{\bullet}$ (resp. $E$) un moule symetrel
(resp. alternel) et $\Exp^{\bullet}$ le moule exponentiel d\'efini
par
\begin{equation}
\Exp^{\emptyset} =1,\ \ \Exp^{s_1 \dots s_r} =1/r! \ \ \forall \
\bs =s_1 \dots s_r \not= \emptyset .
\end{equation}
Le moule
\begin{equation}
Sa^{\bullet} =Se^{\bullet} \circ \Exp^{\bullet} ,\ (resp. \ \
A^{\bullet} =E^{\bullet} \circ \Exp^{\bullet} ),
\end{equation}
est un moule symetral (resp. alternal).
\end{thm}

Autrement dit, il existe une {\it dualit\'e} entre les alg\`ebre $\aA$ et $\eE$ via
la composition par le moule exponentielle.\\

La d\'emonstration repose sur le lemme suivant:

\begin{lem}
Soit $\Omega$ un semi-groupe. Le moule exponentielle $\Exp^{\bullet}$ d\'efinit un alphabet $Y$ tel que chaque lettre $Y_{\omega}$
est groupe-like par rapport au co-produit $\Delta$, i.e.
\begin{equation}
\Delta (Y_{\omega} )=\di\sum_{\omega_1 +\omega_2 =\omega} Y_{\omega_1} \otimes Y_{\omega_2} .
\end{equation}
\end{lem}

\begin{proof}
Soit $X$ un alphabet de type $\Delta$ et
$\bx$ un mot de $X^*$. On associe \`a chaque lettre $x_i$ un poids
$p_i \in \nN$. Le poids d'un mot est la somme des poids de ces
lettres. On note $\parallel \bx \parallel$ le poids. On \'etend
cette d\'efinition
\`a $\aA$ par lin\'earit\'e. \\

On a
\begin{equation}
\Exp =\di\sum_{\bx \in X^*} \Exp^{\bx} \bx =\di\sum_{i\in \nN} y_i
,
\end{equation}
o\`u les $y_i$ sont les composantes homog\`enes de poids $i$ de $\Exp$.\\

Or, l'alphabet $(y_i)$ v\'erifie
\begin{equation}
\Delta (y_i )=\di\sum_{l+k=i} y_l \otimes y_k .
\end{equation}
En effet, on a
\begin{equation}
\label{nitia} \Delta (y_i)=\Delta \left ( \di\sum_{\bx ,\
\parallel \bx \parallel =i} \Exp^{\bx} \bx \right ) .
\end{equation}
On a
\begin{equation}
\label{decoA} \Delta (\bx )= \sum_{\bx^1 , \bx^2 \mid\ \bx \in \sh
(\bx^1 ,\bx^2 )} \bx^1 \otimes \bx^2 .
\end{equation}
Comme la somme (\ref{nitia}) fait intervenir tous les $\bx$ tels
que $\parallel \bx\parallel =i$, on peut la re\'ecrire, en
utilisant (\ref{decoA}) :
\begin{equation}
\Delta (y_i)=\di\sum_{\bx^1 ,\bx^2 \mid \ \parallel \bx^1
\parallel +\parallel \bx^2 \parallel =i} Exp^{\bx^1} \Exp^{\bx^2}
\bx^1 \otimes \bx^2 ,
\end{equation}
o\`u nous avons implicitement utilis\'e le fait que le moule
$\Exp^{\bullet}$ d\'epend
seulement de la longueur des suites. \\

En regroupant les termes, on a finalement :
\begin{equation}
\Delta (y_i )= \di\sum_{l+k=i} y_l \otimes y_k ,
\end{equation}
ce qui termine la d\'emonstration.
\end{proof}

La d\'emonstration du th\'eor\`eme \ref{dualalg} s'en d\'eduit comme suit:\\

Soit $Sa$ la s\'erie formelle non commutative associ\'ee au moule $Sa^{\bullet}$ sur $X$. Par d\'efinition de la composition
des moules, on a:
\begin{equation}
Sa=\di\sum_{x\in X^*} Sa^x x =\di\sum_{Y\in Y^*} Se^{Y} Y .
\end{equation}
Comme l'alphabet est constitu\'e de lettres groupe-like pour le coproduit $\Delta$ et le moule $Se^{\bullet}$
est symmetrel, alors la s\'erie $\di\sum_{Y\in Y^*} Se^{Y} Y$ est groupe-like pour le coproduit $\Delta$. On en d\'eduit
que le moule $Sa^{\bullet}$ est symmetral, puisqu'il d\'efinit un \'el\'ement group-like sur un alphabet primitif.\\

Le d\'emonstration du passage entre moule alternel et alternal est analogue.

%% file: moule3.tex
\setcounter{section}{0}
\setcounter{thm}{0}
\setcounter{lem}{0}
\setcounter{defi}{0}
\setcounter{rema}{0}
\setcounter{equation}{0}

\newpage
\part{Th\'eorie des formes normales d'objets locaux}

Le but de cette partie est de montrer le langage des moules et comoules en action sur le probl\`eme de la
recherche des formes normales d'objets analytiques locaux, comme les champs de vecteurs
et les diff\'eomorphismes de $\cC^{\nu}$. On d\'emontre les versions mouliennes de th\'eor\`emes classiques :
th\'eor\`eme de Poincar\'e, forme normale r\'esonante de Poincar\'e-Dulac,
th\'eor\`eme de Bruyno, et ceci, aussi bien pour les diff\'eomorphismes que les champs de vecteurs analytiques locaux.
Outre leurs int\'er\^ets conceptuels, ces d\'emonstrations donnent des expressions explicites des
normalisateurs\footnote{C'est \`a dire des changements de variables.} et permettent de mettre en \'evidence
des coefficients {\it universels}\footnote{Nous donnerons un sens pr\'ecis \`a cette terminologie dans le texte.} dans ces
probl\`emes, qui sont justement des moules. Au passage, on donne une pr\'esentation originale de la m\'ethode
d'arborification, qui restaure la convergence des s\'eries mouliennes, et dont le domaine d'applications d\'epasse
largement celui de la th\'eorie des formes normales.

\section{Objets locaux : champs de vecteurs et diff\'eomorphismes}

\subsection{Champs de vecteurs}

Soit $\nu \in \nN^*$. On consid\`ere un champ de vecteur de
$\cC^{\nu }$ de la forme
\begin{equation}
X=\di\sum_{i=1}^{\nu} X_i (x) \partial_{x_i} ,\ X_i (x) \in \cC \{
x\} ,
\end{equation}
avec la notation simplifi\'ee $\partial_{x_i } =\partial /\partial_{x_i }$. \\

Nous \'etudions les champs de vecteurs {\it locaux}, i.e. tels que
\begin{equation}
X_i (0)=0, \ i=1,\dots ,\nu .
\end{equation}

On note $X_{\rm lin}$ la partie lin\'eaire de $X$. En suivant Brjuno \cite{br}, on d\'ecompose le champ $X$ sous
la forme
\begin{equation}
\label{decom}
X=X_{\rm lin} +\di\sum_{\un{n} \in A(X)} D_{\un{n}} ,
\end{equation}
o\`u les $D_{\un{n}}$ sont des {\it op\'erateurs homog\`enes} de
degr\'e $\un{n}$, avec $\un{n} =(n_1 ,\dots ,n_{\nu} )$, et tous
les $n_i \geq 0$, sauf au plus un qui peut valoir $-1$, i.e.
\begin{equation}
D_{\un{n}} (x^{\un{m}} )=c_{\un{n} ,\un{m}} x^{\un{n} +\un{m}} ,\
c_{\un{n} ,\un{m}} \in \cC ,\ \un{m} \in \nN^{\nu } .
\end{equation}
On note $A(X)$ l'ensemble des degr\'es des op\'erateurs homog\`enes intervenant dans la d\'ecompo\-si\-tion
(\ref{decom}).\\

On remarque que pour tout $\un{n} \in A(X)$, $D_{\un{n}}$ est une d\'erivation.\\

\subsubsection{Exemple}

Soit $X(x,y)$ un champ de vecteurs de $\cC^2$ de la forme
$$
\left .
\begin{array}{lll}
X(x,y) & = & \lambda x \partial_x +\beta y \partial_y + (a_{20} x^2 +a_{11} xy +a_{02} y^2) \partial_x \\
  & & +(b_{20} x^2 +b_{11} xy +b_{02} y^2 ) \partial_y .
\end{array}
\right .
$$
o\`u $a_{i,j} \in \cC$, $b_{i,j} \in \cC$ pour $i+j =2$, $i,j\in \nN$. \\

La d\'ecomposition (\ref{decom}) s'\'ecrit
$$X=X_{\rm lin} +D_{1,0} +D_{0,1} +D_{-1,2}+D_{2,-1},$$
o\`u $X_{\rm lin} = \lambda x\partial_x +\beta y\partial_y$, et
\begin{equation}
\left .
\begin{array}{lll}
D_{1,0} & = & a_{20} x^2 \partial_x +b_{11} xy\partial_y , \\
D_{0,1} & = & a_{11} yx\partial x +b_{02} y^2 \partial y , \\
D_{-1,2} & = & a_{02} y^2 \partial_x , \\
D_{2,-1} & = & b_{20} x^2 \partial_y .
\end{array}
\right .
\end{equation}
On a donc $A(X)=\{ (1,0),(0,1), (-1,2),(2,-1) \}$. \\

Dans la suite, on supposera toujours que la partie lin\'eaire du champ est sous forme {\it diagonale},
i.e.
$$X_{\rm lin} =\di\sum_{i=1}^{\nu} \lambda_i x_i \partial_{x_i} ,$$
o\`u $\blambda =(\lambda_1 ,\dots ,\lambda_{\nu} )$ est le {\it spectre} de $X_{\rm lin}$.\\

Le champ est alors dit {\it sous forme pr\'epar\'e} par Ecalle. \\

Cette condition est-elle restrictive ? Non, si on se place dans la classe des champs de vecteurs
formels. En effet, en suivant Martinet (\cite{mar},$\S$.1.1) tout champ de vecteur formel peut, via
un diff\'eomorphisme formel, se mettre sous la forme
$$X=X_{\rm lin} +X_N ,$$
o\`u $X_{\rm lin}$ est lin\'eaire diagonale et $X_N$ est nilpotent, avec $[X_{\rm lin} ,X_N ]=0$.
Evidemment, le champ $X_N$ n'est d\'etermin\'e que modulo l'action du groupe des diff\'eomorphismes
formels laissant $X_{\rm lin}$ invariant.

\subsection{Diff\'eomorphismes}

On consid\`ere un diff\'eomorphisme locale de $\cC^{\nu }$ de la
forme
$$f : (x_1 ,\dots ,x_{\nu } )\mapsto (f_1 (x) ,\dots ,f_{\nu } (x)),\ f_i (x)\in \cC \{ x\},$$
tel que
$$f_i (0) =0.$$
On associe \`a $f$ son {\it op\'erateur de substitution} :
$$F:\phi \rightarrow \phi \circ f.$$

On d\'emontre que $F$ peuvent se mettre sous la forme
\begin{equation}
\label{comp} \left .
\begin{array}{lll}
F & = & F_{lin} \left ( 1+\di\sum_n B_n \right ) , \\
F_{lin} : \phi & \mapsto & \phi (\lambda_1 x_1 ,\dots
,\lambda_{\nu } x_{\nu } ) ,
\end{array}
\right .
\end{equation}
o\`u les $B_{\un{n}}$ sont des op\'erateurs homog\`enes de degr\'e
$\un{n}$, $n=(n_1 ,\dots ,n_{\nu} )$, $n_i \geq 0$, sauf au plus
un qui peut valoir $-1$. On note $A(F)$ l'ensemble des degr\'es
des
op\'erateurs $B_{\un{n}}$ obtenus dans la d\'ecomposition (\ref{comp}). \\

On remarque que $B_{\un{n} }$, $\un{n} \in A(F)$, est un
op\'erateur {\it diff\'erentiel} (cela provient de la formule de
Taylor). Son coproduit est donc
$$\Delta_* (B_{\omega } )=\di\sum_{\omega_1 +\omega_2 =\omega } B_{\omega_1 } \otimes B_{\omega_2 } .$$

\begin{exemp}
On consid\`ere le diff\'eomorphisme de $\cC$ d\'efini par
$$f(x)=\lambda x +x^2.$$
Soit $\phi (x)=\di\sum_{n} a_n x^n$ un \'el\'ement de $\cC \{
x\}$. On veut expliciter l'op\'erateur de substitution $F : \phi
\rightarrow \phi \circ f$. On a
$$\phi \circ f (x)=\di\sum_n a_n (\lambda x +x^2 )^n.$$
Comme
$$(\lambda x +x^2 )^n =\di\sum-k C_n^k (x^2)^k (\lambda x)^{n-k},$$
on a, en regroupant correctement les termes,
$$\phi \circ f (x)=\di\sum_k \di {(x^2)^k \over k!} \di\sum_{n\geq k} a_n \di {n!\over (n-k)!}
(\lambda x)^{n-k}.$$ Par ailleurs, on a
$$\di {d^k \phi \over d x^k} (x) = \di\sum_{n\geq k}
a_n \di {n!\over (n-k)!} x^{n-k},$$ ce qui donne $F =F_{lin}
(1+\di\sum_{k\geq 1} B_k )$, en posant
$$B_k =\di {(x^2)^k \over k!} \di {d^k \over dx^k}.$$
\end{exemp}

\section{Conjugaison des objets analytiques locaux}

\subsection{Conjugaison}

Soit $X$ (resp. $F$) un champ de vecteur (resp. diff\'eomorphisme) sous forme pr\'epar\'ee. On
regarde l'effet d'un changement de variable (formel ou non) sur ces objets. On note $x$ les
variables initiales dans lesquelles sont explicit\'es $X$ et $F$.\\

On consid\`ere un changement de variable $x=h(y)$. On note $\Theta$ l'op\'erateur de substitution
associ\'e \`a $h$, d\'efini pour tout $\phi \in \cC [[x]]$ par
\begin{equation}
\Theta \phi =\phi \circ h ,\ \ \Theta^{-1} \phi =\phi \circ h^{-1} .
\end{equation}
La remarque importante est :

\begin{pro}
L'op\'erateur de changement de variable $\Theta$ est un automorphisme de $\cC [[x]]$.
\end{pro}

\begin{proof}
L'application $\Theta$ est $\cC$-lin\'eaire, et $h$ \'etant un
diff\'eomorphisme de $\cC [[x]]$, l'application $\phi \mapsto \phi
\circ h^{-1}$ est bien d\'efinie et correspond \`a $\Theta^{-1}$. Par
ailleurs, on a
$$
\left .
\begin{array}{lll}
\Theta (\phi \psi ) & = & (\phi \psi )\circ h ,\\
 & = & (\phi \circ h )(\psi \circ h) ,\\
  & = & \Theta (\phi ) \Theta (\psi ) ,
\end{array}
\right .
$$
et $\Theta$ est bien un morphisme de $\cC [[x]]$.
\end{proof}

\begin{rema}
Connaissant l'automorphisme $\Theta$, on retrouve le changement de
variable associ\'e en appliquant $\Theta$ \`a l'application
identit\'e de $\cC^{\nu}$.
\end{rema}

\begin{defi}
Soient $X$ et $F$ un  champ de vecteur (resp. diff\'eomorphisme)
analytique local de $\cC^{\nu}$. Un champ $X_{conj}$ (resp. un
diff\'eomorphisme $F_{conj}$) est dit conjugu\'e \`a $X$ (resp.
$F$), si il existe un changement de variable $h$ tel que
\begin{equation}
X_{\rm conj} =\Theta X \Theta^{-1} , (\mbox{\rm resp}.\ F_{\rm
conj} =\Theta F \Theta^{-1} ),
\end{equation}
o\`u $\Theta$ est l'op\'erateur de substitution associ\'e \`a $h$, i.e. si le diagramme suivant commute
\begin{equation}
\left .
\begin{array}{lll}
\cC \{ x\} & \stackrel{X}{\longrightarrow} & \cC \{ x\} ,\\
\uparrow \Theta & & \uparrow \Theta \\
\cC \{ y\} & \stackrel{X_{\rm conj}}{\longrightarrow} & \cC \{ y\} .
\end{array}
\right .
\end{equation}
\end{defi}

L'origine de ces \'equations de conjugaison est la suivante : \\

-- Un champ de vecteur $X$ en $x$ est une d\'erivation sur les germes de fonctions en $x$.
Soit $h$ un diff\'eomorphisme. L'image de $X$ par $h^{-1}$\footnote{Rappelez vous que l'on cherche
un changement de variable $x=h(y)$, o\`u $x$ est le syst\`eme de coordonn\'ees initial.}
est un champ de vecteur en $y=h^{-1}(x)$,
donc une d\'erivation sur les germes de fonctions en $y$. Notons $X_{\rm conj}$ ce champ
image. Comment lui donner un sens ? Soit $\phi$ un germe de fonction en $h^{-1}(x)$, alors
$\phi \circ h^{-1}$ est un germe de fonction en $x$. On peut faire agir $X$ dessus, et
on obtient le germe $X(\phi \circ h^{-1} )$ en $x$. Comme $X_{\rm conj} (\phi )$ doit \^etre
un germe en $h^{-1}(x)$, on transporte $X(\phi \circ h^{-1} )$ par $h$ \`a droite, soit
$X(\phi \circ h^{-1} ) \circ h$. Cette fonction est d\'efinie sur un voisinage de $y=h^{-1} (x)$, c'est
donc un germe de fonctions en $y$. Une d\'efinition est donc $X_{\rm conj} (\phi ) =
\Theta X\Theta^{-1} (\phi )$ pour tout germe $\phi$ en $h(x)$, o\`u $\Theta$ est l'automorphisme
de substitution associ\'e \`a $h$. On renvoie \`a (\cite{lafon},p.96) pour plus de d\'etails.\\

-- Pour les automorphismes, un raisonnement analogue au pr\'ec\'edent conduit au r\'esultat.\\

La conjugaison est dite {\it formelle} (resp. {\it analytique}) si le changement de variables
associ\'e est formel (resp. analytique).\\

Lorsque $X_{\rm conj} =X_{\rm lin}$ ou $F_{\rm conj} =F_{\rm
lin}$, on parle de {\it lin\'earisation}.

\begin{rema}
On peut imposer des contraintes sur la forme des objets
conjugu\'es. Si le champ conjugu\'e ne contient que des termes
{\it r\'esonants}, on parle de {\it pr\'enormalisation}. Si de
plus, le nombre de ces termes est minimal parmis toutes les formes
pr\'enormales, on parle de {\it normalisation}. On renvoie au
$\S$.\ref{prenormalisation} pour plus de d\'etails.
\end{rema}

\subsection{\'Equation de conjugaison}

Soit $X$ un champ de vecteur analytique local, d'alphabet $A(X)$ et $F$ l'op\'erateur de
substitution d'un diff\'eomorphisme analytique local $f$, d'alphabet $A(F)$. Soit $Na$ et $Ne$ des
automorphismes de conjugaison de $X$ et $F$ respectivement.\\

Comme $Na$ et $Ne$ sont des automorphismes de $\cC [[x]]$, on peut les chercher sous la forme
\begin{equation}
Na = \di\sum_{\un{s} \in A^* (X)} Na^{\un{s}} D_{\un{s}},
\end{equation}
avec $Na^{\bullet}$ un moule sym\'etral pour $X$, et
\begin{equation}
Ne =\di\sum_{\un{s} \in A^* (F)} Ne^{\un{s}} B_{\un{s}} ,
\end{equation}
o\`u $Ne^{\bullet}$ est un moule sym\'etrel pour $F$.\\

Via ces changements de variables, on obtient des objets conjugu\'es de la forme
\begin{equation}
X_{\rm conj} =X_{\rm lin} +\di\sum_{\un{s} \in A^* (X)} Ca^{\un{s}} D_{\un{s}} ,\ \
F_{\rm conj} =X_{\rm lin} +\di\sum_{\un{s} \in A^* (F)} Ce^{\un{s}} B_{\un{s}} ,
\end{equation}
o\`u les moules $Ca^{\bullet }$ et $Ce^{\bullet}$ sont alternal et sym\'etrel respectivement.\\

Dans ce cas, $X_{\rm conj}$ est bien une d\'erivation de $\cC [[x]]$, et $F_{\rm conj}$ un
automorphisme de $\cC [[x]]$. \\

L'\'equation de conjugaison pour un champ de vecteur et un diff\'eomorphisme s'\'ecrit donc en
terme moulien sous la forme
\begin{equation}
\label{equaconj}
\left .
\begin{array}{l}
X_{\rm lin} +\di\sum_{\bullet} Ca^{\bullet } D_{\bullet } = \left (
\di\sum_{\bullet} Na^{\bullet} D_{\bullet} \right )
\left ( X_{\rm lin} +\di\sum_{\bullet} I^{\bullet} D_{\bullet} \right )
\left ( \di\sum_{\bullet}( Na^{\bullet})^{-1} D_{\bullet } \right ) , \\
F_{\rm lin} +\di\sum_{\bullet} Ce^{\bullet } B_{\bullet } = \left (
\di\sum_{\bullet} Ne^{\bullet}  B_{\bullet} \right )
\left ( F_{\rm lin} +\di\sum_{\bullet} I^{\bullet} B_{\bullet} \right )
\left ( \di\sum_{\bullet} ( Ne^{\bullet} )^{-1} B_{\bullet } \right ) ,
\end{array}
\right .
\end{equation}
respectivement, o\`u $I^{\bullet}$ est le moule \'el\'ement neutre pour la composition.\\

\begin{thm}[Conjugaison]
\label{conjexpl}
Les \'equations de conjugaison (\ref{equaconj}) sont \'equivalentes aux \'equations mouliennes
\begin{equation}
Na(e)^{\bullet}  \times \nabla (Na(e) ^{-1} )^{\bullet} + Na(e)^{\bullet}
\times I^{\bullet} \times (Na(e)^{-1} )^{\bullet} = Ca(e)^{\bullet} ,
\end{equation}
o\`u $\nabla$ est la d\'erivation moulienne d\'efinie par
\begin{equation}
\left ( \nabla M^{\bullet} \right ) ^{\bs} =\parallel \bs.\blambda \parallel M^{\bs} .
\end{equation}
\end{thm}

La d\'emonstration de ce th\'eor\`eme repose sur le lemme technique suivant :

\begin{lem}
\label{tech1}
Soient $\phi (x)=\di\sum_{m\in \nN^{\nu} } a_m x^m \in \cC [[x]]$, et $B_{n^i}$, $i=1,\dots ,r$,
une famille d'op\'erateurs diff\'erentiels homog\`enes de degr\'e $n^i$ satisfaisant
\begin{equation}
B_{n^i} (x^m )=\beta_m^{n^i} x^{m+n^i},\ \beta_m^{n^i} \in \cC .
\end{equation}
. On a :
\begin{equation}
\label{tech11}
B_{\bn} \phi (x)=\di\sum_m a_m (\beta^{n^r}.m) (\beta^{n^{r-1}} .(m+n^r )) \dots
(\beta^{n^1} .(m+n^r +\dots +n^2 )) x^{m+n^r +\dots +n^1} .
\end{equation}
\end{lem}

\begin{proof}
On a
$$
\left .
\begin{array}{lll}
B_{\bn} \phi (x) & = & \di\sum_m a_m B_{\bn} (x^m ) ,\\
 & & \di\sum_m a_m B_{n^1} \dots B_{n^r} (x^m ) .
\end{array}
\right .
$$
Par ailleurs, on a pour tout $n^i$,
$$
\left .
\begin{array}{lll}
B_{n^i} \phi (x) & = & \di\sum_m a_m (B_{n^i} x^m ) ,\\
 & = & \di\sum_m a_m (\beta^{n^i} .m) x^{m+n^i} .
\end{array}
\right .
$$
Une r\'ecurrence imm\'ediate termine la d\'emonstration.
\end{proof}

On en d\'eduit le corollaire essentiel suivant :

\begin{cor}
\label{ess1} Pour toute suite $\bs$, on a
\begin{equation}
X_{\rm lin} D_{\bs } =\parallel \blambda. \bs \parallel D_{\bs}
+D_{\bs } X_{\rm lin} .
\end{equation}
\end{cor}

\begin{proof}
On a, en utilisant (\ref{tech11}) :
$$B_{\bn} \phi (x)=\di\sum_m a_m (\beta^{n^r}.m) (\beta^{n^{r-1}} .(m+n^r )) \dots
(\beta^{n^1} .(m+n^r +\dots +n^2 )) x^{m+n^r +\dots +n^1} .$$
Comme $X_{\rm lin} (x^m )=(\blambda.m) x^m$, on en d\'eduit
$$X_{\rm lin} B_{\bn} \phi (x) =\di\sum_m a_m
(\beta^{n^r}.m) \dots (\beta^{n^1} .(m+n^r +\dots +n^2 )) (\blambda .(m+n^r +\dots +n^1 ))
x^{m+n^r +\dots +n^1} ,$$
soit
$$
\left .
\begin{array}{lll}
X_{\rm lin} B_{\bn} \phi (x) & = & \di\sum_m a_m
(\beta^{n^r}.m) \dots (\beta^{n^1} .(m+n^r +\dots +n^2 )) (\blambda .\parallel \bn\parallel )
x^{m+n^r +\dots +n^1} \\
 & & +\di\sum_m a_m
(\beta^{n^r}.m) \dots (\beta^{n^1} .(m+n^r +\dots +n^2 )) (\blambda .m )
x^{m+n^r +\dots +n^1} .
\end{array}
\right .
$$
\end{proof}

\begin{proof}[D\'emonstration du th\'eor\`eme \ref{conjexpl}]
On a, en utilisant le corollaire \ref{tech1}
\begin{equation}
\label{ina1} X_{\rm lin} \left ( \di\sum_{\bullet} M^{\bullet}
D_{\bullet} \right ) = \di\sum_{\bullet} \nabla M^{\bullet}
D_{\bullet} + \di\sum_{\bullet} M^{\bullet} D_{\bullet} X_{\rm
lin} .
\end{equation}
La multiplication \`a gauche par $\Theta^{-1}$ donne
$$
\Theta^{-1} X_{\rm lin} \Theta
 =
\sum_{\bullet} ( \Theta^{\bullet} )^{-1} \times \nabla
\Theta^{\bullet} D_{\bullet} +\sum_{\bullet} (\Theta^{\bullet}
)^{-1} \Theta^{\bullet} D_{\bullet} X_{\rm lin} .
$$
Comme $(\Theta^{\bullet} )^{-1} \Theta =1^{\bullet}$, on a
$$\sum_{\bullet} (\Theta^{\bullet} )^{-1} \Theta^{\bullet} D_{\bullet} X_{\rm lin} =X_{\rm lin},$$
soit
$$
\Theta^{-1} X_{\rm lin} \Theta = \sum_{\bullet} ( \Theta^{\bullet}
)^{-1} \times \nabla \Theta^{\bullet}  D_{\bullet} +X_{\rm lin} .
$$
On a donc
$$X_{\rm lin} +\di\sum_{\bullet} Ca^{\bullet} D_{\bullet} =
X_{\rm lin} +\di\sum_{\bullet} ((Na^{-1} )^{\bullet} ) \nabla Na^{\bullet } D_{\bullet} +
\di\sum_{\bullet} ((Na^{-1} )^{\bullet} ) I^{\bullet} Na^{\bullet } D_{\bullet}.$$
On en d\'eduit le r\'esultat. Pour les diff\'eomorphismes, la d\'emonstration est analogue.
\end{proof}

\section{Lin\'earisation formelle}

Dans cette section, on red\'emontre le th\'eor\`eme de
lin\'earisation formelle de Poincar\'e. L'outil moulien permet de
renouveler son approche, en mettant en \'evidence des coefficients
{\it universels} de lin\'erisation, qui n'apparaissaient pas dans
la litt\'erature classique du sujet.

\subsection{Le th\'eor\`eme de Poincar\'e}

On cherche le normalisateur $Na$ tel que
$$X_{\rm lin} =Na X Na^{-1} ,$$
ou, ce qui revient au m\^eme
$$X=Na^{-1} X_{\rm lin} Na .$$
Par construction, nous supposons que $Na$ est de la forme
$$Na
=\di\sum_{\bullet} Na^{\bullet} D_{\bullet} ,$$
c'est \`a dire dans ${\cal D}$. On a $Na^{-1} =\di\sum_{\bullet}
(Na^{\bullet })^{-1} D_{\bullet }$. L'\'equation de
lin\'earisation donne
\begin{equation}
\left ( \di\sum_{\bullet} Na^{\bullet} D_{\bullet} \right )
X_{\rm lin} \left ( \di\sum_{\bullet} (Na^{\bullet } )^{-1}
D_{\bullet } \right ) =X_{\rm lin} +\di\sum_{\bullet} I^{\bullet}
D_{\bullet } ,
\end{equation}
o\`u $I^{\bullet}$ est le moule \'el\'ement neutre pour la
composition. \\

On en d\'eduit l'\'egalit\'e suivante sur les moules
\begin{equation}
\label{symetraconjugo}
\td (Na^{\bullet} )^{-1} \times Na^{\bullet} =I^{\bullet } ,
\end{equation}
o\`u $\td$ est la d\'erivation sur l'alg\`ebre des moules
d\'efinie par $\td M^{\un{\omega }} =\parallel \un{\omega}
\parallel M^{\un{\omega}}$. On a donc les formules de r\'ecurrence\footnote{Comme nous allons le voir, ces formules
d\'efinissent le moule $Na^{\bullet}$ par r\'ecurrence sur la {\it longueur} des mots, ce qui n'est pas \'evident a
priori.} suivantes:
\begin{equation}
\left .
\begin{array}{lll}
\td (Na^{\bullet})^{-1} & = & I^{\bullet} \times
(Na^{\bullet})^{-1} , \\
\td Na^{\bullet} & = & -Na^{\bullet} \times I^{\bullet } .
\end{array}
\right .
\end{equation}

Nous avons donc la version explicite\footnote{Le normalisateur est
donn\'e explicitement.} suivante du th\'eor\`eme de Poincar\'e :

\begin{thm}[Poincar\'e]
\label{thmpoinc}
Soit $X=X_{\rm lin} +\di\sum_{\un{n} \in A(X)}
D_{\un{n}}$ un champ de vecteur local de $\cC^{\nu}$, avec $X_{\rm lin}
=\di\sum_{i=1}^{\nu} \lambda_i x_i \partial_{x_i}$, $\un{\lambda}
=(\lambda_1 ,\dots ,\lambda_{\nu} )$. Pour tout $\un{s} =(s_1
,\dots ,s_r )\in A^* (X)$, $s_i \in A(X)$, on note $\omega (\un{s}
)=(\omega_1 , \dots ,\omega_r ) \in \cC^{\nu }$ le vecteur
d\'efini par $\omega_i =s_i .\un{\lambda }$, $i=1,\dots ,\nu$, et
$\parallel \omega (\un{s}) \parallel =\omega_1 +\dots +\omega_r$.\\

On suppose que le champ est non r\'esonant, i.e.
$$
\parallel \omega (\un{s} ) \parallel \not =0\ \mbox{\rm pour tout}\ \un{s} \in A^* (X) .
$$
Alors, le champ $X$ est formellement lin\'earisable, par un
automorphisme formel de $\cC [[x]]$ de la forme
$$Na = \di\sum_{\bullet} Na^{\bullet} D_{\bullet} ,$$
avec pour tout $\un{s} = (s_1 ,\dots ,s_r ) \in A^* (X)$, $s_i \in A(X)$,
$$
Na^{\un{s}} =\di {1 \over \omega_1 (\omega_1 +\omega_2
)\dots (\omega_1 +\omega_2 +\dots +\omega_r )},
$$
o\`u $\omega_i = s_i .\un{\lambda}$.
\end{thm}

\begin{proof}
On calcule $Na^{\bullet }$ par r\'ecurrence sur la longueur des
suites. On a $Na^{\emptyset } =1$ par hypoth\`ese.

Pour $l(\omega )=1$, on a
$$\omega Na^{\omega } =Na^{\omega }I^{\emptyset } + Na^{\emptyset } I^{\omega } =1,$$
d'o\`u
$$Na^{\omega } =\di {1\over \omega },$$
si $\omega \not= 0$.

Pour $l(\un{\omega })=r$, $\un{\omega } =(\omega_1 ,\dots
,\omega_r )$, on a la relation de r\'ecurrence
\begin{equation}
\parallel \un{\omega } \parallel Na^{\un{\omega }} =
Na^{\omega_1 ,\dots ,\omega_{r-1}} I^{\omega_r } ,
\end{equation}
d'o\`u, pour $\parallel \un{\omega } \parallel \not= 0$, on a
\begin{equation}
Na^{\un{\omega } } =\di {1\over \parallel \un{\omega}
\parallel } Na^{\omega_1 ,\dots ,\omega_{r-1} } .
\end{equation}
Une simple r\'ecurrence donne la forme g\'en\'erale de
$Na^{\bullet}$, \`a savoir
\begin{equation}
Na^{\omega_1 ,\dots ,\omega_r } =\di {1\over \omega_1
(\omega_1 +\omega_2 )\dots (\omega_1 +\omega_2 +\dots +\omega_r )} .
\end{equation}
Comme $Na^{\bullet}$ v\'erifie l'\'equation (\ref{symetraconjugo}), on sait d'apr\`es le th\'eor\`eme \ref{dermou} que le
moule $Na^{\bullet}$ est sym\'etral, i.e. que l'objet $Na$ est un automorphisme formel de $\cC [[x]]$. Ceci
termine la d\'emonstration du th\'eor\`eme.
\end{proof}

\begin{rema}
i. L'op\'erateur $Na$ est un objet formel. Pour \'etudier son \'eventuelle convergence,
Ecalle a introduit la m\'ethode d'{\it arborification}. Nous renvoyons \`a l'article
d'Ecalle \cite{E3} pour plus de d\'etails. \\

ii. Si le champ de vecteur est hamiltonien, le changement de variable est, par construction,
automatiquement symplectique.
\end{rema}

\subsection{Cas des diff\'eomorphismes}

L'\'equation de lin\'earisation est
$$Ne^{-1} F Ne =F_{\rm lin} ,$$
soit
$$Ne F_{lin} Ne^{-1} =F.$$
Via le calcul moulien, on a donc
$$\left
( \sum_{\bullet} Ne^{\bullet} B_{\bullet } \right ) F_{lin} \left
( \sum_{\bullet} (Ne^{\bullet} )^{-1} B_{\bullet} \right ) =\left
( \sum_{\bullet} (1^{\bullet} +I^{\bullet} )B_{\bullet} \right )
F_{lin} ,$$ o\`u $1^{\bullet }$ est le moule \'el\'ement neutre
pour la multiplication. On en d\'eduit la relation suivante sur
les moules
\begin{equation}
e^{\td} (Ne^{\bullet} )^{-1} \times Ne =1^{\bullet} +I^{\bullet} ,
\end{equation}
o\`u $e^{\td }$ est un automorphisme de l'alg\`ebre des moules
d\'efini par
$$(e^{\td} M^{\bullet } )^{\un{\omega}} =e^{\parallel
\un{\omega} \parallel} M^{\un{\omega}} .$$
On a donc la relation suivante
\begin{equation}
\label{symetrelconjugo}
(Ne^{\bullet} )^{-1} =(1^{\bullet} +I^{\bullet }
)\times (Ne^{\bullet } )^{-1} .
\end{equation}

\begin{lem}
Soit $F=F_{lin} (1+\di\sum_{\un{n} \in A(F)} B_{\un{n}} )$ un
automorphisme local de $\cC^{\nu }[[x]]$, $F_{lin} : \cC [[x]]
\rightarrow \cC [[x]]$, d\'efini par $F_{in} (\phi )=\phi
(\lambda_1 x_1 ,\dots ,\lambda_{\nu} x_{\nu} )$ pour tout $\phi
\in \cC^{\nu } [[x]]$. On note $\un{\lambda }=(\lambda_1 ,\dots
,\lambda_{\nu} ) \in \cC^{\nu }$. Pour toute suite $\un{s} =(s_1
,\dots ,s_r ) \in A^* (F)$, $s_i \in A(F)$, on note $\omega
(\un{s} ) =(\omega_1 ,\dots ,\omega_r )$ le vecteur de $\cC^{\nu}$
d\'efini par $\omega_i =s_i .\un{\lambda}$, et $\parallel \omega
(\un{s} )\parallel =\omega_1 +\dots +\omega_r$.\\

Si $F$ est non r\'esonnant, i.e.
\begin{equation}
\parallel \omega (\un{s} )\parallel \not= 0 \ \mbox{\rm pour tout}\ \un{s} \in A^* (F),
\end{equation}
alors, il existe un automorphisme de lin\'earisation formelle de
la forme
$$ Ne=\di\sum_{\bullet} Ne^{\bullet} B_{\bullet} ,$$
avec, pour tout $\un{s} =(s_1 , \dots ,s_r ) \in A^* (F)$, $s_i
\in A(F)$,
\begin{equation}
Ne^{\un{s}} =\di \di {e^{\omega_1 +\dots +\omega_r } \over
(e^{-\omega_1} -1) \dots (e^{-(\omega_1 +\dots +\omega_r)} -1)} ,
\end{equation}
o\`u $\omega_i =s_i .\un{\lambda}$.
\end{lem}

\begin{proof}
On calcule le moule $(Ne^{\bullet } )^{-1}$ par r\'ecurrence
sur la longueur des suites. \\

Pour $l(\omega )=1$, on a
$$Ne^{\omega } (Ne^{\omega })^{-1} =1+(Ne^{\omega } )^{-1} ,$$
d'o\`u
$$(Ne^{\omega} )^{-1} =\di {1\over e^{\omega} -1} .$$

Pour $l(\un{\omega} )=r$, $\un{\omega } =(\omega_1 ,\dots
,\omega_r )$, on a
$$e^{\parallel \un{\omega} \parallel } (Ne^{\un{\omega}} )^{-1} = (Ne^{\un{\omega}} )^{-1}
+(Ne^{\omega_2 ,\dots ,\omega_r } )^{-1} ,$$ soit
$$(Ne^{\un{\omega}} )^{-1} =\di
{1\over e^{\parallel \un{\omega} \parallel} -1} (Ne^{\omega_2
,\dots ,\omega_r } )^{-1} .$$
Une simple r\'ecurrence donne la formule suivante
\begin{equation}
(Ne^{\un{\omega}} )^{-1} =\di {1\over (e^{\omega_1 +\dots
+\omega_r} -1)(e^{\omega_2 +\dots +\omega_r} -1)\dots
(e^{\omega_r} -1)} .
\end{equation}
On d\'eduit alors de la relation
$$(Ne^{\bullet} )\times (Ne^{\bullet} )^{-1} =1^{\bullet } ,$$
la formule de $Ne^{\bullet}$ :
\begin{equation}
Ne^{\un{\omega}} =\di {e^{\omega_1 +\dots +\omega_r } \over
(e^{-\omega_1} -1) \dots (e^{-(\omega_1 +\dots +\omega_r)} -1)} .
\end{equation}
Comme $Ne^{\bullet}$ v\'erifie l'\'equation (\ref{symetrelconjugo}), on d\'eduit du th\'eor\`eme \ref{?} que le moule
$Ne^{\bullet}$ est symetrel, i.e. que $Ne$ est bien un automorphisme de $\cC [[x]]$.
\end{proof}

\subsection{Universalit\'e du moule de lin\'earisation}

On peut qualifier le moule $Na^{\bullet}$ ou $Ne^{\bullet}$ de coefficient {\it universel}, et
ceci pour au moins deux raisons :\\

-- deux champs de vecteurs de m\^eme partie lin\'eaire et de m\^eme alphabet ont {\it exactement}
le m\^eme moule de lin\'earisation.\\

-- tous les champs de vecteurs non r\'esonants se lin\'earisent via un moule dont
l'expression formelle est fixe.\\

Le calcul moulien permet donc d'isoler dans le probl\`eme de lin\'earisation, ce qui
est intrins\`eque, de ce qui ne l'est pas\footnote{Pour une utilisation de cette propri\'et\'e d'universalit\'e dans
un probl\`eme de bifurcation de champs de vecteurs, on renvoie \`a \cite{cs2} et \cite{Cr1}.}.

\section{Pr\'enormalisation}
\label{prenormalisation}

\subsection{Formes pr\'enormales}

On d\'efinit, en suivant Ecalle-Vallet \cite{EV}, la notion de
{\it forme pr\'enormale continue}.

\begin{defi}
Une forme pr\'enormale d'un champ $X$ dont la partie semi-simple
est diagonale est la donn\'ee d'un champ de vecteur $X_{\rm pran}$ de
la forme
\begin{equation}
X_{\rm pran} =X_{\rm lin} +X_r ,\ \ \mbox{\rm avec}\ \ [X_{\rm lin} ,X_r ]=0.
\end{equation}
\end{defi}

Le champ $X_r$ est donc uniquement constitu\'e de mon\^omes r\'esonnants.

\begin{rema}
La classique forme normale de {\it Poincar\'e-Dulac} est une forme pr\'enor\-male.
\end{rema}

\begin{defi}
Soit $X$ un champ de vecteur local de $\cC^{\nu}$, sous forme
pr\'epar\'ee. Une forme pr\'enormale est dite continue, si elle
est continue par rapport aux op\'erateurs $D_n$, $\un{n} \in
A(X)$, le spectre de $X$ \'etant fix\'e.
\end{defi}

Nous avons le r\'esultat suivant :

\begin{lem}
Soit $X$ un champ de vecteur local de $\cC^{\nu}$, de partie
lin\'eaire diagonale $X_{\rm lin}$, d'alphabet $A(X)$. Le champ
\begin{equation}
\label{pran} X_{\rm pran} =X_{\rm lin} +\di\sum_{\bullet \in A^* (X)}
\mbox{\rm Pran}^{\bullet } D_{\bullet } ,
\end{equation}
o\`u le moule $\mbox{\rm Pran}^{\bullet}$ est alternal et v\'erifie
\begin{equation}
\label{res}
\mbox{\rm Pran}^{\omega} =0\ \ \mbox{\rm si}\ \ \parallel
\omega \parallel \not=0 ,
\end{equation}
est une forme pr\'enormale continue.
\end{lem}

\begin{proof}
Le champ $X_{\rm pran}$ est une forme pr\'enormale car $X_{\rm pran} -X_{\rm lin}$
ne contient que des termes r\'esonnants via (\ref{res}). Par
construction, $X_{\rm pran}$ est evidemment continue par rapport aux
$D_{\un{n}}$, $\un{n} \in A(X)$, le spectre de $X$ \'etant fix\'e.
\end{proof}

\begin{rema}
Le moule $\mbox{\rm Pran}^{\bullet }$ d\'epend seulement de
l'alphabet $A(X)$. Autrement dit, si deux champs de vecteurs
locaux $X_1$ et $X_2$, de m\^eme partie lin\'eaire, d\'efinissent
le m\^eme alphabet, i.e. $A(X_1 )=A(X_2 )$, alors le moule
$\mbox{\rm Pran}^{\bullet}$ d\'efinissant la forme pr\'enormale
(\ref{pran}) est le m\^eme pour $X_1$ et $X_2$.
\end{rema}

\subsection{Forme normale de Poincar\'e-Dulac}

La forme normale de {\it Poincar\'e-Dulac} est construite par
it\'eration. On rappelle ici le proc\'ed\'e de construction:

On associe \`a $X$ un champ dit {\it simplifi\'e} de la forme
\begin{equation}
X_{sam} = \left ( \mbox{\rm exp} \sum_n^* \di {D_n \over n.\lambda }
\right ) X \left ( \mbox{\rm exp} \sum_n^* \di {D_n \over
n.\lambda } \right ) ,
\end{equation}
o\`u $\di\sum_n^*$ porte sur tous les multi-entiers $n=(n_1 ,\dots
,n_{\nu} )\in \Omega$ tels que $\omega \not=0$. C'est le
proc\'ed\'e classique de supression des termes non r\'esonnants du
champ.

On peut it\'erer ce processus et passer \`a la limite. On appelle
forme normale de Poincar\'e-Dulac le champ limite et on le note
$X_{tram}$. On a donc
\begin{equation}
X\rightarrow X_{sam}^1 \rightarrow X_{sam}^2 \rightarrow \dots
\rightarrow X_{tram}.
\end{equation}
Nous avons le th\'eor\`eme suivant:

\begin{thm}
La forme normale de Poincar\'e-Dulac est une forme pr\'enormale
continue, appel\'ee forme pr\'enormale \'elagu\'ee et not\'ee
$X_{tram}$:
\begin{equation}
X_{tram} =X_{lin} + \sum_{\bullet} Tram^{\bullet} D_{\bullet} .
\end{equation}
\end{thm}

La d\'emonstration repose sur une r\'eecriture de la construction
classique en terme de moules et comoules, faisant apparaitre ainsi
les constantes universelles $Tram^{\bullet }$.

\begin{rema}
i. La forme normale de Poincar\'e-Dulac n'est pas une forme {\it
normale} dans le sens de Baider, Sanders ou Ecalle. C'est pour
eviter toute confusion qu'Ecalle appelle cet objet forme
pr\'enormale \'elagu\'ee. Dans la suite de cet article, nous
utiliserons la terminologie
d'Ecalle.\\

ii. Les constantes universelles $Tram^{\bullet }$ sont explicites
(voir le lemme \ref{expli}).
\end{rema}

\begin{proof}
Elle repose en partie sur le lemme suivant:

\begin{lem}
Le champ simplifi\'e $X_{sam}$ associ\'e \`a $X$ poss\`ede un
d\'eveloppement moulien de la forme
\begin{equation}
\label{simp} X_{sam} =X_{lin} +\sum Sam^{\bullet} D_{\bullet } ,
\end{equation}
o\`u le moule $Sam^{\bullet }$ est d\'efini par

-- $Sam^{\emptyset} =0$

-- $Sam^{\omega} =0$ si $\omega\dot \not=0$ et $Sam^{\omega} =1$ si $\omega =0$.

-- si $l(\bomega )\geq
2$, et tous les $\omega_i$, $1\leq i\leq r$ sont non nuls, alors
\begin{equation}
Sam^{\bomega } =\di {1\over \omega_1 \dots \omega_r}\di\sum_{k=1}^r
\di {(-1)^{r-k} (\omega_k (r-k) -\omega_{k+1} -\dots -\omega_r ) \over (k-1)! (r-k+1)!} .
\end{equation}

-- Si un seul $\omega_i$ s'annule,
\begin{equation}
Sam^{\bomega } =\di {(-1)^{r-1} \over (i-1)! (r-i)! \omega_1 \dots \omega_{i-1} \omega_{i+1}
\dots \omega_r} .
\end{equation}

-- Enfin, si plus d'un $\omega_i$ s'annule, alors $Sam^{\bomega} =0$.
\end{lem}

\begin{rema}
Le champ simplifi\'e poss\`ede un d\'eveloppement moulien mais
n'est pas une forme pr\'enormale. Il contient des op\'erateurs non
r\'esonnants.
\end{rema}

\begin{proof}
On commence par noter que $\Theta =\di\sum_{n\in {\cal M} }^* {D_n
\over \omega }$ peut s'\'ecrire comme
\begin{equation}
\Theta =\sum_{n\in S_{\cal M}} J^{\bullet} D_{\bullet} ,
\end{equation}
o\`u le moule $J^{\bullet}$ est d\'efini par
$$J^{\bullet} =
\left \{
\begin{array}{c}
\di
{1\over \omega} \ \mbox{\rm si}\ \bullet=\omega\ \mbox{\rm et} \ \omega \not=0 ,\\
0\ \mbox{\rm sinon.}
\end{array}
\right .
$$
L'exponentielle de $\Theta$ poss\`ede un d\'eveloppement moulien
de la forme
$$\exp \Theta =\sum_{n\in S_{\cal M}} \exp (J^{\bullet} ) D_{\bullet } .$$
Par d\'efinition, on a
$$\exp (\Theta )^{\bomega} = (1^{\bullet} )^{\bomega} +\di {1\over 2!}
(J^{\bullet} \times J^{\bullet} )^{\bomega} + \dots .$$
Or, un simple calcul donne
$$
\underbrace{J^{\bullet} \times \dots \times J^{\bullet}}_{r facteurs} =
\left \{
\begin{array}{c}
\di {1\over \omega_1 \dots \omega_r} \ \mbox{\rm si}
\ \bomega=\omega_1 \dots \omega_r\ \mbox{\rm avec tous les}\ \omega_i \not= 0,\ i=1,\dots ,r ,\\
0\ \mbox{\rm sinon.}
\end{array}
\right .
$$
On a donc
$$
\exp \Theta^{\bullet} =
\left \{
\begin{array}{c}
1\ \mbox{\rm si}\ \bullet =\emptyset,\\
\di {1\over r! \omega_1 \dots \omega_r} \ \mbox{\rm si}\ \bullet =\omega_1 \dots \omega_r,\
\mbox{\rm avec}\ \mbox{\rm tous les}\ \omega_i \not= 0,\\
0\ \mbox{\rm sinon}.
\end{array}
\right .
$$
Par d\'efinition de $X_{sam }$, nous avons
$$X_{sam} =(\sum \exp (J^{\bullet } )B_{\bullet } ) X (\sum \exp (-J^{\bullet } )
B_{\bullet } ).$$ Or $X=X_{lin} +\di \sum_{n\in \Omega } D_n$, qui
poss\`ede l'\'ecriture moulienne
$$X=X_{lin} +\sum_{\bullet } I^{\bullet} D_{\bullet } ,
$$
o\`u $I^{\bullet }$ est le moule d\'efini par $I^{\omega_1 } =1$
et $I^{\omega } =0$ si $l(\omega )\geq 2$. On a donc l'expression
moulienne compl\`ete de $X_{sam }$ :
$$
X_{sam} = \left (\di\sum_{\bullet } \exp (J^{\bullet } )D_{\bullet
} \right ) \left ( \di X_{lin} +\sum I^{\bullet } D_{\bullet }
\right ) \left ( \di\sum_{\bullet } \exp (-J^{\bullet } )
D_{\bullet } \right ) .
$$
Comme
$$\left ( \di\sum_{\bullet } \exp (J^{\bullet } )D_{\bullet } \right ) (X_{lin})
\left ( \di\sum_{\bullet } \exp (-J^{\bullet } ) D_{\bullet }
\right ) = \di\sum_{\bullet } \exp (J^{\bullet } ) \times \td \exp
(-J^{\bullet } ) D_{\bullet } ,$$ avec $\td$ la d\'erivation sur
l'alg\`ebre des moules d\'efinie par
$$
\td M^{\omega } =\omega M^{\omega } ,
$$
on obtient l'\'equation moulienne suivante pour $Sam^{\bullet }$ :
\begin{equation}
Sam^{\bullet } =\exp J^{\bullet } \times \td\exp (-J^{\bullet } )
+\exp J^{\bullet } \times I^{\bullet } \times \exp (-J^{\bullet }
) .
\end{equation}
Cette \'equation permet de calculer l'expression du moule
$Sam^{\bullet }$ par r\'ecurrence sur la longueur des s\'equences.\\

Le moule $C^{\bullet } =\exp J^{\bullet } \times I^{\bullet }$ est
d\'efini par
$$C^{\bullet } =
\left \{
\begin{array}{c}
0\ \mbox{\rm si au moins un des $\omega_i$, $1\leq i\leq r-1$ est nul},\\
\di {1\over (r-1)! \omega_1 \dots \omega_{r-1}}\ \mbox{\rm sinon.}
\end{array}
\right .
$$
Le moule $D^{\bullet } =C^{\bullet } \times \exp (-J^{\bullet })$
est donc d\'efini par $D^{\emptyset } =0$, $D^{\omega } =1$. Pour une suite $\bomega$ de
longueur $r\geq 2$, on a
$$D^{\omega_1 \dots \omega_r} =C^{\omega_1} (\exp (-J^{\bullet} ))^{\omega_2 \dots \omega_r}
+C^{\omega_1 \omega_2} (\exp (-J^{\bullet} ))^{\omega_3 \dots \omega_r} +\dots +C^{\omega_1 \dots
\omega_r} .$$
On a plusieurs cas :

-- si au moins un des $\omega_i$ est nul, $1\leq i\leq r-1$, alors tous les $C^{\omega_1 \dots \omega_j}$,
avec $j\geq i+1$ son nuls, ainsi que tous les $(\exp (-J^{\bullet} ))^{\omega_k \dots \omega_r}$
pour $k\leq i$. Donc, on a
$$D^{\omega_1 \dots \omega_r} =C^{\omega_1 \dots \omega_i} (\exp (-J^{\bullet} )^{\omega_{i+1}
\dots \omega_r} .$$ On a donc :
$$
D^{\omega_1 \dots \omega_r} =
\left \{
\begin{array}{c}
0\ \mbox{\rm si au moins deux $\omega_i$ sont nuls.}\\
\di {(-1)^{r-i}\over (i-1)!(r-i)! [ \omega_1 \dots \omega_{i-1}][\omega_{i+1}
\dots \omega_r ]}\ \mbox{\rm si un seul $\omega_i$, $1\leq i\leq r-1$, est nul.}
\end{array}
\right .
$$

-- si seul $\omega_r$ est nul, alors $D^{\omega_1 \dots \omega_r} =C^{\omega_1 \dots \omega_r}$
et on a
$$D^{\omega_1 \dots \omega_r} =\di {1\over (r-1)! \omega_1 \dots \omega_{r-1}} .$$

-- si aucun $\omega_i$ n'est nul alors
$$
\left .
\begin{array}{ll}
D^{\omega_1 \dots\omega_r} =\di {(-1)^{r-1} \over (r-k)! \omega_2 \dots \omega_r} & +
\di {1\over \omega_1} \times \di {(-1)^{r-2}\over (r-2)! \omega_3 \dots \omega_r} +\dots  \\
 & +\di {1\over (r-2)! \omega_1 \dots \omega_{r-2}} \times \di {-1\over \omega_r}
 +\di {1\over (r-1)! \omega_1 \dots \omega_{r-1}} ,
\end{array}
\right .
$$
soit
$$D^{\omega_1 \dots \omega_r} =
\di {1\over \omega_1 \dots \omega_r} \di\sum_{k=1}^r \di {(-1)^{r-k} \omega_k \over (k-1)! (r-k)!} .$$

De m\^eme, on montre facilement que $E^{\bullet } = \td \exp
J^{\bullet }$ est d\'efini par $E^{\emptyset } =0$, et plus g\'en\'eralement
$$E^{\omega_1 \dots \omega_r } =
\left \{
\begin{array}{c}
0\ \mbox{\rm si au moins un des $\omega_i$, $1\leq i\leq r$, et nul;}\\
\di{(\omega_1 +\dots +\omega_r)(-1)^r  \over r! \omega_1 \dots \omega_r } .
\end{array}
\right .
$$
On en d\'eduit l'expression du moule $F^{\bullet } =\exp(J^{\bullet })\times E^{\bullet}$. On a
$F^{\emptyset } =0$ et $F^{\omega}
=E^{\omega}$ donc $F^{\omega} =0$ si $\omega =0$ et $F^{\omega} =-1$ sinon; pour un mot de
longueur $r\geq 1$, on a :
$$F^{\omega_1 \dots \omega_r} =(\exp (J^{\bullet} ))^{\emptyset} E^{\omega_1 \dots \omega_r}
+(\exp (J^{\bullet} ))^{\omega_1} E^{\omega_2 \dots \omega_r} +\dots + (\exp (J^{\bullet} ))^{\omega_1
\dots \omega_{r-1}} E^{\omega_r} ,$$
donc $F^{\omega_1 \dots \omega_r} =0$ si un au moins des $\omega_i$ est nul. Si aucun
$\omega_i$ n'est nul, alors
$$
\left .
\begin{array}{rl}
F^{\omega_1 \dots \omega_r} = & \di {(\omega_1 +\dots +\omega_r )(-1)^r \over r! \omega_1 \dots \omega_r} +
\dots +\di {(-1)\omega_r \over (r-1)! \omega_1 \dots \omega_r},\\
 = & \di {1\over \omega_1 \dots \omega_r}\di\sum_{k=1}^r \di {(-1)^{r-k+1} (\omega_k +
 \dots +\omega_r )\over (r-k+1)! (k-1)!} .
\end{array}
\right .
$$

$Sam^{\bullet} =F^{\bullet }+D^{\bullet}$, on a
$$Sam^{\emptyset} =F^{\emptyset}+D^{\emptyset} =0,$$
et
$$Sam^{\omega} =\di \left \{
\begin{array}{l}
1\ \mbox{\rm si $\omega =0$.},\\
0\ \mbox{\rm si $\omega\not= 0$.}
\end{array}
\right .
$$
Pour toute suite de longueur $r\geq 2$, on a:\\

-- si au moins deux $\omega_i$ sont nuls, alors $F^{\bomega} =D^{\bomega} =0$, donc
$Sam^{\bomega} =0$.

-- si un seul $\omega_i$ est nul, $F^{\bomega} =0$ et
$$Sam^{\bomega} =\di {1\over (i-1)! \omega_1 \dots \omega_{i-1}}
\di {(-1)^{r-1} \over (r-i)! \omega_{i+1} \dots \omega_r} .$$

-- si tous les $\omega_i$, $1\leq i\leq r$, sont non nuls, alors
$$Sam^{\omega_1 \dots \omega_r} =
\di {1\over \omega_1 \dots \omega_r}
\di\sum_{k=1}^r \di {(-1)^{r-k} (\omega_k (r-k) -\omega_{k+1} -\dots -\omega_r )\over
(k-1)! (r-k+1)!} .$$
On en d\'eduit le lemme.
\end{proof}

Pour terminer la d\'emonstration du th\'eor\`eme, il suffit de remarquer que par it\'eration les champs $X_{sam}^r$
garde une forme de type (\ref{simp}), de m\^eme pour $X_{tram}$:
\begin{equation}
X_{tram} =X_{lin} + \sum_{\bullet} Tram^{\bullet} D_{\bullet} .
\end{equation}
On en d\'eduit le th\'eor\`eme.
\end{proof}

Les moules de la forme \'elagu\'ee ont une expression alg\'ebrique simple:

\begin{lem}
\label{expli} Le moule de la forme \'elagu\'ee est d\'efinie par
\begin{equation}
Tram^{\omega } =(Sam ^{\bullet } )^{\circ l(\omega )},
\end{equation}
o\`u $(Sam ^{\bullet} )^{\circ n} =Sam^{\bullet } \circ \dots
\circ Sam^{\bullet }$, $n$ fois.
\end{lem}

\begin{rema}
Le lemme pr\'ec\'edent traduit simplement le fait qu'\`a la $r$-i\`eme \'etape, on ne touche pas les composantes
de degr\'e inf\'erieur \`a $r$ du champ.
\end{rema}

Le calcul explicite des moules $Tram^{\bullet }$ est possible car on a une expression exacte des moules du
champ simplifi\'e. Il n'est pas n\'ecessaire d'it\'erer la composition du moule $Sam$. En effet, nous avons les deux relations
suivantes:
\begin{eqnarray}
Tram^{\bullet} =Tram^{\bullet} \circ Sam^{\bullet} ,\\
Tram^{\bullet} =Sam^{\bullet} \circ Tram^{\bullet},
\end{eqnarray}
qui proviennent toutes les deux de la stabilisation des compos\'es it\'er\'es du moule $Sam^{\bullet}$ sur $Tram^{\bullet}$.\\

De ces deux formules, seule la premi\`ere fournit une relation de r\'ecurrence. En effet, on a:
\begin{eqnarray}
Tram^{\omega_1 \dots \omega_n} =\di\sum_{u^1 \dots u^r =\omega} Tram^{\parallel u^1 \parallel \dots \parallel u^r\parallel}
Sam^{u^1} \dots Sam^{u^r},\\
Tram^{\omega_1 \dots \omega_n} =\di\sum_{u^1 \dots u^r =\omega} Sam^{\parallel u^1 \parallel \dots \parallel u^r\parallel}
Tram^{u^1} \dots Tram^{u^r}.
\end{eqnarray}

Si au moins un $\omega_i \not= 0$, alors le moule $Tram^{\omega_1 \dots \omega_r}$ est d\'efini par
\begin{equation}
Tram^{\omega_1 \dots \omega_n} =\di\sum_{u^1 \dots u^r =\omega ;\ r\leq n-1} Tram^{\parallel u^1 \parallel \dots \parallel u^r\parallel}
Sam^{u^1} \dots Sam^{u^r} .
\end{equation}
Par ailleurs si tous les $\omega_i$ sont nuls, on a par alternalit\'e pour toute suite de longueur $\geq 2$:
\begin{equation}
Tram^{0\dots 0} =0 .
\end{equation}

\begin{rema}
Le moule $Tram$ est associ\'e \`a une forme pr\'enormale. Il est donc nul sur toute suite $\omega$ telle que $\parallel \omega \parallel
\not= 0$, i.e. sur les suites non r\'esonantes. Malheureusement, si $\parallel \omega\parallel =0$, alors le moule $Tram^{\omega}$
disparait dans la seconde \'equation.
\end{rema}

\section{Arborification : une introduction}

Un probl\`eme important dans la normalisation des champs de vecteurs ou diff\'eomor\-phis\-mes, est celui de la
convergence/divergence de l'automorphisme de conjugaison. Ce probl\`eme est \`a la source des travaux de Siegel, Bruyno,
Russman, Arnold, Moser et Yoccoz. La m\'ethode d'arborification/coarborification permet de d\'emontrer la
convergence (lorsque c'est le cas) des s\'eries formelles construites via le calcul moulien.\\

Pourquoi a-t-on besoin d'arborifier/corarborifier la s\'erie pour \'etablir sa convergence ?\\

Les s\'eries obtenues par le calcul moulien se prettent mal \`a l'analyse. Une semi-norme sur les op\'erateurs
de $\cC\{ x\}$ \'etant donn\'ee, on obtient, en g\'en\'eral, de tr\`es mauvaises
estimations sur la semi-norme de la s\'erie. Cet artefact est d\^u \`a la majoration directe d'op\'erateurs de
la forme $B_1 \dots B_n$ qui ne sont pas homog\`enes. Une id\'ee est donc de
re\'ecrire la s\'erie en faisant appara\^\i tre des op\'erateurs
homog\`enes. Le codage de ce proc\'ed\'e est la m\'ethode d'{\it arborification}.\\

Commen\c{c}ons par d\'efinir une norme:

\begin{defi}
Soient $U$ et $V$ deux voisinages compacts de $0$ dans $\cC^{\nu}$, tels que $V\subset U$. Pour tout germe de fonction
$\phi$ de $\cC \{ x\}$ en $0$, on d\'efinit
\begin{equation}
\parallel \phi \parallel_U =\sup_{x\in U} \mid \phi (x)\mid .
\end{equation}
De m\^eme, pour tout op\'erateur $P$ de $\cC \{ x\}$ dans lui-m\^eme, on d\'efinit
\begin{equation}
\parallel P\parallel_{U,V} =\sup_{\phi \mid \, \parallel \phi\parallel_U \leq 1}
\parallel P.\phi \parallel_V .
\end{equation}
\end{defi}

On dit que la s\'erie d'op\'erateurs $\di\sum_{\bn} P_{\bn}$ est {\it normalement convergente} si la famille
$(\parallel P_{\bn} \parallel )_{\bn}$ est sommable pour une paire $(U,V)$ au moins.

\subsection{Convergence sans arborification : le th\'eor\`eme de Poincar\'e}

La d\'emons\-tra\-tion de la convergence normale des s\'eries mouliennes ne n\'ecessite pas toujours le recours \`a
la m\'ethode d'arborifcation. C'est le cas des s\'eries du th\'eor\`eme de lin\'earisation de Poincar\'e. Vue la
forme du moule intervenant dans le probl\`eme de lin\'earisation, cette convergence ne peut avoir lieu que sous une
contrainte sur la vitesse \`a laquelle les $\omega (\bn )$ peuvent s'approcher de $0$ lorsque $\parallel \bn
\parallel$ augmente. Cette ``vitesse" d\'epend essentiellement de la disposition des valeurs propres du spectre de
la partie lin\'eaire du champ.

\begin{defi}
Soit $\blambda =(\lambda_1 ,\dots ,\lambda_{\nu} )$ une collection de valeurs propres dans $\cC^{\nu}$.
On d\'efinit:\\

-- le domaine de Poincar\'e ${\cal P}$ comme l'ensemble des $\blambda$ dont l'enveloppe convexe ne
contient pas $0$.\\

-- le domaine de Siegel ${\cal S}$, comme  le compl\'ementaire du pr\'ec\'edent.
\end{defi}

Une condition de contr\^ole tr\`es forte du spectre est l'{\it absence de petits diviseurs}.

\begin{defi}
On dit que le champ $X$ ne contient pas de petits diviseurs s'il existe une constante $C>0$, telle que
$$\forall \omega \in \Omega ,\ \mid \omega \mid \geq C .$$
\end{defi}

Si le champ ne poss\`ede pas de petits diviseurs, il est n\'ecessairement non r\'esonant, donc formellement
lin\'earisable d'apr\`es le th\'eor\`eme de Poincar\'e (voir th\'eor\`eme \ref{thmpoinc}). La convergence de la
s\'erie normalisante est assur\'ee par le th\'eor\`eme suivant:

\begin{thm}
Soit $X$ un champ de vecteur ne contenant pas de petits diviseurs,
et dont le spectre est dans le domaine de Poincar\'e. Alors,
l'automorphisme de lin\'earisation du th\'eor\`eme de Poincar\'e
(th\'eor\`eme \ref{thmpoinc}) est analytique.
\end{thm}

La d\'emonstration suit la d\'emarche usuelle, mais sur les
moules. Elle va nous permettre de d\'egager des majorations
importantes pour la suite.

\begin{proof}
La premi\`ere \'etape consiste \`a se ramener, via un changement de variables analytique (en fait, m\^eme alg\'ebrique)
\`a une situation o\`u l'alphabet est tel que les $\omega$ sont toutes de partie r\'eelle positive. C'est finalement,
uniquement ce fait qui assure la convergence de la normalisation\footnote{Nous allons le voir, ceci implique que les
combinaisons lin\'eaires des $\omega \in \Omega$, intervenant dans le moule de lin\'earisation, ne peuvent pas \^etre
trop petites.}.\\

En effet, comme les valeurs propres $\lambda_i$ sont dans $\cal P$, il existe un $\theta \in \rR$ tel que
$$\lambda_i \in P_{\theta} =\{ z\in \cC , \mbox{\rm Re} (ze^{i\theta} ) >0 \} .$$
On peut donc, via une rotation, se ramener au cas o\`u toutes les valeurs propres
sont dans le demi-plan $\{ \mbox{\rm Re} (z) >0 \}$. Il existe donc une constante
$\rho >0$ telle que pour tout $i$, $\mbox{\rm Re} (\lambda_i ) >\rho $. Une cons\'equence
importante est qu'il n'existe qu'un nombre fini de $\omega$ dans $\Omega$ tels que
$\mbox{\rm Re}(\omega ) <0$. Il existe donc un changement de variable polynomial tel que
$X$ s'\'ecrive
$$X=X_{\rm lin} +\di\sum_{n \in P(X)} B_n ,$$
o\`u $P(X)$ est le nouvel alphabet tel que
$$\forall n\in P(X),\ \mbox{\rm Re} (\omega(n)) >0 .$$
Nous avons le lemme suivant :

\begin{lem}
\label{inega1}
Pour toute suite $\bn$ de longueur $r$, on a
\begin{equation}
\label{inegamou}
\mid Na^{\bn} \mid \leq \di {1\over r! C_1^r } ,\
\mid (Na^{-1})^{\bn} \mid \leq \di {1\over r! C_1^r } ,\
C_1 >0 ,
\end{equation}
\begin{equation}
\label{inegaop}
\parallel B_{\bn} \parallel_{U,V} \leq r! C_2^{N(\un{n})} \parallel B_{n^1} \parallel_{U,V}
\dots \parallel B_{n^r} \parallel_{U,V} ,\ C_2 >0 ,
\end{equation}
\begin{equation}
\label{inegaindi}
\parallel B_n \parallel_{U,V} \leq (C_{U,V} )^{\mid n\mid} ,\ C_{U,V} >0 ,
\end{equation}
\begin{equation}
\label{ineganom}
c(N)\leq (C_3 )^{N(\bn )},\ C_3 >0 ,
\end{equation}
o\`u $\mid n\mid =n_1 +\dots +n_{\nu}$ et $N(\bn )=\mid n^1 \mid +\dots +\mid n^r \mid$ est
le poids de $\bn$; $c(N)$ est le nombre de mots de poids $N$.
\end{lem}

La d\'emonstration est donn\'ee dans la prochaine section.\\

On d\'emontre la convergence normale du normalisateur $Na$. On a
$$
\left .
\begin{array}{lll}
\parallel Na^{\bn} B_{\bn} \parallel  & \leq & \di {1\over r! C_1^r} r! \parallel B_{n^1}
\parallel_{U,V} \dots \parallel B_{n^r} \parallel_{U,V} C_2^{N} ,\\
 & \leq & \di {1\over C_1^r} (C_{U,V} )^{rN} C_2^N .
\end{array}
\right .
$$
Par cons\'equent, on a
$$
\left .
\begin{array}{lll}
\left \|
\di\sum_{l(\bn )=r ,\ N(\bn )=N } Na^{\bn} B_{\bn}
\right \| _{U,V} & \leq & \di\sum_{l(\bn )=r,\ N(\bn )=N} \parallel Na^{\bn}
B_{\bn} \parallel_{U,V} , \\
 & \leq & \di {c(N)\over C_1^r } C_{U,V}^{rN} C_2^N ,\\
 & \leq & (C_3 C_2 )^N \di {(C_{U,V}^r )^N \over C_1^r} .
\end{array}
\right .
$$
Comme $r\leq N$, on a
$$\di {1\over C_1^r} \leq \di {1\over C_1^N}.$$
On en d\'eduit
$$
\left \|
\di\sum_{l(\bn )=r ,\ N(\bn )=N } Na^{\bn} B_{\bn}
\right \| _{U,V} \leq \left ( \di {C_3 C_2 C_{U,V}^r \over C_1} \right ) ^N  .
$$
Pour un bon choix de $(U,V)$, on peut rendre $C_{U,V}$ aussi petit que l'on veut. La
s\'erie est donc normalement convergente.
\end{proof}

\subsubsection{D\'emonstration du lemme \ref{inega1}}

Par hypoth\`ese, tous les $\omega (n)$ sont de parties r\'eelles strictement positives. Par ailleurs,
l'absence de petits diviseurs impose que pour tout $\omega$, $\mid \omega \mid \geq C >0$. On a
donc
$$\forall \omega \in P(X),\ \mbox{\rm Re}(\omega )\geq C .$$
On en d\'eduit pour tout $r$,
$$
\left .
\begin{array}{lll}
\mid \omega_1 +\dots +\omega_r \mid & \geq & \mbox{\rm Re}(\omega_1 +\dots +\omega_r ) ,\\
 & \geq & r C.
\end{array}
\right .
$$
On a donc
$$
\mid Na^{\bn} \mid \leq \di {1\over r! C^r} .
$$

\subsection{Le probl\`eme des petits diviseurs et le th\'eor\`eme de Bruyno}

L'existence des petits diviseurs conduit \`a des difficult\'es analytiques s\'erieuses. Nous avons
le th\'eor\`eme de Siegel :

\begin{thm}
\label{siegel}
En pr\'esence de petits diviseurs, les s\'eries mouliennes de $\Theta$ et $\Theta^{-1}$ sont
g\'en\'eriquement divergentes.
\end{thm}

N\'eanmoins, on imagine bien qu'un contr\^ole de la vitesse de convergence des $\omega (\bn )$ vers
$0$ lorsque $\bn $ augmente doit permettre de r\'etablir la convergence de la s\'erie. C'est
effectivement ce qui se passe, sous une condition, appel\'ee {\it condition diophantienne de Bruyno}.\\

On note $\omega (k )$ la quantit\'ee d\'efinie par
\begin{equation}
\omega (k)=\inf \left \{
\blambda .m =\di\sum_{i=1}^{\nu} \lambda_i m_i ,\ \mid m\mid \leq 2^{k+1}\ \mbox{\rm et}\
\blambda .m \not= 0 \right \} ,
\end{equation}
o\`u les $m_i$ sont tous positifs, sauf au plus un qui peut valoir $-1$, de somme $\mid m \mid
=\di\sum_{i=1}^{\nu} m_i \geq 0$.\\

La condition diophantienne de Bruyno est :\\
$$
\mbox{\rm La s\'erie}\ S=\di\sum_{k=0}^{\infty} \di {\log (1/\omega (k))\over 2^k }\ \mbox{\rm
est convergente.}
\eqno{(B)}
$$

Le th\'eor\`eme de Bruyno s'\'enonce alors comme suit :

\begin{thm}[Bruyno]
Soit $X$ un champ de vecteurs analytique dont le spectre de la partie lin\'eaire v\'erifie la
condition (B). Alors, ce champ est analytiquement lin\'earisable.
\end{thm}

La d\'emonstration originale de Bruyno est longue et difficile. Une pr\'esentation claire et soign\'ee
de son travail est donn\'ee par Martinet \cite{mar}.\\

Ce que peut apporter le calcul moulien dans ce probl\`eme, c'est un cadre conceptuel clair qui guide
les diff\'erents calculs et estimations n\'ecessaires \`a la d\'emonstration.\\

En tout premier lieu, il faut comprendre pourquoi ce probl\`eme est beaucoup plus difficile que le
th\'eor\`eme de convergence sous la condition de Poincar\'e.\\

Les majorations du lemme \ref{inega1} concernant les op\'erateurs $B_{\bn}$ sont les m\^emes pour
le th\'eor\`eme de Bruyno, car elles ne d\'ependent pas du spectre de la partie lin\'eaire. Le
seul changement est dans l'estimation de la taille du moule de lin\'earisation $Na^{\bullet}$,
qui d\'epend fortement des propri\'et\'es arithm\'etiques des $\omega$.

\begin{lem}
Pour toute suite $\bn$, on a
\begin{equation}
\mid Na^{\bn} \mid \leq Q^{N(\bn )} ,\  Q>0 .
\end{equation}
\end{lem}

Une estimation directe de la norme de $Na$ donne, en gardant les in\'egalit\'es (\ref{inegaop}),
(\ref{inegaindi}), (\ref{ineganom}) du lemme \ref{inega1},
$$
\left \|
\di\sum_{l(\bn )=r ,\ N(\bn )=N } Na^{\bn} B_{\bn}
\right \| _{U,V} \leq \left ( \di r! {Q C_3 C_2 C_{U,V}^r } \right ) ^N  ,
$$
ce qui ne permet pas de conclure quand \`a la convergence de la s\'erie.\\

Les majorations pr\'ec\'edentes ne peuvent pas \^etre am\'eliorer. Autrement dit, le
th\'eor\`eme de Bruyno est inaccessible via l'expression moulienne initiale du normalisateur
$Na$. Si convergence il y a, il faut affiner l'analyse de cette s\'erie. \\

La m\'ethode d'{\it arborification} donne un proc\'ed\'e alg\'ebrique pour \'etudier la convergence de
ces s\'eries. N\'eanmoins, et c'est un point important, cette m\'ethode {\it prend sa source dans un
probl\`eme d'analyse}, et sa mise au point n'est pas un probl\`eme alg\'ebrique\footnote{Du fait
de l'utilisation d'arbres et autres suites arborifi\'ees, les sp\'ecialistes des alg\`ebres de
Hopf y ont souvent vu un relicat des bases de Hall. Or, l'arborification n'a strictement rien
\`a voir avec ce probl\`eme, qui lui, est purement alg\'ebrique. Bien entendu, une fois la
m\'ethode formalis\'ee, rien n'emp\'eche une \'etude purement alg\'ebrique de ses propri\'et\'es,
comme dans \cite{ev2} par exemple.}. Essentiellement, la question
\`a laquelle nous allons r\'epondre est la suivante : comment affiner notre estimation de la
norme de la s\'erie moulienne $Na$ ?

\subsection{Premier pas : homog\'en\'eit\'e et sym\'etrie}

La construction du normalisateur $Na$ utilise des op\'erateurs diff\'erentiels homog\`enes en
degr\'e, $B_n$, $n\in A(X)$. Les diff\'erentes combinaisons $B_{\bn}$ de ces op\'erateurs sont
encore des op\'erateurs diff\'erentiels homog\`enes de degr\'e $\parallel \bn\parallel$. Peut-on
pr\'eciser ce point ? Il suffit de faire un calcul en longueur $2$ pour avoir une id\'ee du
r\'esultat g\'en\'eral. \\

Soit $B_1$ et $B_2$ deux op\'erateurs de degr\'e $n^1$ et $n^2$ respectivement de la forme
$$B_i =\di\sum_{j=1}^{\nu} B_i (x_j )\partial_{x_j} ,\ i=1,2 .$$
L'op\'erateur compos\'e $B_{(1,2)}=B_1 B_2$ s'\'ecrit
\begin{equation}
B_{(1,2)} =\di\sum_{i,j} B_1 (x_i )  \partial_{x_i} [ B_2 (x_j ) ] \partial_{x_j} +
\di\sum_{i ,j} B_1 (x_i ) B_2 (x_j ) \partial^2_{x_i x_j} .
\end{equation}
De la m\^eme mani\`ere, on a :
\begin{equation}
B_{(2,1)} =\di\sum_{i,j} B_2 (x_i )  \partial_{x_i} [ B_1 (x_j ) ] \partial_{x_j} +
\di\sum_{i ,j} B_1 (x_i ) B_2 (x_j ) \partial^2_{x_i x_j} .
\end{equation}
On voit que le second terme du d\'eveloppement de $B_{(2,1)}$ est le m\^eme que
celui de $B_{(1,2)}$. Devant ce terme, le coefficient est $Na^{(1,2)} +Na^{(2,1)}$. Comme
le moule $Na^{\bullet}$ est sym\'etral, on a $Na^{(1,2)} +Na^{(2,1)}=Na^1 Na^2$.\\

Que faut-il retenir de ce calcul ?\\

On peut d\'ecomposer les $B_{\bn}$ de fa\c{c}on \`a profiter des sym\'etries du moule $Na^{\bullet}$
pour obtenir des majorations plus fines.\\

Avant de pr\'eciser la d\'ecomposition, on introduit la notion d'{\it ordre} d'un op\'erateur
diff\'erentiel.

\begin{defi}
Pour tout $\nu$-uplet $\delta =(\delta_1 ,\dots ,\delta_{\nu} )$ d'entiers positifs, on note
$\partial_x^{\delta} =\partial_{x_1}^{\delta_1} \dots \partial_{x_{\nu}}^{\delta_{\nu}}$. On
dit que $\partial_x^{\delta}$ est d'ordre $\mid \delta \mid =\delta_1 +\dots +\delta_{\nu}$.
Un op\'erateur diff\'erentiel est homog\`ene en ordre, d'ordre $d$, si il est de la forme
$$
B=\di\sum_{\delta ,\ \mid \delta \mid =d} b_{\delta} (x) \partial_x^{\delta} .
$$
\end{defi}

Quelle-est la d\'ecomposition ad\'equate des op\'erateurs $B_{\bn}$ ? On peut imposer deux contraintes
naturelles dans le cadre de ce probl\`eme :\\

i) On isole dans les $B_{\bn}$ des op\'erateurs diff\'erentiels homog\`enes en {\it ordre},
ceux-ci \'etant d\'eja homog\`enes en degr\'e. \\

ii) Les coefficients de ces op\'erateurs sont des produits de termes d\'ependants des $B_{n^i}$.\\

La condition ii) est n\'ecessaire si l'on veut facilement majorer la norme de ces op\'erateurs, ce
qui est essentiel pour notre propos.

\subsection{Codage du proc\'ed\'e de d\'ecomposition}

Le proc\'ed\'e de d\'ecomposition peut se coder via deux op\'erations sur les op\'erateurs
homog\`enes en ordre et en degr\'e, et donne naissance, une fois formalis\'e, \`a la
m\'ethode d'arborification.

\begin{defi}
\label{opelem}
Soient $B_1$ et $B_2$ deux op\'erateurs diff\'erentiels homog\`enes en degr\'e et en ordre,
$$B_i =\di\sum_{\delta ,\ \mid \delta\mid =d_i } b_{i,\delta} (x) \partial_x^{\delta} ,\ i=1,2,$$
o\`u les $b_{i,\delta}$ sont dans $\cC [[x]]$. On note
$$B_{(12)^{<}} \di\sum_{\delta ,\gamma }  b_{1,\delta} (x) \partial_x^{\delta} [b_{2,\gamma} (x)]
\partial_x^{\gamma} ,$$
et
$$B_{1\oplus 2} =\di\sum_{\delta ,\gamma} b_1^{\delta} (x) b_2^{\gamma} (x)
\partial_x^{\delta +\gamma } .$$
\end{defi}

Par r\'ecurrence sur la longueur des suites, on peut coder l'ensemble des parties homog\`enes
en degr\'e d'un op\'erateur $B_{\bn}$.\\

Notons $B_i =\di\sum_{i=1}^{\nu} b_i (x) \partial_{x_i}$,
$i=1,2,3$. On a
$$
\left .
\begin{array}{lll}
B_{1,2,3} & = & B_1 \left ( \di\sum_{j,k} B_2 (x_j ) \partial_{x_j} [B_3 (x_k )] \partial_{x_k} +
\sum_{j,k} B_2 (x_j ) B_3 (x_k ) \partial^2_{x_j x_k} \right ) ,\\
 & = & \di\sum_{i,j,k} B_1 (x_i) \left ( \partial_{x_i} [B_2 (x_j ) ]\partial_{x_j} [B_3 (x_k )]
 \partial_{x_k} \right . \\
 & & + B_2 (x_j ) \partial^2 _{x_i x_j} [B_3 (x_k )] \partial_{x_k} \\
 & & + B_2 (x_j ) \partial_{x_j} [B_3 (x_k )] \partial^2_{x_i x_k} \\
 & & + \partial_{x_i} [B_2 (x_j )] B_3 (x_k ) \partial^2_{x_j x_k} \\
 & & + B_2 (x_j ) \partial_{x_i} [B_3 (x_k ) ] \partial^2_{x_j x_k} \\
 & & \left .+ B_2 (x_j ) B_3 (x_k ) \partial^3_{x_i x_j x_k} \right ) .
\end{array}
\right .
$$
On a deux op\'erateurs diff\'erentiels d'ordre $1$, trois d'ordre $2$ et un d'ordre $3$. Soit, en
reprenant les notations de la d\'efinition \ref{opelem} :
$$B_{1,2,3} =B_{(1,2,3)^{<}} +B_{(1\oplus 2) 3)^{<}} +B_{(2,3)^{<} \oplus 1} +
B_{(1,2)^{<} \oplus 3} +B_{(1,3)^{<}\oplus 2} +B_{1\oplus 2\oplus 3} .$$

La mise en forme de cette d\'ecomposition n\'ecessite un alphabet nouveau, dit {\it arborifi\'e}.

\begin{defi}
Une appelle suite arborescente, et on note $\bomega^{<}$, une suite construite
sur $A(X)$ avec les symboles $<$ et $\oplus$. On note $\mbox{\rm Arb} (X)$ l'ensemble des
suites arborescentes.
\end{defi}

Si $\bn \in A(X)^*$, on note $\mbox{\rm Arb}(\bn )$ l'ensemble des $\ba^{<} \in \mbox{\rm Arb} (X)$
intervenant dans la d\'ecomposition de $B_{\bn}$. On a donc
$$B_{\bn } =\di\sum_{\ba^{<} \in \mbox{\rm Arb} (\bn )} B_{\ba^{<}} .$$
On peut donc r\'eecrire la s\'erie moulienne de $Na$ sous la forme
$$Na =\di\sum_{\bn \in A^* (X)} Na^{\bn} B_{\bn} =\di\sum_{\ba^{<} \in \mbox{\rm Arb} (X)}
Na^{\ba^{<}} B_{\ba^{<}} .$$
Les coefficients $Na^{\ba^{<}}$ d\'efinissent un moule sur $\mbox{\rm Arb} (X)$, not\'e
$Na^{{\bullet}^{<}}$.

\begin{defi}
Le moule $Na^{{\bullet}^{<}}$ est appel\'e moule arborifi\'e de $Na^{\bullet}$.
\end{defi}

Le moule arborifi\'e est bien entendu une combinaison lin\'eaire du moule initial. Si on note
$X(\ba^{<} )$ l'ensemble des suites $\bn$ de $A(X)^*$ telles que $\ba^{<} \in \mbox{\rm Arb} (\bn )$,
on a :
$$Na^{\ba^{<}} =\di\sum_{\bn \in X(\ba^{<} )} Na^{\bn} .$$
Il nous reste \`a formaliser cette construction, afin de simplifier sa pr\'esentation.

\subsection{D\'efinition formelle de l'arborification}

Les d\'efinitions pr\'ec\'eden\-tes sont li\'ees au codage de la construction des op\'erateurs diff\'erentiels
homog\`enes en ordre et degr\'e. On peut evidemment donner une d\'efinition purement alg\'ebrique de
ces objets, comme par exemple les suites arborescentes.

\begin{defi}
Une suite arborescente sur $A(X)$ est une suite $\bn^{<}$ d'\'el\'ements de $A(X)$, avec sur les indices
un ordre partiel appel\'e ordre arborescent : chaque $i$ de $\{ 1,\dots ,r\}$ poss\`ede au plus un
cons\'equent not\'e $i_+$.

On note $\bn_1^{<} \oplus \bn_2^{<}$ l'union disjointe de
$\bn_1^{<}$ et $\bn_2^{<}$; dans cette suite l'ordre partiel interne est conserv\'e, mais les
\'el\'ements de $\bn_1^{<}$ et $\bn_2^{<}$ sont incomparables.

Un $\bn^{<}$ est dit irr\'eductible s'il ne poss\`ede pas de d\'ecomposition non triviale
$\bn_1^{<} \oplus \bn_2^{<}$, autrement dit s'il poss\`ede un plus grand \'el\'ement.
\end{defi}

On peut toujours repr\'esenter une suite arborescente par un arbre.\\

Les op\'erateurs arborifi\'es peuvent alors se d\'efinir directement par
r\'ecurrence.

\begin{defi}
Pour une suite arborescente donn\'ee $\bn^{<}=(n^1 ,\dots ,n^r )^{<}$, on d\'efinit
$B_{\bn^{<}}$ comme \'etant l'unique op\'erateur v\'erifiant les trois propri\'et\'es
suivantes :\\

-- $B_{\bn^{<}} (\phi \psi )=\di\sum_{\bn_1^{<} \oplus \bn_2^{<} =\bn^{<}} \left (
B_{\bn_1^{<}} \phi \right ) \left ( B_{\bn_2^{<}} \psi \right ) $.\\

-- Si la suite $\bn^{<}$ se d\'ecompose en eaxctement $d$ suites irr\'eductibles non vides :
$$\bn^{<} =\bn_1^{<} \oplus \dots \oplus \bn_d^{<} ,$$
alors $B_{\bn^{<}}$ est un op\'erateur diff\'erentiel homog\`ene d'ordre $d$.\\

-- Si $\bn^{<} =\bn_1^{<} n_0$ alors
$$B_{\bn^{<}} =B_{\bn_1^{<}} B_{n_0} .$$
\end{defi}

Ces trois propri\'et\'es d\'efinissent bien $B_{\bn^{<}}$ qui se calcule par r\'ecurrence de la
mani\`ere suivante :
$$
\left .
\begin{array}{lll}
B_{\bn^{<}} & = & \di\sum_{i=1}^{\nu} (B_n (x_i )) \partial_{x_i}\ \mbox{\rm si}\ l(\bn^{<} ) =1,\\
B_{\bn^{<}} & = & \di\sum_{i=1}^{\nu} B_{\bn_1^{<}} (B_{n_0} (x_i )) \partial_{x_i} \ \mbox{\rm si}\
\bn^{<} =\bn_1^{<} n_0 .\\
B_{\bn^{<}} & = & \di {1\over d_1 ! \dots d_s !} \di\sum_{1\leq i_1 \dots i_d \leq \nu , 1\leq
j_1 \dots j_d \leq d}
\left ( B_{\bn_{j_1}^{<}} (x_{i_1} ) \right ) \dots \left ( B_{\bn_{j_d}^{<}} (x_{i_d} ) \right )
\partial_{x_{i_1}} \dots \partial_{x_{i_d}} ,
\end{array}
\right .
$$
o\`u $d_i$ est le nombre de suites arborescentes identiques $\bn_i^{<}$ qui interviennent dans la
d\'ecomposition de $\bn^{<} =\bn_1^{<} \oplus \dots \bn_d^{<}$.\\

Pour une suite arborescente $\ba^{<} =(a^1 ,\dots ,a^r )^{<}$ et une suite $\bn =(n^1 ,\dots ,n^{r'} )$,
on note $\mbox{\rm proj} \di\left (
\begin{array}{c}
\ba^{<} \\
\bn
\end{array}
\right ) $ le nombre de bijection de $\{ 1,\dots ,r \}$ dans $\{ 1,\dots ,r' \}$ (nul si $r\not= r'$)
v\'erifiant :
$$\mbox{\rm si}\ i_1 <i_2 \ \mbox{\rm dans}\ \ba^{<} ,\ \mbox{\rm alors}\ \sigma (i_1 )<\sigma (i_2 )
\ \mbox{\rm dans}\ \bn\ \mbox{\rm et}\ n^j =a^i \ \mbox{\rm si}\ j=\sigma (i).$$
On a alors les relations
$$S^{\ba^{<}} =\di\sum_{\bn \in A(X)^*}
\mbox{\rm proj} \di\left (
\begin{array}{c}
\ba^{<} \\
\bn
\end{array}
\right )
S^{\bn} ,
$$
pour tout moule $S^{\bullet}$, et
$$B_{\bn} =\di\sum_{\ba^{<} \in \mbox{\rm Arb} (X)} \mbox{\rm proj}
\di\left (
\begin{array}{c}
\ba^{<} \\
\bn
\end{array}
\right )
B_{\ba^{<}}
.$$

\subsection{D\'emonstration du th\'eor\`eme de Bruyno}

Par construction, on a
$$
\parallel B_{\bn^{<}} \parallel_{U,V} \leq C^r Q_1 ^{N(\bn^{<} )} ,\ Q_1 >0 ,
$$
pour une constante $C$ d\'ependant de $U$ et $V$, et
$$q(N) \leq Q_2^{N(\bn^{<} )} , Q_2 >0 .$$
La disparition du $r!$ dans la majoration des op\'erateurs arborifi\'e est
\'evidente. On devrait s'attendre \`a la retrouver dans la majoration des
moules arborifi\'es. Mais, et c'est l\`a que joue \`a fond les sym\'etries, on
a:

\begin{lem}
Pour toute suite arborescente $\bn^{<}$, on  a
\begin{equation}
\label{majmoularb}
\mid Na^{\bn^{<}} \mid \leq Q_3^{N(\bn^{<} )} ,\ Q_3 >0 .
\end{equation}
\end{lem}

La d\'emonstration de cette in\'egalit\'e n'est pas, contra\^\i rement au cas des comoules $B_{\bn^{<}}$, une cons\'equence
directe de la m\'ethode d'arborification. Elle r\'esulte du fait que l'\'equation diff\'erentielle satisfaite par le moule
$(\rm Na)^{-1}$ s'arborifie, i.e. que l'on a
\begin{equation}
\left [ \Delta ({\rm Na})^{-1} \right ]^{\bullet^{<}}=I^{\bullet^{<}} \times ({\rm Na}^{-1} )^{\bullet^{<}} .
\end{equation}
La forme du moule $({\rm Na}^{-1})^{\bullet^{<}}$ est donc la m\^eme que celle de ${\rm Na}^{-1}$ et les estimations
classiques sur les petits diviseurs permettent de conclure\footnote{Je ne ferai pas ces calculs ici, qui ne sont pas
simplifi\'es par la m\'ethode d'arborification ou l'utilisation du formalisme des moules. Jean Ecalle souligne qu'il faut utiliser les
estimations obtenues par Bruyno (\cite{br},p.207-224).}.\\

La convergence normale de la s\'erie arborifi\'ee s'en d\'eduit sans peine.

\begin{rema}
\`A ma connaissance, il n'y a pas de th\'eor\`eme g\'en\'eral concernant la m\'ethode d'arborification qui permet de
pr\'eciser son domaine d'application. Le fait que la m\'ethode restaure la convergence se fait au coup par coup sur les
exemples.
\end{rema}